\documentclass[opre,nonblindrev]{informs_TL}

\DoubleSpacedXI 

\usepackage{subcaption}
\usepackage[flushleft]{threeparttable}
\usepackage{booktabs,rotating}
\usepackage{appendix}
\usepackage{mathtools}
\usepackage{bbm}
\usepackage{mathrsfs}
\usepackage{yfonts}
\usepackage{amsmath}
\usepackage{graphicx}
\usepackage{amssymb}
\usepackage[utf8]{inputenc}
\usepackage[english]{babel}
\usepackage{etoolbox}
\usepackage{lscape} 
\usepackage{tcolorbox}
\usepackage{multicol} 
\usepackage{multirow} 
\usepackage{natbib}
\usepackage[ruled,linesnumbered]{algorithm2e}
\usepackage{wasysym} 
\usepackage{tikz}
\usepackage{lmodern}
 \tikzset{
    invisible/.style={opacity=0,text opacity=0},
    visible on/.style={alt={#1{}{invisible}}},
    alt/.code args={<#1>#2#3}{%
      \alt<#1>{\pgfkeysalso{#2}}{\pgfkeysalso{#3}} 
    },
  }
\usetikzlibrary{patterns}

\usepackage{pgfplots}
\pgfplotsset{compat=newest}
\usetikzlibrary{arrows,decorations.pathreplacing,backgrounds,automata}

\newcommand{\bA}{ \mathbf{A} }
\newcommand{\ba}{ \mathbf{a} }
\newcommand{\bbb}{ \mathbf{b} } 

\newcommand{\bc}{ \mathbf{c} }
\newcommand{\bD}{ \mathbf{D} }
\newcommand{\bd}{ \mathbf{d} } 
\newcommand{\beee}{ \mathbf{e} }

\newcommand{\bI}{ \mathbf{I} }

\newcommand{\bu}{ \mathbf{u} }

\newcommand{\bv}{ \mathbf{v} }

\newcommand{\bw}{ \mathbf{w} }

\newcommand{\bx}{ \mathbf{x} }

\newcommand{\by}{ \mathbf{y} }

\newcommand{\bz}{ \mathbf{z} }
\newcommand{\bbR}{ \mathbb{R} }
\newcommand{\cI}{ \mathcal{I} }

\newcommand{\cJ}{ \mathcal{J} }
\newcommand{\cX}{ \mathcal{X} }
\newcommand{\cK}{ \mathcal{K} }

\newcommand{\cS}{ \mathcal{S} }
\newcommand{\cT}{ \mathcal{T} }
\newcommand{\cR}{ \mathcal{R} }
\newcommand{\cN}{ \mathcal{N} }
\newcommand{\cC}{ \mathcal{C} }
\newcommand{\cA}{ \mathcal{A} }
\newcommand{\cG}{ \mathcal{G} }

\newcommand{\cH}{ \mathcal{H} }

\newcommand{\cV}{ \mathcal{V} }
\newcommand{\bxhat}{\mathbf{\hat{x}}}

\newcommand{\bzero}{ \mathbf{0} }
\newcommand{\balpha}{ \boldsymbol{\alpha} }

\newcommand{\bepsilon}{ \boldsymbol{\epsilon} }
\newcommand{\bvartheta}{ \boldsymbol{\vartheta} }

\newcommand{\blambda}{ \boldsymbol{\lambda} }
\newcommand{\bgamma}{ \boldsymbol{\gamma} }

\def\rev#1{{\color{black}#1}}

\usepackage{natbib}
 \bibpunct[, ]{(}{)}{,}{a}{}{,}%
 \def\bibfont{\small}%
\TheoremsNumberedThrough     
\ECRepeatTheorems

\EquationsNumberedThrough    

\begin{document}
\TITLE{Quantile Inverse Optimization: Improving Stability in Inverse Linear Programming}

\RUNAUTHOR{Shahmoradi and Lee}

\RUNTITLE{Quantile Inverse Optimization}

\ARTICLEAUTHORS{
\AUTHOR{Zahed Shahmoradi}
\AFF{Industrial Engineering, University of Houston, 4722 Calhoun Road, Houston, TX 77204, USA, \EMAIL{zshahmoradi@uh.edu}} 
\AUTHOR{Taewoo Lee}
\AFF{Industrial Engineering, University of Houston, 4722 Calhoun Road, Houston, TX 77204, USA, \EMAIL{tlee6@uh.edu}}
} 

\ABSTRACT{
Inverse linear programming (LP) has received increasing attention due to its potential to infer efficient optimization formulations that can closely replicate the behavior of a complex system. However, inversely inferred parameters and corresponding forward solutions from the existing inverse LP methods can be highly sensitive to noise, errors, and uncertainty in the input data, limiting their applicability in data-driven settings. We introduce the notion of inverse and forward stability in inverse LP and propose a novel inverse LP method that determines a set of objective functions that are stable under data imperfection and generate forward solutions close to the relevant subset of the data. We formulate the inverse model as a large-scale mixed-integer program (MIP) and elucidate its connection to biclique problems, which we exploit to \rev{develop efficient algorithms that solve much smaller MIPs instead to construct a solution to the original problem}. We numerically evaluate the stability of the proposed method and demonstrate its use in the diet recommendation and transshipment applications.
}


\maketitle

%

\section{Introduction}\label{sec:Intro}
Given a set of observed decisions as input data, inverse optimization infers parameters of a ``forward'' optimization problem, e.g., objective function coefficients, that make the given decisions optimal or near-optimal. By doing so, inverse optimization allows the forward problem to capture the preferences or utilities of the decision maker (DM) and reproduce the decisions accordingly. Inverse optimization has recently found numerous applications including finance \citep{bertsimas2012inverse}, transportation \citep{bertsimas2015}, electricity market \citep{saez2017,birge2017inverse}, incentive design \citep{aswani2019}, and healthcare \citep{erkin2010eliciting,lee2013predicting}. 

Various inverse optimization models have been developed for different types of forward problems including network optimization (e.g., \cite{ahuja2001inverse,heuberger2004inverse}), linear programs \citep{ahuja2001inverse,troutt2008linear,chan2014,ghate2015}, conic programs \citep{iyengar2005inverse}, convex programs \citep{keshavarz2011imputing}, integer and mixed-integer programs (MIPs) \citep{schaefer2009inverse,wang2009cutting,lamperski2015polyhedral}, and multicriteria optimization \citep{chan2018trade,Naghavi2019}. \rev{The underlying assumption of these studies is that the choice of the inverse model depends on which type of the forward problem is to be used for reproducing the decisions; e.g., if the user wants to derive a linear programming (LP) formulation that can replace a complex decision model, inverse LP is used to fit a linear objective function to the observed decisions.}

\rev{
As the goal of inverse optimization is to fit a model to given data, recent studies have naturally revealed and leveraged the connection between inverse optimization and regression. \cite{aswani2018inverse} highlight that the residuals used in their inverse model are similar to those in regression and propose a loss function analogous to the sum of squared errors in ordinary least squares. \cite{bertsimas2015} propose nonparametric inverse optimization motivated by kernel methods in nonparametric regression. \cite{chan2019goodness} study the similarity between inverse LP and linear regression and propose a goodness-of-fit measure for inverse LP, analogous to the R-squared measure in regression.
}

\rev{
Inverse LP has received a particular attention among different types of inverse problems due to its potential to find an LP formulation---whether the true decision-generating system is linear or not---which is easy to solve while closely replicating the system. While inverse optimization for more complex (e.g., non-convex) forward problems might lead to a better fit, the resulting forward problems can be computationally burdensome.} Inverse LP has also been used for inferring LP formulations customized to different DMs or problem instances, which can replace a complex formulation in a personalized or distributed manner, e.g., 
different LP cancer therapy planning formulations inferred for different types of patients, thus enabling personalized treatment modeling \citep{boutilier2015models}.

However, such potential of inverse LP often does not translate well into situations where the observed decisions are subject to noise, measurement errors, and uncertainty. In particular, objective functions (or cost vectors) obtained by the existing inverse LP methods can be very sensitive to small changes in the data. For example, suppose that the DM's true decisions turn out slightly different from the observed ones due to measurement errors, or some of them are simply outliers due to the DM's inconsistent behavior. The cost vector inferred by the existing methods often changes substantially in response to such errors and outliers, leading to substantially different forward solutions; also, the forward solutions are often far removed from the observed decisions, making it hard to predict the performance of the inferred cost vector (see Section~\ref{sec:prelim} for detailed illustration). Such instability can limit the applicability of inverse LP in data-driven settings.

In this study, we propose a novel inverse LP framework that addresses the instability issues caused by noise, errors, and uncertainty in data. We formally define stability measures in inverse optimization and propose a new inverse model that improves on the previous methods. 
Inspired by least quantile linear regression, our model infers a cost vector based on a quantile statistic of optimality errors associated with the observed decisions. Furthermore, we aim to find a set of cost vectors that are guaranteed to generate 
solutions within a specified distance from the relevant 
subset of the data. The presentation of our method in this paper focuses on the setting where the constraint matrix of the forward problem remains the same while observed decisions may vary. This setting can find various real world applications, e.g., diet problems where nutritional factors for each food do not change and network optimization problems where the geographical configuration (i.e., nodes and arcs) 
remains the same. 
We also discuss how this method can be extended when such assumptions do not hold.
\vspace{-0.1in}
\subsection{Relevant Literature}\label{sec:Related_Works}
Recent studies in inverse optimization focus on data-driven settings where a large (potentially noisy) dataset collected over a period of time or from many DMs is used for inferring the objective functions. 
\cite{keshavarz2011imputing} formulate an inverse model based on relaxed KKT conditions and impute a convex objective function that minimizes the KKT residuals with respect to the input solutions. \cite{bertsimas2015} consider inverse variational inequality with noisy data and find model parameters that minimize the optimality gap associated with the data. Similarly, \cite{esfahani2018data} develop a distributionally robust inverse optimization model to infer an objective function from noisy data
. \cite{aswani2018inverse} introduce the notion of risk consistency in inverse optimization 
and propose a model that finds an objective function that replicates the data in a statistically consistent manner. Inverse optimization has also been used for online learning where the inferred objective function is updated adaptively as new data are available over time \citep{barmann2017emulating,dong2018generalized}. Although these inverse convex programming frameworks can be specialized to inverse LPs, they are rather focused on generic convex programs with assumptions that preclude the above-mentioned instability issues (e.g., strictly convex feasible regions) and thus are not designed to address these particular issues in inverse LP.

\rev{
Although several recent inverse LP methods accommodate data that may not be optimal, these works focus on the development of closed-form solutions or efficient algorithms under often limiting assumptions \citep{chan2014, chan2018multiple, ghobadi2018robust}. As a result, how the obtained cost vectors actually work in the presence of data imperfection has received little attention. In particular, the previous studies exploit the polyhedral nature of the underlying forward LP and develop algorithms that find a cost vector orthogonal to one of the constraint vectors. While such algorithms are efficient, the resulting cost vector is highly sensitive to outliers or small data shifts.
}

\rev{
Our modeling approach to address the instability issues follows in spirit the line of work in the robust linear regression literature. 
Robust regression aims to infer a model that is stable against outliers or data shifts by increasing stability measure of the fitted model, known as the breakdown point, which is defined as the fraction of the data that can be altered arbitrarily without making the model arbitrarily bad. Various robust regression techniques have been proposed, including least median of squares \citep{rousseeuw1984least}, least trimmed squares \citep{hubert2008high,rousseeuw2006}, and least quantile of squares methods \citep{bertsimas2014least}. In these techniques, a parameter indicating the fraction of the data to be ``trimmed'' is pre-specified based on the user's application-specific knowledge and preferences, and model-fitting is done by minimizing a certain quantile error statistic, excluding the residuals associated with data points that are deemed to be outliers. 
}

\subsection{Contributions and Organization of the Paper} \label{sec:Contribution}
Our contributions are as follows.
\begin{enumerate}
\item We introduce the notion of forward and inverse stability in inverse optimization and show that the previous inverse LP methods are often unstable in the face of data noise, uncertainties, and outliers. We then propose a new inverse LP method that improves on the previous methods in terms of inverse and forward stability measures and provides bounds for both measures. 

\item We formulate the new inverse LP model as a large-size MIP and formally characterize the set of all feasible cost vectors. By establishing a new connection between the inverse model and a class of biclique problems, we develop efficient algorithms that solve many smaller MIPs (often in parallel) instead of directly solving the large MIP. We show that one of the algorithms is exact with a condition that is easy to check, while other heuristics are useful for very large instances. 
	
\item 
\rev{
We demonstrate two settings where the proposed inverse model is relevant: (i) the user observes decisions from an unknown system and aims to infer an LP formulation that can make new decisions similar to the observations
; (ii) the user provides an LP formulation to the DMs as a signal and collects their decisions as a response, which are then used for inferring their objective functions for the LP
. For the former, we demonstrate our method in the diet recommendation problem to quantify an individual's diet preferences from noisy and inconsistent data. For the latter, we use the transshipment problem to infer costs perceived by the DMs from a sequence of noisy datasets.
}

\end{enumerate}

The rest of the paper is organized as follows. In Section \ref{sec:prelim}, we illustrate how the previous inverse LP can be sensitive to data noise and outliers and formalize the notion of stability in inverse optimization. In Section \ref{sec:model}, we present the new inverse LP model that addresses the instability issues as well as its MIP reformulation. In Section~\ref{sec:model_algorithms}, we discuss the connection between the proposed model and biclique problems and propose efficient solution approaches. In Section \ref{sec:result}, we examine the performance of our model in terms of stability using various LP instances and demonstrate its use in the diet recommendation and transshipment problems. We conclude in Section~\ref{sec:conclusion}. Unless otherwise stated, proofs are in the appendix. 
\rev{
Throughout this paper, vectors are column vectors and ${\ba}'$ denotes the transpose of vector $\ba$, $\beee$ denotes the vector of ones, cone$(\cdot)$ denotes the set of conic combinations of given vectors, and cone$_+(\cdot)$ denotes the set of strict conic combinations of given vectors.
}

\section{Preliminaries}\label{sec:prelim}
\subsection{Forward Linear Program}
We consider the following forward optimization (FO) problem:
\begin{equation}\label{eq:FO}
\textrm{\textbf{FO}}(\bc):
\quad \underset{\bx}{\text{minimize}}\;\{\bc'\bx\,|\, \bA\bx \geq \bbb\},
\end{equation}
where $\bc \in \bbR^n, \bx \in \bbR^n$, $\bA \in \bbR^{m\times n}$, and $\bbb \in \bbR^m$. Let $\cI=\{1,\ldots,m\}$ index the constraints,  $\cJ=\{1,\ldots,n\}$ index the variables, and $\ba^i \in\mathbb{R}^n$ be a (column) vector corresponding to the $i$-th row of $\bA$. Let $\cX$ be the set of feasible solutions for the FO problem, assumed bounded,  full-dimensional and free of redundant constraints, and $\cX_i=\{\bx\in \cX \,|\, {\ba^i}'\bx=b_i\}$, $i\in\cI$. Without loss of generality, we assume that $\ba^i$ for each $i\in\cI$ is normalized {\it{a priori}} such that $\|\ba^i\|_p=1$ for some $p\ge 1$.

\subsection{Previous Inverse Linear Programming Method}
Consider a set of $K$ data points (or observations) $\hat\cX=\{\bxhat^1,\ldots,\bxhat^K\}$ with the index set $\cK=\{1,\ldots,K\}$. We make no assumption on the feasibility or optimality of the observations for the forward problem~\eqref{eq:FO}. Previous inverse LP methods with sub-optimal or infeasible observations aim to find a $\bc$ vector that can generate a forward optimal solution that is closest to the observations $\{\bxhat^k\}_{k\in\cK}$ under some distance metric, thus closely ``reproducing'' the data \citep{bertsimas2015, aswani2018inverse, chan2019goodness}. In LP, using general $\ell$-norm as a distance metric, such an inverse model can be written as follows \citep{aswani2018inverse,chan2018multiple}:
\begin{subequations}\label{eq:IO_prev}
	\begin{align}
	\quad \underset{\bc,\{\bepsilon^k\}_{k\in\cK},\by}{\text{minimize}}
	&\quad \sum_{k \in \cK}\lVert\bepsilon^k\rVert_{\ell} \label{eq:IO_prev.1}\\
	 \text{subject to}
     	& \quad \bA'\by=\bc, \label{eq:IO_prev.2}\\
	& \quad \by \geq \bzero, \label{eq:IO_prev.3}\\
	& \quad \bA(\bxhat^k-\bepsilon^k)\geq \bbb, \quad \forall k \in \cK, \label{eq:IO_prev.4}\\
   	& \quad \bc'(\bxhat^k-\bepsilon^k)=\bbb'\by, \quad \forall k \in \cK, \label{eq:IO_prev.5}\\
	& \quad \lVert\bc\rVert_{p}=1, \label{eq:IO_prev.6}
	\end{align}
\end{subequations}
where $\ell\ge1$, $p\ge1$, and $\bepsilon^k$ denotes a perturbation vector for observation $\hat\bx^k$. Given observations $\{\hat\bx^k\}_{k\in\cK}$, the above problem finds a $\bc$ vector that makes each perturbed solution $\hat\bx^k - \bepsilon^k$ satisfy dual feasibility \eqref{eq:IO_prev.2}--\eqref{eq:IO_prev.3}, primal feasibility \eqref{eq:IO_prev.4}, and strong duality \eqref{eq:IO_prev.5}, thus rendering it optimal, while the perturbations (i.e., ``optimality errors'' associated with the observations) are minimized in $\ell$-norm in the objective function \eqref{eq:IO_prev.1}. The normalization constraint \eqref{eq:IO_prev.6} prevents a trivial, all-zero cost vector from being feasible.

Note that existing inverse convex programming models (e.g., \cite{aswani2018inverse}) can be written equivalently as \eqref{eq:IO_prev} when the underlying forward problem is exactly \eqref{eq:FO}. While the above problem is non-convex, recent works propose an efficient, exact algorithm for the problem by exploiting the solution structure \citep{chan2018multiple, chan2019goodness}, namely that an optimal cost vector $\bc^*$ is orthogonal to one of the hyperplanes defining polyhedron $\cX$ (i.e., $\bc^*=\ba^i$ for some $i\in\cI$).

\subsection{Instability Issues in Inverse LP} \label{secsub:Instability_Issues}
\begin{figure}[t]\centering
\begin{subfigure}[b]{.25\textwidth}
\centering
\vfill
\begin{tikzpicture}[scale=1.3]
   \draw[->] (0,0) -- coordinate (x axis mid) (3,0) node[pos=1,below=0.1cm] {$x_1$};
    \foreach \x in {0,0.5,...,2.5}
		\draw  (\x,1pt)--(\x,-1pt);
    \foreach \x in {0,...,2.5}
		\draw  (\x,1pt)--(\x,-1pt)
		node[anchor=north] {\x};
    \draw[->] (0,0) -- coordinate (y axis mid) (0,3) node[pos=1,left] {$x_2$};
    \foreach \y in {0,0.5,...,2.5}
		\draw  (1pt,\y)--(-1pt,\y);
    \foreach \y in {0,...,2.5}
		\draw  (1pt,\y)--(-1pt,\y)
		node[anchor=east] {\y};
		

	\draw (2.5,0) -- (2.5,2.5);
	\draw[gray,dotted] (2.5,0)--(2.5,-0.5);
	\draw[gray,dotted] (2.5,2.5)--(2.5,3);

	\draw (2.5,2.5) -- (0,2.5);
	\draw[gray,dotted] (-0.5,2.5) -- (0,2.5);
	\draw[gray,dotted] (2.5,2.5) -- (3,2.5);

	
	\draw (0.75,3) -- (1.75,3); 
	\draw[->] (1.25,3) coordinate -- (1.25,2.7) 
	     node[anchor=west] {$\ba^1$};
    
    \draw (3,0.75) -- (3,1.75); 
	\draw[->] (3,1.25) coordinate -- (2.75,1.25)  node[anchor=south] {$\ba^2$};
    
    \draw[fill=black] (2,2.3) circle (0.04);
    \node [black,anchor= east] at (2.1,2.4) {\scriptsize$\hat\bx^1$};
    
    \draw[fill=black] (2.2,2.3) circle (0.04);
    \node [black,anchor= west] at (2.15,2.4) {\scriptsize$\hat\bx^2$};
    
    \draw[fill=black] (2.2,2) circle (0.04);
    \node [black,anchor= west] at (2.15,2.07) {\scriptsize$\hat\bx^3$};
    
    \draw[fill=black] (2,2) circle (0.04);
    \node [black,anchor= east] at (2.1,2.1) {\scriptsize$\hat\bx^4$};
    
    \draw[gray,fill=gray] (2.3,2.2) circle (0.04);
   
    \draw[gray,fill=gray] (2.3,1.9) circle (0.04);

\end{tikzpicture}
\caption{Inverse instability due to data shift.}
\label{fig:instab_inv_shift}
\end{subfigure}
\hspace{0.4in}
\begin{subfigure}[b]{.25\textwidth}
\centering
\vfill
\begin{tikzpicture}[scale=1.3]
\draw[->] (0,0) -- coordinate (x axis mid) (3,0) node[pos=1,below=0.1cm] {$x_1$};
    \foreach \x in {0,0.5,...,2.5}
		\draw  (\x,1pt)--(\x,-1pt);
    \foreach \x in {0,...,2.5}
		\draw  (\x,1pt)--(\x,-1pt)
		node[anchor=north] {\x};
    \draw[->] (0,0) -- coordinate (y axis mid) (0,3) node[pos=1,left] {$x_2$};
    \foreach \y in {0,0.5,...,2.5}
		\draw  (1pt,\y)--(-1pt,\y);
    \foreach \y in {0,...,2.5}
		\draw  (1pt,\y)--(-1pt,\y)
		node[anchor=east] {\y};
		

	\draw (2.5,0) -- (2.5,2.5);
	\draw[gray,dotted] (2.5,0)--(2.5,-0.5);
	\draw[gray,dotted] (2.5,2.5)--(2.5,3);

	\draw (2.5,2.5) -- (0,2.5);
	\draw[gray,dotted] (-0.5,2.5) -- (0,2.5);
	\draw[gray,dotted] (2.5,2.5) -- (3,2.5);

	
	\draw (0.75,3) -- (1.75,3); 
	\draw[->] (1.25,3) coordinate -- (1.25,2.7) 
	     node[anchor=west] {$\ba^1$};
    
    \draw (3,0.75) -- (3,1.75); 
	\draw[->] (3,1.25) coordinate -- (2.75,1.25)  node[anchor=south] {$\ba^2$};
	
    \draw[fill=black] (2,2.3) circle (0.04);
    \node [black,anchor= east] at (2.1,2.4) {\scriptsize$\hat\bx^1$};
    
    \draw[fill=black] (2.2,2.3) circle (0.04);
    \node [black,anchor= west] at (2.15,2.4) {\scriptsize$\hat\bx^2$};
    
    \draw[fill=black] (2.2,2) circle (0.04);
    \node [black,anchor= west] at (2.15,2.1) {\scriptsize$\hat\bx^3$};
    
    \draw[fill=black] (2,2) circle (0.04);
    \node [black,anchor= east] at (2.1,2.1) {\scriptsize$\hat\bx^4$};
    
    \draw[red,fill=red] (2.2,0.3) circle (0.04);
    \node [red,anchor= east] at (2.2, 0.3) {\scriptsize$\hat\bx^5$};
\end{tikzpicture}
\vfill
\caption{Inverse instability due to outliers.}
\label{fig:instab_inv_outlier}
\end{subfigure}
\hspace{0.4in}
\begin{subfigure}[b]{.25\textwidth}
\centering
\vfill
\begin{tikzpicture}[scale=1.3]
    \draw[->] (0,0) -- coordinate (x axis mid) (3,0) node[pos=1,below=0.1cm] {$x_1$};
    \foreach \x in {0,0.5,...,2.5}
		\draw  (\x,1pt)--(\x,-1pt);
    \foreach \x in {0,...,2.5}
		\draw  (\x,1pt)--(\x,-1pt)
		node[anchor=north] {\x};
    \draw[->] (0,0) -- coordinate (y axis mid) (0,3) node[pos=1,left] {$x_2$};
    \foreach \y in {0,0.5,...,2.5}
		\draw  (1pt,\y)--(-1pt,\y);
    \foreach \y in {0,...,2.5}
		\draw  (1pt,\y)--(-1pt,\y)
		node[anchor=east] {\y};
		
	\draw (2.5,0) -- (2.5,2.5);
	\draw[gray,dotted] (2.5,0)--(2.5,-0.5);
	\draw[gray,dotted] (2.5,2.5)--(2.5,3);

	\draw (2.5,2.5) -- (0,2.5);
	\draw[gray,dotted] (-0.5,2.5) -- (0,2.5);
	\draw[gray,dotted] (2.5,2.5) -- (3,2.5);


    
    \draw[fill=black] (2,2.3) circle (0.04);
    \node [black,anchor= east] at (2.1,2.4) {\scriptsize$\hat\bx^1$};
    
    \draw[fill=black] (2.2,2.3) circle (0.04);
    \node [black,anchor= west] at (2.15,2.4) {\scriptsize$\hat\bx^2$};
    
    \draw[fill=black] (2.2,2) circle (0.04);
    \node [black,anchor= west] at (2.15,2.1) {\scriptsize$\hat\bx^3$};
    
    \draw[fill=black] (2,2) circle (0.04);
    \node [black,anchor= east] at (2.1,2.1) {\scriptsize$\hat\bx^4$};
    
    \draw[fill=black] (2.2,0.3) circle (0.04);
    \node [black,anchor= east] at (2.2, 0.3) {\scriptsize$\hat\bx^5$};
    

    \draw (3,0.75) -- (3,1.75); 
	\draw[->] (3,1.25) coordinate -- (2.75,1.25)  node[anchor=south] {$\ba^2$};
	
    \draw[red,fill=red] (2.5,0) circle (0.04);
  	\node[red, anchor=south west] at (2.5,0) (coord) {\scriptsize$\bx^{*}$};
\end{tikzpicture}
\vfill
\caption{Forward instability.
\newline}
\label{fig:instab_FW}
\end{subfigure}
\caption{
\rev{
Instability issues of the previous inverse LP method. 
}
}
\label{fig:instability}
\end{figure}
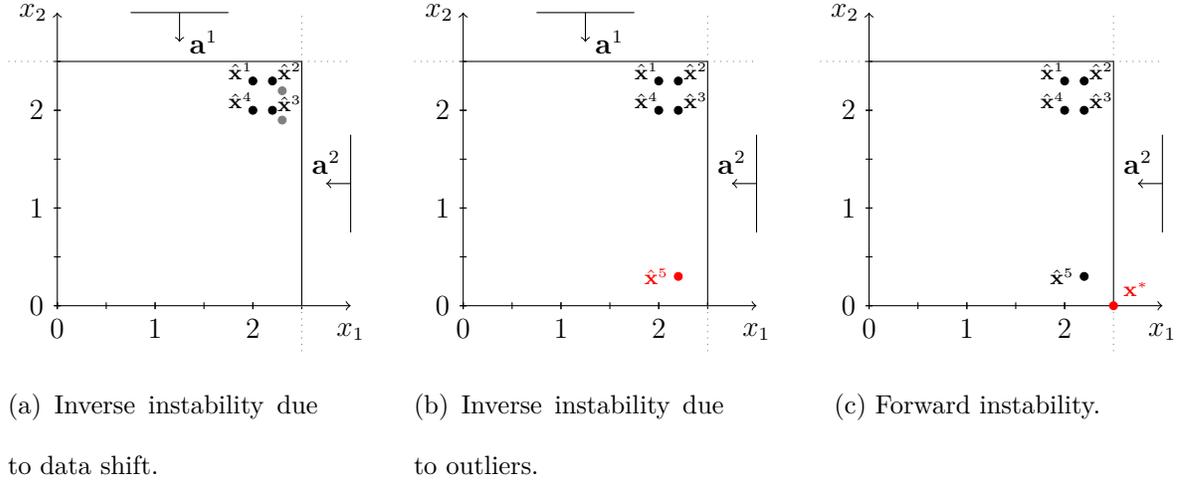

When the dataset $\hat\cX$ contains noise or measurement errors, the previous inverse model~\eqref{eq:IO_prev} can be unstable in terms of both the cost vector $\bc^*$ it produces as well as the forward solution that $\bc^*$ generates (i.e., $\bx^*\in\argmin \textbf{FO}(\bc^*)$). 
\rev{The following examples illustrate such instability issues.} 
\rev{
\begin{example}\label{Example:Inv_Instability_shift}
Consider the forward LP: $\underset{\bx}{\text{minimize}}\;\{\bc'\bx\,|\, \bA\bx \geq \bbb\}$
where $\ba^1=\begin{bsmallmatrix}0\\-1\end{bsmallmatrix}, \; \ba^2=\begin{bsmallmatrix}-1\\0\end{bsmallmatrix}, \; \ba^3=\begin{bsmallmatrix}0\\1\end{bsmallmatrix},$ and $\ba^4=\begin{bsmallmatrix}1\\0\end{bsmallmatrix}$ are the rows of $\bA$ (written as column vectors) and $\bbb=[-2.5,-2.5,0,0]'$, with four data points $\hat\bx^1=\begin{bsmallmatrix}2\\2.3\end{bsmallmatrix}$, $\hat\bx^2=\begin{bsmallmatrix}2.2\\2.3\end{bsmallmatrix}$, $\hat\bx^3=\begin{bsmallmatrix}2.2\\2\end{bsmallmatrix}$, and $\hat\bx^4=\begin{bsmallmatrix}2\\2\end{bsmallmatrix}$ (see Figure \ref{fig:instab_inv_shift}). The previous inverse model \eqref{eq:IO_prev} with  $\lVert\cdot\rVert_\ell=\lVert\cdot\rVert_\infty$ finds $\bc^*=\ba^1=\begin{bsmallmatrix}0\\-1\end{bsmallmatrix}$ as the unique optimal cost vector with the objective value of $1.4$ ($\bepsilon^1=\begin{bsmallmatrix}0\\-0.2\end{bsmallmatrix}$, $\bepsilon^2=\begin{bsmallmatrix}0\\-0.2\end{bsmallmatrix}$, $\bepsilon^3=\begin{bsmallmatrix}0\\-0.5\end{bsmallmatrix}$, and $\bepsilon^4=\begin{bsmallmatrix}0\\-0.5\end{bsmallmatrix}$). Suppose that $\hat\bx^2$ and $\hat\bx^3$ are shifted to $\begin{bsmallmatrix}2.3\\2.2\end{bsmallmatrix}$ and $\begin{bsmallmatrix}2.3\\1.9\end{bsmallmatrix}$, respectively (shown in gray). Even with these small shifts, the optimal cost vector changes from $\bc^*=\begin{bsmallmatrix}0\\-1\end{bsmallmatrix}$ (which now leads to the objective value of 1.6) to a substantially different vector $\bc^*=\ba^2=\begin{bsmallmatrix}-1\\0\end{bsmallmatrix}$ with the objective value of 1.4 ($\bepsilon^1=\begin{bsmallmatrix}-0.5\\0\end{bsmallmatrix}$, $\bepsilon^2=\begin{bsmallmatrix}-0.2\\0\end{bsmallmatrix}$, $\bepsilon^3=\begin{bsmallmatrix}-0.2\\0\end{bsmallmatrix}$, and $\bepsilon^4=\begin{bsmallmatrix}-0.5\\0\end{bsmallmatrix}$).
\end{example}
\begin{example}\label{Example:Inv_Instability_outlier}
Consider the same initial data points from Example~\ref{Example:Inv_Instability_shift}. Suppose that an outlier $\hat\bx^5=\begin{bsmallmatrix}2.2\\0.3\end{bsmallmatrix}$ is introduced (see Figure \ref{fig:instab_inv_outlier}). The optimal cost vector then changes from $\bc^*=\ba^1=\begin{bsmallmatrix}0\\-1\end{bsmallmatrix}$ to a substantially different vector $\bc^*=\ba^2=\begin{bsmallmatrix}-1\\0\end{bsmallmatrix}$ with the objective value of $1.9$ ($\bepsilon^1=\begin{bsmallmatrix}-0.5\\0\end{bsmallmatrix}$, $\bepsilon^2=\begin{bsmallmatrix}-0.3\\0\end{bsmallmatrix}$, $\bepsilon^3=\begin{bsmallmatrix}-0.3\\0\end{bsmallmatrix}$, $\bepsilon^4=\begin{bsmallmatrix}-0.5\\0\end{bsmallmatrix}$, and $\bepsilon^5=\begin{bsmallmatrix}-0.3\\0\end{bsmallmatrix}$). The previous solution now has the objective value of $1.4\,+\,\rVert\begin{bsmallmatrix}0\\-2.2\end{bsmallmatrix}\lVert_\infty=3.6$, hence no longer optimal. 
\end{example}
\begin{example}\label{Example:FW_Instability}
Consider the five data points from Example~\ref{Example:Inv_Instability_outlier} where the optimal cost vector is $\bc^*=\ba^2=\begin{bsmallmatrix}-1\\0\end{bsmallmatrix}$. While a desirable forward solution is supposed to be close to the relevant subset of the data without the outlier, Figure~\ref{fig:instab_FW} shows that using this cost vector for the forward LP can lead to a solution $\bx^*$ that is far from the majority of the dataset.
\end{example}
}

Examples \ref{Example:Inv_Instability_shift} and \ref{Example:Inv_Instability_outlier} show that the previous inverse model is sensitive to outliers or small data shifts, which we refer to as being \textit{inverse-unstable}. On the other hand, Example \ref{Example:FW_Instability} shows that solving the forward problem with a cost vector 
from the previous model can lead to a solution that is (unexpectedly) far from the relevant subset of the data, which we refer to as being \textit{forward-unstable}. 

\rev{
We note that these instability issues are common in general inverse LP settings. For example, if a data-generating system produces errors 10\% of the time, the inferred cost vector can still be heavily dragged toward these errors spread far from the majority of the data, regardless of the size of the dataset, leading to a cost vector that would have not been chosen if only the pristine data (the remaining 90\%) was used. Also, when small shifts occur in multiple data points 
according to the same distribution (e.g., the DM's behavior changing gradually over time), the cost vector can also suddenly change 
to a substantially different vector, leading to a substantially different forward solution. In general, these instability issues depend heavily on the spatial distribution of the data in relation to 
polyhedron $\cX$, e.g., the distance between each data point and each facet of $\cX$.
}

\subsection{Stability Measures for Inverse LP}\label{secsub:Instability_Measures}
In this subsection, we formally define a notion of stability in inverse LP and propose measures that we use to assess the stability of an inverse LP model. Given a dataset $\hat\cX$ and a certain inverse LP model, let $\hat\cC$ be the set of cost vectors obtained by the model.

\textbullet \textbf{ Inverse stability:} We propose to measure inverse stability of an inverse LP model (or solutions to the model) by the minimum tolerable data shift in $\hat\cX$ until the model loses all of its initial solutions $\hat \cC$. That is, if we let $\tilde\cX$ denote the shifted data and $\tilde\cC$ denote the set of cost vectors obtained by the model with $\tilde\cX$, inverse stability is measured by  
\begin{align} \label{eq:inv_stab}
\underset{\tilde\cX}{\min} \{d(\hat\cX,\tilde\cX) \, | \, \hat\cC \cap \tilde\cC = \emptyset\},
\end{align}
where $d$ is some distance function; e.g., $d(\hat\cX,\tilde\cX)= \max\{\|\hat\bx - \tilde\bx\|_{\ell}\,|\,\hat\bx\in\hat\cX,\tilde\bx\in\tilde\cX\}$ or $d(\hat\cX,\tilde\cX)=\sum_{k\in\cK} \|\hat\bx^k - \tilde\bx^k\|_{\ell},\ell \ge1,$ if $\hat\cX$ and $\tilde \cX$ are in one-to-one correspondence. That is, ideally, a stable inverse model should maintain some cost vectors when reasonably small changes occur in the data, as quantified by the above measure. \rev{Note that this measure is analogous to the stability measure in regression, i.e., the finite-sample breakdown point, defined as the fraction of the data that can be altered without spoiling the inferred model completely \citep{yohai1987high}.}

\textbullet \textbf{ Forward stability:} 
Suppose we select $\hat\bc\in\hat\cC$ and find a set of forward optimal solutions $\cX^*(\hat\bc)=\argmin\, \textrm{\textbf{FO}}(\hat\bc)$. How \textit{unstable} this cost vector $\hat\bc$ can be is assessed by how far a forward solution $\bx\in\cX^*(\hat\bc)$ can be from the given observations $\hat\cX$%
, i.e., the worst-case distance between $\cX^*(\hat\bc)$ and $\hat\cX$:
\begin{align}\label{eq:for_stab}
\, \underset{\bx \in \cX^*\!(\hat\bc)}\max  \{d(\hat\cX,\bx)\},
\end{align}
where the distance function $d$ can be defined similarly as above. That is, to guarantee forward stability in inverse LP, an inverse model should identify a cost vector that is guaranteed to produce an optimal solution that is reasonably close to the observations.

\rev{
The inverse and forward stability measures often conflict with each other. For example, a cost vector whose elements are all zeros is perfectly inverse-stable, yet using this vector for the forward problem likely leads to a solution far from the given data (i.e., highly forward-unstable). Also, one might attempt to include as many cost vectors as possible in $\hat\cC$ just to make it less vulnerable to data shifts (i.e., more inverse-stable), yet which of these vectors is truly forward-stable then becomes less clear. Another challenge is that neither of the above measures lends itself to a tractable inverse optimization problem. In the next section, we provide a tractable inverse LP framework that improves both inverse and forward stability as quantified by the above measures.
}

\section{Models}\label{sec:model}
In this section, we first propose a new inverse LP model that improves the inverse stability of the previous model \eqref{eq:IO_prev}. We then analyze the solution structure of the model, which we exploit to reformulate the model as an MIP and characterize the set of all feasible cost vectors. Finally, we introduce a specific objective function for the MIP which can serve as a proxy for the forward stability measure and leads to a cost vector with improved forward stability.


\vspace{-0.1in}
\subsection{Inverse LP for Quantile Statistics: Improving Inverse Stability}\label{sec:model1}
As discussed, inverse instability is caused by the vulnerability of inferred cost vectors to data shifts or outliers. Our modeling strategy to address this is inspired by the close relationship between the least squares method in  
regression and inverse optimization: the previous inverse model \eqref{eq:IO_prev} is similar to the least squares method in that it fits the LP model by minimizing the sum 
of optimality errors \citep{chan2019goodness}. In regression, one way to address data imperfection 
is to minimize the $\theta$ quantile statistic where $\theta\in[0,1]$ (i.e., ($\theta\times100$)-th percentile) instead of the sum of the squared errors, which is known as the least quantile of squares method \citep{Koenker2001,bertsimas2014least} or least trimmed quantile regression \citep{rousseeuw1984least, rousseeuw2006}.
  
We adopt this idea and propose a generalized inverse framework where a cost vector is inferred based on the $\theta$ quantile optimality error, i.e., 
associated with a data point that induces the $\lceil \theta K \rceil$-th smallest optimality error with respect to the cost vector. Instead of minimizing this $\theta$ quantile error as was done in \cite{bertsimas2014least}, we find a set of cost vectors such that 
the $\lceil \theta K \rceil$-th smallest optimality error is no greater than a certain threshold $\tau$
. That is, such cost vectors render a relevant subset of the data (of cardinality $\lceil\theta K\rceil$) within the threshold optimality error. This problem, which we call the quantile inverse optimization (QIO) problem, can be written as follows:
\vspace{-0.1in}
\begin{subequations}\label{eq:IO1}
	\begin{align}
	\textrm{\textbf{QIO}}(\cK,\tau,\theta):
	\quad \underset{\bc,\{\bepsilon^k\}_{k\in\cK},\by,\cS}{\text{minimize}} \quad & 0 \label{eq:IO1.1}\\
	\text{subject to} \quad & \bA'\by=\bc, \quad \label{eq:IO1.2}\\
	&  \by \geq \bzero, \label{eq:IO1.3}\\
	&  \bA(\bxhat^k-\bepsilon^k)\geq \bbb, \quad \forall k \in \cS, \label{eq:IO1.4}\\
	&  \bc'(\bxhat^k-\bepsilon^k)=\bbb'\by,  \quad \forall k \in \cS, \label{eq:IO1.5}\\
	&  \|\bepsilon^k\|_{\ell} \leq \tau, \quad \forall k \in \cS, \label{eq:IO1.6}\\
	&  |\cS| \ge \theta K, \label{eq:IO1.7}\\
	&  \cS \subseteq \cK, \label{eq:IO1.8}\\
	&  \|\bc\|_{p}=1. \label{eq:IO1.9}
	\end{align}
\end{subequations}
%
%
Constraints \eqref{eq:IO1.2}--\eqref{eq:IO1.3} represent dual feasibility and \eqref{eq:IO1.4}--\eqref{eq:IO1.5} enforce primal feasibility and strong duality associated with the perturbed solutions $\hat\bx^k - \bepsilon^k$ for a subset of the data $\cS\subseteq\cK$ whose cardinality is enforced to be 
no less than $ \theta K$ by 
\eqref{eq:IO1.7}. Constraint \eqref{eq:IO1.6} then ensures that optimality error for each chosen observation 
is within $\tau$. Clearly, these constraints ensure that the $\lceil \theta K \rceil$-th smallest optimality error is no greater than $\tau$. Note that this problem is written as a feasibility problem---we later introduce an objective function that leads to a specific subset of the feasible cost vectors (see Section~\ref{sec:model2}). We call cost vectors that are feasible for 
\eqref{eq:IO1.2}--\eqref{eq:IO1.9}  ``inverse-feasible'' cost vectors. \rev{Similar to quantile-based regression, the choice of $\theta$ depends on the application context or user preference; e.g., $\theta$ can be chosen based on the user's belief about the fraction of outliers and can also be adjusted post-hoc depending on the result from the model. Or, it can be set to 50\% if one is interested in the median error-based inverse problem.} The above problem will be referred to as $\textbf{QIO}(\cK,\tau,\theta)$ or $\textbf{QIO}(\hat\cX,\tau,\theta)$ interchangeably, depending on which is more convenient for the context.

\rev{
A strong advantage of this modeling framework is that it generalizes the previous model \eqref{eq:IO_prev}---given the same dataset $\hat\cX$, there exists $\tau$ such that $\textbf{QIO}(\hat\cX,\tau,\theta)$ with $\theta = 100\%$ produces the same set of feasible cost vectors as \eqref{eq:IO_prev}. Furthermore, given such $\tau$, as $\theta$ decreases (i.e., the feasible region expands), the set of feasible cost vectors from $\textbf{QIO}(\hat\cX,\tau,\theta)$ contains the cost vectors from the previous model. That is, in the presence of data shift or additional outliers, $\textbf{QIO}(\hat\cX,\tau,\theta)$ does not lose all of its initial inverse-feasible solutions until the previous model does; i.e., it is at least as stable as the previous model as quantified by inverse stability measure~\eqref{eq:inv_stab}. We formalize this idea as follows. 
\begin{proposition}\label{prop:ISLB_Nondecreasing}
Let $\hat\cC(\theta)$ be the set of inverse-feasible cost vectors for $\mathbf{QIO}(\hat\cX,\tau,\theta)$. Given $\hat\cX$ and $\tau$, inverse stability measure \eqref{eq:inv_stab} for $\hat\cC(\theta)$ is non-decreasing as $\theta$ decreases. 
\end{proposition}
We later further show that a solution to $\textbf{QIO}(\hat\cX,\tau,\theta)$ provides a lower bound on the inverse stability measure. In the following, we first analyze the solution structure of $\textbf{QIO}(\hat\cX,\tau,\theta)$.
}

\noindent\subsubsection{Solution Structure and MIP Reformulation.}\label{secsub:QIO_sol_struct}
The QIO problem is non-convex due to constraint \eqref{eq:IO1.5} and thus is hard to solve. We analyze the solution structure of this problem, which leads to an MIP reformulation 
and a lower bound on the inverse stability measure.
\begin{proposition}\label{prop:IO1}
If $\mathbf{QIO}(\cK,\tau,\theta)$ is feasible, then there exists a feasible ${\bc}$ for $\mathbf{QIO}(\cK,\tau,\theta)$ such that ${\bc}=\ba^i$ for some $i \in \cI$.
\end{proposition}
Proposition \ref{prop:IO1} implies that the feasibility of the QIO model can be checked by evaluating at most $m$ constraint vectors, i.e., $\ba^i, i=1,\ldots,m$. 
\rev{
For example, revisiting Example~\ref{Example:Inv_Instability_shift} (Figure \ref{fig:instab_inv_shift}), for the QIO problem with $\|\cdot\|_\ell=\|\cdot\|_\infty$, $\theta=0.5$, and $\tau=0.3$, both $\bc=\ba^1$ and $\bc=\ba^2$ are feasible: $\bc\!=\!\ba^1$ leads to $\|\bepsilon^1\|_\infty\!=\!\|\bepsilon^2\|_\infty\!=\!0.2$ and $\bc\!=\!\ba^2$ leads to $\|\bepsilon^2\|_\infty\!=\!\|\bepsilon^3\|_\infty\!=\!0.3$, both satisfying the threshold $\tau=0.3$. 
}
Proposition \ref{prop:IO1}  also suggests that if there is no $\bc=\ba^i$ feasible for \eqref{eq:IO1} for any $i \in \cI$, then there exists no cost vector that can make at least $ \lceil \theta K \rceil$ observations within $\tau$-optimality. In this case, users could decrease $\theta$ or increase $\tau$ to make the model feasible. Or, $\tau$ can be set by adding a reasonable margin to the minimum possible value of $\tau$ that keeps the model feasible \rev{(see Appendix \ref{appx:finding_E} for more details)}. With an appropriate threshold $\tau$, the set of inverse-feasible cost vectors that are orthogonal to some hyperplanes defining $\cX$ can be found efficiently by evaluating each hyperplane $i$ (see Algorithm~\ref{alg:IO1}).

{\OneAndAHalfSpacedXI
\begin{algorithm}[H]\label{alg:IO1}
	\caption{Finding inverse-feasible solutions of the form $\bc=\ba^i$, $i\in\cI$}
	\SetAlgoLined
	\KwResult{
	$\cC = \{\ba^i\}_{i\in\cA}$
	}
	\KwIn{$\bA, \bbb, 
	\hat\cX, \tau,\theta$}
	$ \cS_i 
	\leftarrow \emptyset \ \forall i \in \cI$; $\cA
	\leftarrow\emptyset$\\
	\For{$i \in \cI$}{
		\For{$k \in \cK$}{
			\If{$\exists\ \bx \in \cX_i$ such that $\lVert \bxhat^k-\bx \rVert_{\ell} \leq \tau$}{$\cS_i	\leftarrow\cS_i \cup \{k \}$;}			
		}    
			\If{$|\cS_i|\geq \theta K$ }{$\cA
			\leftarrow\cA \cup \{i\}$} 
	}
\end{algorithm}
}

Note that Algorithm \ref{alg:IO1} returns the set of $\ba^i$'s that are inverse-feasible, i.e., $\cA$, which is only a subset of the set of all inverse-feasible solutions. 
\rev{
We extend the results in Proposition \ref{prop:IO1} and Algorithm \ref{alg:IO1} and derive the following result, which further identifies other inverse-feasible cost vectors.
\begin{lemma}\label{lemma:Feasible_conic}
If there exists $\bar\bx^k \in \cap_{i\in\bar\cA}\cX_i$ for some set $\bar\cA\subseteq\!\cA$ such that $\|\bar\bx^k-\hat\bx^k\|_\ell\!\le \tau$ for all $k\in {\bar\cS}\subseteq\cK$ where $|{\bar\cS}|\ge \theta K$, then $\bc\in \textup{cone}(\{\ba^i\}_{i\in \bar\cA})$ such that $\|\bc\|_p=1$ is inverse-feasible.
\end{lemma}
}

\rev{
Lemma~\ref{lemma:Feasible_conic} suggests that we can check each subset $\bar\cA$ of $\cA$ to see if a conic combination of $\ba^i$'s for $i\in\bar\cA$ is also inverse-feasible. Using this observation, we propose the following MIP:
}
\vspace{-0.15in}
\begin{subequations}\label{eq:IO_mip}	
	\begin{align}	
	\underset{\bv,\bu,\{\bepsilon^k\}_{k\in\cK}}{\text{minimize}} 
	&\quad 0 \label{eq:IO_mip.1}\\	
	 \text{subject to} 
	&\quad b_i \leq {\ba^i}'(\bxhat^k-\bepsilon^k) \leq b_i +M_1(1-v_i), \quad \forall i  \in \cI,\forall k \in  \cK,  \label{eq:IO_mip.2}\\
	&\quad \lVert\bepsilon^k\rVert_{\ell}\leq \tau+M_2(1-u_k), \quad \forall k \in \cK,  \label{eq:IO_mip.3}\\
	&\quad \sum_{k=1}^{ K} u_k \ge \theta K, \label{eq:IO_mip.4} \\
	&\quad  v_i,u_k \in \{0,1\}, 
	\quad \forall i  \in \cI ,\forall k \in \cK, \label{eq:IO_mip.5}
	\end{align}	
\end{subequations}
%
where $v_i=1$ if ${\ba^i}'(\bxhat^k-\bepsilon^k)\ge b_i$ holds with equality (i.e., $\bxhat^k-\bepsilon^k\in\cX_i$) for all $k \in \cK$ and $v_i=0$ otherwise, 
and $u_k=1$ if observation $\hat\bx^k$ is ``chosen'' and $u_k=0$ otherwise. Parameters $M_1$ and $M_2$ denote sufficiently large constants; Appendix \ref{appx:Big-M} shows how to find appropriate values for them.   
\rev{
In formulation \eqref{eq:IO_mip}, inverse-feasible $\ba^i$'s are identified by letting $v_i = 1$
; the hypothesis 
of Lemma \ref{lemma:Feasible_conic} is then explicitly written as constraints 
to further identify a conic combination of the chosen $\ba^i$'s that is also inverse-feasible (the set $\cI$ in 
\eqref{eq:IO_mip.2} and \eqref{eq:IO_mip.5} can be replaced by $\cA$ without losing any feasible solution because $v_i$ can be 1 only for $i\in\cA\subseteq\cI$ anyways; we use $\cI$ in \eqref{eq:IO_mip} and subsequent formulations 
for notational clarity).
} 
The following result establishes the equivalence between problem \eqref{eq:IO_mip} and \textbf{QIO}$(\cK,\tau,\theta)$ in terms of the set of achievable inverse-feasible cost vectors given the same dataset $\hat\cX$.
\begin{theorem} \label{prop:IO_mip}
	A cost vector $\bc$ is feasible for $\mathbf{QIO}(\cK,\tau,\theta)$ if and only if there exists a feasible solution $(\bar\bv$ $\ne$ $\mathbf{0},\bar\bu,\{\bar{\bepsilon}^k\}_{k\in \cK})$ for model \eqref{eq:IO_mip} such that $\bc \in \textup{cone}(\{\ba^i\}_{ i:  \bar{v}_i=1})$ and $\lVert \bc \rVert_p=1$.
\end{theorem}
%
\rev{Thus, an inverse-feasible cost vector for \textbf{QIO}$(\cK,\tau,\theta)$ can be obtained by finding a feasible solution $(\bar\bv,\bar\bu,\{\bar{\bepsilon}^k\}_{k\in \cK})$ for model~\eqref{eq:IO_mip} and creating a conic combination of $\ba^i$'s for $i$ such that $\bar{v}_i=1$.} 
In fact, without having to normalize $\bc$ post-hoc, for any conic combination of such $\ba^i$'s, there is an inverse-feasible cost vector for $\mathbf{QIO}(\cK,\tau,\theta)$ that generates the same forward optimal solutions. That is, if we let $\cC_{QIO}$ be the set of all inverse-feasible $\bc$ vectors for $\mathbf{QIO}(\cK,\tau,\theta)$ and $\cC_{MIP} = \bigcup_{\bv \in \hat\cV}\textup{cone}(\{\ba^i\}_{  i:  {v}_i=1})$ where $\hat\cV$ is the set of feasible non-zero $\bv$ vectors for model~\eqref{eq:IO_mip}, we have $\bigcup_{\bc\in\cC_{QIO}} \text{argmin}\,\mathbf{FO}(\bc) = \bigcup_{\bc\in\cC_{MIP}} \text{argmin}\,\mathbf{FO}(\bc)$. The following result formally characterizes the set of all inverse-feasible solutions for $\mathbf{QIO}(\cK,\tau,\theta)$ (proof is straightforward from the proof of Theorem~\ref{prop:IO_mip} and thus is omitted).
\begin{corollary}
\label{remark:IO_mip}
Let $\Pi(\hat{\mathcal{X}})$ denote the feasible region of model \eqref{eq:IO_mip}. Then the set of all inverse-feasible cost vectors, $\hat\cC$, can be characterized by $\hat\cV = \text{Proj}_{\bv}(\Pi(\hat\cX))\setminus \{\bzero\}$. That is, $\hat\cC = \bigcup_{\bv\in\hat\cV} \textup{cone}(\{\ba^i\}_{i:v_i=1}).$ 
\end{corollary}

\subsubsection{Lower Bound on the Inverse Stability Measure.}
\rev{
Finally, we show that a solution to \eqref{eq:IO_mip} provides a lower bound on inverse stability measure  \eqref{eq:inv_stab} for $\mathbf{QIO}(\cK,\tau,\theta)$. To do so, we make the following assumptions:
(i) data noise in $\hat\cX$ is in the form of point-wise shift, i.e., a shifted data point is expressed as $\tilde\bx^k=\hat\bx^k+\boldsymbol{\delta}^k$ for each $k\in \cK$; (ii) $\displaystyle d(\hat\cX,\tilde\cX)=\sum_{k\in \cK} \|\hat\bx^k- \tilde\bx^k\|_\ell$ for some $\ell \ge 1$. 
\begin{proposition}\label{prop:IS_LB}
Let $\xi^{*}$ be the inverse stability value for $\mathbf{QIO}(\cK,\tau,\theta)$. Let $(\bar\bv\ne\mathbf{0},\bar\bu,\{\bar{\bepsilon}^k\}_{k\in \cK})$ be a feasible solution to \eqref{eq:IO_mip}, $\bar\cI=\{i\in\cI\,|\,\bar{v_i}=1\}$, $\displaystyle d_{ik}^{*} = \min_{\bx\in\cX_i}\{\|\hat\bx^{k}-\bx\|_\ell
\},i\in\bar\cI,k\in\cK$, and $\displaystyle d_{i[k]}^{*}$ denote the $k$-th largest value of $\displaystyle d_{ik}^{*},k\in\cK$. Then we have $\displaystyle \max_{i\in\bar\cI}\bigg\{ \sum_{k=1}^{\lfloor(1-\theta)K\rfloor+1} \max(0,\tau-d_{i[k]}^{*})\bigg\} \le \xi^{*}$.
\end{proposition}
Given $\cK$ and $\tau$, we note that the lower bound increases as $\theta$ decreases, which reinforces the role of $\theta$ in increasing inverse stability. Our numerical results show that the lower bound is non-trivial and close to the numerically estimated inverse stability measure.
}

\rev{
While problem \eqref{eq:IO_mip} leads to the set of cost vectors with improved inverse stability, which cost vector among them to use for the forward problem is an important consideration to improve forward stability. In the next subsection, we introduce an objective function for problem \eqref{eq:IO_mip} to find a $\bc$ vector with improved forward stability.
}

\vspace{-0.1in}
\subsection{Finding the Maximal Dimension Inverse-Feasible Set: Improving Forward Stability}\label{sec:model2}
%
Because it is hard to model forward stability measure \eqref{eq:for_stab} as a tractable objective function, we propose a surrogate, tractable function that can be used as an objective function for the QIO model. Recall that the forward stability measure for a cost vector $\bc$ represents the worst-case distance between $\cX^*\!(\bc)=\underset{\bx}{\text{argmin}}\,\textbf{FO}(\bc)$ and the data. \rev{We compute forward stability only in terms of the data points that are actually used by the model and consider $\bc$ to be forward-stable if the worst-case distance is no greater than the error threshold $\tau$ for the relevant subset of the data of cardinality $\lceil\theta K\rceil$ (i.e., excluding those deemed to be outliers).} We first make the following observation.
\rev{
\begin{proposition}\label{prop:FS_bound_E}
If a feasible solution $(\bar\bv,\bar\bu,\{\bar\bepsilon^k\}_{k\in\cK})$ for problem~\eqref{eq:IO_mip} satisfies $\displaystyle\sum_{i\in\cI} \bar{v}_i = n$, then $\bar\bc\in\textup{cone}_+(\{\ba^i\}_{i:\bar{v}_i=1})$ satisfies $\underset{\bx \in \cX^*\!(\bar\bc)}\max  \{\|\hat\bx^k-\bx\|_{\ell}\} \le \tau, \forall k:\bar{u}_k=1$.
\end{proposition}
Proposition~\ref{prop:FS_bound_E} implies that if $\bc$ is a strict conic combination of $n$ $\ba^i$'s 
chosen by \eqref{eq:IO_mip} (i.e., such that $\bar{v}_i=1$), then $\mathbf{FO}(\bc)$ has a unique optimal solution that is within $\tau$-distance from all chosen data points \citep{mangasarian1979uniqueness,tavasli2018}, and thus the forward stability measure is guaranteed to be no greater than the pre-specified threshold. Furthermore, the following result shows that as we assign more 1's to the elements of $\bv$ in problem \eqref{eq:IO_mip} (i.e., adding more basis vectors to the conic combination), the resulting $\bc$ vector improves {\it{both}} forward  and inverse stability. 
\begin{proposition}\label{prop:MQIO-FWStab}
Given $\hat\cX$, let $(\bar\bv, \bar\bu, \{\bar\bepsilon^k\}_{k\in \cK})$ and $(\tilde \bv, \bar\bu, \{\tilde\bepsilon^k\}_{k\in \cK})$ be two distinct feasible solutions to problem \eqref{eq:IO_mip}. Let $\bar\cC=\textup{cone}(\{\ba^i\}_{i: \bar v_i=1})$ and $\tilde\cC=\textup{cone}(\{\ba^i\}_{i: \tilde v_i=1})$. If $\bar\bv \ge \tilde\bv$, then
\begin{enumerate}
\item [(i)] $\bar\cC$ is at least as forward-stable as $\tilde\cC$; i.e., for any $\tilde\bc \in \tilde\cC$ there exists $\bar\bc\in \bar\cC$ such that $\underset{\bx\in \cX^*(\bar\bc)}{\max}\{\|\hat\bx^k-\bx\|_\ell\}\le \underset{\bx\in \cX^*(\tilde\bc)}{\max}\{\|\hat\bx^k-\bx\|_\ell\},\forall k\in\cK$.
\item [(ii)] $\bar\cC$ is at least as inverse-stable as $\tilde\cC$; i.e., $\underset{\cX\subset\mathbb{R}^{K\times n}}{\min}\{d(\hat\cX,\cX) \,|\, \bar\cC\cap \cC(\cX)=\emptyset\} \ge \underset{\cX\subset\mathbb{R}^{K\times n}}{\min}\{d(\hat\cX,\cX) \,|\, \tilde\cC\cap \cC(\cX)=\emptyset\}$ where $\cC(\cX)$ denotes the set of inverse-feasible cost vectors for $\cX$.
\end{enumerate}
\end{proposition}
}

\rev{
Propositions~\ref{prop:FS_bound_E} and \ref{prop:MQIO-FWStab} imply that we can add the objective function of maximizing the number of non-zero $v_i$'s to problem~\eqref{eq:IO_mip} to find which cost vector to eventually use for the forward problem to further improve forward stability without compromising inverse stability. This leads us to the following mixed-integer quantile inverse optimization (MQIO) model:}
\begin{equation}\label{eq:IO_mip_maxV}
\mathbf{MQIO}(\cK,\tau,\theta):\quad
\underset{\bv,\bu,\{\bepsilon^k\}_{k\in\cK}}{\text{maximize}} \; \Big\{\,\sum_{i\in\cI}v_i \, \Big{|} \,\text{\eqref{eq:IO_mip.2}--\eqref{eq:IO_mip.5}}\Big\}.
\end{equation}
Our numerical results show that even when the optimal value of \eqref{eq:IO_mip_maxV} is not $n$ but close to $n$, the MQIO model is more forward-stable than the previous inverse model. \rev{
Computing forward stability is inherently challenging as it involves finding the maximum distance between  the dataset and some face of polyhedron $\cX$. An upper bound for the forward stability measure may be computed if some geometric information about the polyhedron is given, e.g., the minimum volume ellipsoid or other size measure for some faces of $\cX$, yet obtaining such information itself is not easy in general. In Appendix~\ref{appx:FS_upperbound}, we show that an optimal solution to the MQIO model leads to an upper bound on the forward stability measure if such information is available. 
} 

\rev{Appendix~\ref{appx:Extensions} provides extensions to the MQIO model that can accommodate situations where $(\bA,\bbb)$ varies over different data points, $(\bA,\bbb)$ is subject to uncertainty, or the cost vector is constrained.}  
The MQIO problem is a large-size MIP in general and thus is computationally challenging. In the next section, we provide efficient solution approaches by exploiting the problem structure.

\section{Solution Approaches}\label{sec:model_algorithms}
\rev{
In this section, we develop efficient solution approaches to the MQIO problem. In Section~\ref{sec:Problem_Complexity}, we first show that the MQIO problem is NP-hard by establishing its equivalence to a biclique problem. We then provide a formal biclique representation of a feasible solution to the QIO (hence MQIO) problem in Section~\ref{sec:Biclique}, which indicates that algorithms designed for biclique problems can be adapted to solve the MQIO problem. Finally, we propose an exact algorithm and heuristics in Section~\ref{sec:Algorithms}
.
}
\vspace{-0.15in}
\subsection{Problem Complexity and Connection to Biclique Problems}\label{sec:Problem_Complexity}
We show that \textbf{MQIO}$(\cK,\tau,\theta)$ is NP-hard by establishing its connection to biclique problems, which we later use to develop alternative algorithms for the problem. To do so, we first show that there exists an optimal solution to $\mathbf{MQIO}(\cK,\tau,\theta)$ with exactly $\lceil \theta K \rceil$ observations selected. 
\begin{lemma}\label{lemma 1}
 	If $\mathbf{MQIO}(\cK,\tau,\theta)$ is feasible, then there exists an optimal solution $( \bv^*,\bu^*,\{\bepsilon^{k*}\}_{k\in \cK})$ such that $\sum_{k=1}^{ K} u^*_k = \lceil \theta K \rceil$.
\end{lemma}
\begin{proposition} \label{prop:np-hard}
 	The problem $\mathbf{MQIO}(\cK,\tau,\theta)$ is NP-hard. 
\end{proposition}
To prove Proposition~\ref{prop:np-hard} we first introduce the maximum $\kappa$-subset intersection (MSI) problem. Consider a ground set $\cG=\{g_1,\ldots,g_N\}$ and a set of its subsets $\mathscr{R}=\{\cR_1,\ldots,\cR_Q\}$, i.e., $\cR_q \subseteq\cG, q=1,\ldots,Q$. Given a positive integer $\kappa$, MSI finds $\kappa$ subsets in $\mathscr{R}$ whose intersection has maximum cardinality.  
MSI problems are known to be NP-hard \citep{dawande2001,xavier2012}.
\proof{Proof of Proposition \ref{prop:np-hard}:}
Consider the following feasibility problem for each $k \in \cK$.
 %
\vspace{-0.15in}
\begin{subequations}\label{eq:IO_MIP_Reduced}	
\begin{align}	
\underset{\bv,\bepsilon}{\text{maximize}} &\quad  0 \label{eq:IO_MIP_Reduced1}\\	
\text{subject to} 
&\quad b_i \leq {\ba^i}'(\bxhat^k-\bepsilon) \leq b_i +M(1-v_i), \quad \forall i  \in \cI,\label{eq:IO_MIP_Reduced2}\\	
&\quad \lVert \bepsilon \rVert_{\ell}\leq \tau, \label{eq:IO_MIP_Reduced3}\\
&\quad  v_i \in \{0,1\}, \quad\forall i  \in \cI. \label{eq:IO_MIP_Reduced4}
\end{align}	
\end{subequations}

Let $(\bv^{k},\bepsilon^k)$ 
be a feasible solution to \eqref{eq:IO_MIP_Reduced} with respect to $\hat\bx^k$ and $\cA_k=\{i \in \cI \,|\, v_i^{k}=1\}$ be a subset of $\cI$ for each $k \in \cK$, and construct the set of subsets $\mathscr{A}=\{\cA_1,\ldots,\cA_{ K}\}$. Finding exactly $\lceil \theta K \rceil$ subsets in $\mathscr{A}$ such that their intersection has maximum cardinality is equivalent to solving MSI with $\kappa=\lceil \theta K \rceil$. Note that for each $k\in\cK$ there can be multiple solutions for $\bv^k$ satisfying \eqref{eq:IO_MIP_Reduced}, each leading to a different $\cA_k$. As a result, $\mathscr{A}$ may not be unique. Thus, $\mathbf{MQIO}(\cK,\tau,\theta)$ is equivalent to solving the MSI problem multiple times with different $\mathscr{A}$'s and finding the maximum cardinality. Therefore, $\mathbf{MQIO}(\cK,\tau,\theta)$ is at least as hard as MSI.~\Halmos
\endproof
\cite{dawande2001} show that a general MSI problem can be reformulated as a version of the biclique problem (called the maximum one-sided edge cardinality problem).

\vspace{-0.15in}
\subsection{\rev{Biclique Representation of the Solution}} \label{sec:Biclique} 
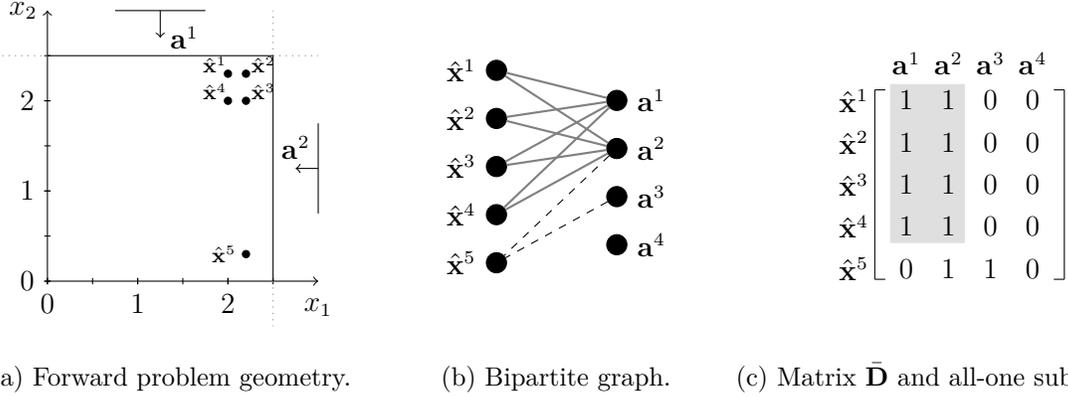
\begin{figure}[t]\centering
\begin{subfigure}[b]{.37\textwidth}
\centering
\vfill
\begin{tikzpicture}[scale=1.2]
\draw[->] (0,0) -- coordinate (x axis mid) (3,0) node[pos=1,below=0.1cm] {$x_1$};
    \foreach \x in {0,0.5,...,2.5}
\draw  (\x,1pt)--(\x,-1pt);
    \foreach \x in {0,...,2.5}
\draw  (\x,1pt)--(\x,-1pt)
node[anchor=north] {\x};
    \draw[->] (0,0) -- coordinate (y axis mid) (0,3) node[pos=1,left] {$x_2$};
    \foreach \y in {0,0.5,...,2.5}
\draw  (1pt,\y)--(-1pt,\y);
    \foreach \y in {0,...,2.5}
\draw  (1pt,\y)--(-1pt,\y)
node[anchor=east] {\y};
\draw (2.5,0) -- (2.5,2.5);
\draw[gray,dotted] (2.5,0)--(2.5,-0.5);
\draw[gray,dotted] (2.5,2.5)--(2.5,3);
\draw (2.5,2.5) -- (0,2.5);
\draw[gray,dotted] (-0.5,2.5) -- (0,2.5);
\draw[gray,dotted] (2.5,2.5) -- (3,2.5);
\draw (0.75,3) -- (1.75,3); 
\draw[->] (1.25,3) coordinate -- (1.25,2.7) 
    node[anchor=west] {$\ba^1$};
    
    \draw (3,0.75) -- (3,1.75); 
\draw[->] (3,1.25) coordinate -- (2.75,1.25)  node[anchor=south] {$\ba^2$};
    \draw[fill=black] (2,2.3) circle (0.04);
    \node [black,anchor= east] at (2.1,2.4) {\scriptsize$\hat\bx^1$};
    
    \draw[fill=black] (2.2,2.3) circle (0.04);
    \node [black,anchor= west] at (2.15,2.4) {\scriptsize$\hat\bx^2$};
    
    \draw[fill=black] (2.2,2) circle (0.04);
    \node [black,anchor= west] at (2.15,2.1) {\scriptsize$\hat\bx^3$};
    
    \draw[fill=black] (2,2) circle (0.04);
    \node [black,anchor= east] at (2.1,2.1) {\scriptsize$\hat\bx^4$};
    
    \draw[fill=black] (2.2,0.3) circle (0.04);
    \node [black,anchor= east] at (2.2, 0.3) {\scriptsize$\hat\bx^5$};
\end{tikzpicture}
\caption{
Forward problem geometry.}
\label{fig:Prop5-LP_Region}
\end{subfigure}
\hspace{0.00in}
\begin{subfigure}[b]{.22\textwidth}
\centering
\vfill
\begin{tikzpicture}[scale=0.8]
  \node [white] at (1, 1) {\textbullet};
  \tikzset{Bullet/.style={circle,draw,fill=black,scale=0.75}}
\node[Bullet,label=left :{$\hat\bx^1$}] (x1) at (0,5) {} ;
\node[Bullet,label=left :{$\hat\bx^2$}] (x2) at (0,4.2) {} ;
\node[Bullet,label=left :{$\hat\bx^3$}] (x3) at (0,3.4) {} ;
\node[Bullet,label=left :{$\hat\bx^4$}] (x4) at (0,2.6) {} ;
\node[Bullet,label=left :{$\hat\bx^5$}] (x5) at (0,1.8) {} ;
\node[Bullet,label=right :{$\ba^1$}] (a1) at (2,4.5) {} ;
\node[Bullet,label=right :{$\ba^2$}] (a2) at (2,3.7) {} ;
\node[Bullet,label=right :{$\ba^3$}] (a3) at (2,2.9) {} ;
\node[Bullet,label=right :{$\ba^4$}] (a4) at (2,2.1) {} ;
\draw[black,thick,gray] (a1)--(x1);
\draw[black,thick,gray] (a1)--(x2);
\draw[black,thick,gray] (a1)--(x3);
\draw[black,thick,gray] (a1)--(x4);
\draw[black,thick,gray] (a2)--(x1);
\draw[black,thick,gray] (a2)--(x2);
\draw[black,thick,gray] (a2)--(x3);
\draw[black,thick,gray] (a2)--(x4);
\draw[black,dashed] (a2)--(x5);
\draw[black,dashed] (a3)--(x5);
\end{tikzpicture}
\vfill
\caption{Bipartite graph.}
\label{fig:Prop5-Bipartite}
\end{subfigure}
\hspace{0.03in}
\begin{subfigure}[b]{.37\textwidth}
\centering
\vfill
\begin{tikzpicture}[scale=0.7]
  \node [white] at (1, 1) {\textbullet};
  \node [black] at (0,5) {$\hat\bx^1$};
  \node [black] at (0,4.2) {$\hat\bx^2$};
  \node [black] at (0,3.4) {$\hat\bx^3$};
  \node [black] at (0,2.6) {$\hat\bx^4$};
  \node [black] at (0,1.8) {$\hat\bx^5$};
   
  \node [black] at (1,5.7)   {$\ba^1$};
  \node [black] at (1.8,5.7) {$\ba^2$};
  \node [black] at (2.6,5.7) {$\ba^3$};
  \node [black] at (3.4,5.7) {$\ba^4$};
  \draw (0.4,5.2) -- (0.4,1.6);
  \draw (4.0,5.2) -- (4.0,1.6);
  
  \draw (0.4,5.2) -- (0.6,5.2);
  \draw (0.4,1.6) -- (0.6,1.6);
  
  \draw (4.0,5.2) -- (3.8,5.2);
  \draw (4.0,1.6) -- (3.8,1.6);
  \draw[gray!25,fill=gray!25]   (0.7,2.3) rectangle (2.1,5.3);
  \node [black] at (1,5)   {$1$};
  \node [black] at (1,4.2) {$1$};
  \node [black] at (1,3.4) {$1$};
  \node [black] at (1,2.6) {$1$};
  \node [black] at (1,1.8) {$0$};
  
  \node [black] at (1.8,5)   {$1$};
  \node [black] at (1.8,4.2) {$1$};
  \node [black] at (1.8,3.4) {$1$};
  \node [black] at (1.8,2.6) {$1$};
  \node [black] at (1.8,1.8) {$1$};
  
  \node [black] at (2.6,5)   {$0$};
  \node [black] at (2.6,4.2) {$0$};
  \node [black] at (2.6,3.4) {$0$};
  \node [black] at (2.6,2.6) {$0$};
  \node [black] at (2.6,1.8) {$1$};
   
  \node [black] at (3.4,5)   {$0$};
  \node [black] at (3.4,4.2) {$0$};
  \node [black] at (3.4,3.4) {$0$};
  \node [black] at (3.4,2.6) {$0$};
  \node [black] at (3.4,1.8) {$0$};
\end{tikzpicture}
\vfill
\caption{Matrix $\bar\bD$ and all-one submatrix.}
\label{fig:Prop5-Matrix}
\end{subfigure}
\caption{
\rev{
Illustrative example for the construction of $\bar\bD$ and its all-one submatrix
}
}
\label{fig:Prop5}
\end{figure}

Motivated by its connection to the MSI problem, we also cast the MQIO problem as a biclique problem with a bipartite graph constructed as follows. We first create a node $k$ in $V_1$ for each data point $\bxhat^k$ and a node $i$ in $V_2$ for each 
$\ba^i$. With every node $k \in V_1$ and $i \in V_2$ we associate an edge $e=(k,i) \in \mathbb{E}$ if $v_i^{k}=1$ where $\bv^{k}$ denotes a feasible solution to problem \eqref{eq:IO_MIP_Reduced} with respect to $\bxhat^k$. Then we define a bi-adjacency matrix $\bar\bD \in 
\{0,1\}^{ K\times m}$ associated with graph $G=(V_1 \cup V_2,\mathbb{E})$ as follows.
\begin{align}\label{eq:D_bar}
  \bar\bD_{ki} = \begin{cases} 1, 
               & \mbox{if } e=(k,i) \in \mathbb{E}, \\ 0, 
               & \mbox{otherwise.} 
                 \end{cases}
\end{align}
Matrix $\bar\bD$ can be built by solving problem \eqref{eq:IO_MIP_Reduced} for each $k\in\cK$ to find a feasible solution $\bv^k$ and letting the $k$-th row of $\bar\bD$ be ${\bv^k}'$. 
The following result shows how an inverse-feasible solution can be found from the $\bar\bD$ matrix above.
\begin{proposition}\label{prop:Dbar-IO}
If there exists $\bar\bD$ that satisfies \eqref{eq:D_bar} and has an all-one submatrix  whose rows and columns correspond to $\bar\cS$ and $\bar\cA$, respectively, where $\bar\cS \subseteq \cK$, $|\bar \cS| \ge \theta K$, and $\bar\cA \subseteq \cI$, then there exists a solution $(\bar{\bc},\{\bar{\bepsilon}^k\}_{k\in \cK},\bar\by,\bar{\cS})$ feasible for $\mathbf{QIO}(\cK,\tau,\theta)$ where $\bar\bc \in \textup{cone}(\{\ba^i\}_{ i \in \bar\cA})$.
\end{proposition}
\rev{
Proposition~\ref{prop:Dbar-IO} suggests that if we can find a matrix $\bar\bD$ that has an all-one submatrix with at least $\lceil \theta K\rceil$ rows, there exists a corresponding inverse-feasible cost vector for the QIO (hence MQIO) problem. The following example further illustrates this.
}
\rev{
\begin{example}\label{ex:InvSol_Clique}
Consider Example~\ref{Example:FW_Instability} again (see Figure~\ref{fig:Prop5-LP_Region}), and the QIO model with the following parameters: $\tau=1$, $\theta=0.8$, $\|\cdot\|_\ell=\|\cdot\|_\infty$, and $\|\cdot\|_p=\|\cdot\|_\infty$. Let $(\bar\bv^k,\bar\bepsilon^k)$ be a feasible solution to problem \eqref{eq:IO_MIP_Reduced} 
for each data point $k$: we have $\bar\bv^k = [1,1,0,0]'$ for $k=1,2,3,4,$ with $\bar\bepsilon^1=\begin{bsmallmatrix}-0.5\\-0.2\end{bsmallmatrix}$, $\bar\bepsilon^2=\begin{bsmallmatrix}-0.3\\-0.2\end{bsmallmatrix}$, $\bar\bepsilon^3=\begin{bsmallmatrix}-0.3\\-0.5\end{bsmallmatrix}$,  and $\bar\bepsilon^4=\begin{bsmallmatrix}-0.5\\-0.5\end{bsmallmatrix}$, respectively, and $\bar\bv^5=[0,1,1,0]'$ with $\bar\bepsilon^5=\begin{bsmallmatrix}-0.3\\0.3\end{bsmallmatrix}$. Then we can construct a bipartite graph, shown in Figure~\ref{fig:Prop5-Bipartite}, where each edge identifies $(i,k)$ such that $\bar v_i^k=1$, and the corresponding $\bar\bD$ matix in Figure~\ref{fig:Prop5-Matrix}. Note that this $\bar\bD$ has an all-one submatrix with rows and columns corresponding to $\bar\cS=\{1,2,3,4\}$ and $\bar\cA=\{1,2\}$, respectively (the shaded rectangle), and $|\bar\cS| = 4 \ge \theta K$. Now consider a cost vector $\bar\bc = 0.5\ba^1+0.5\ba^1 = \begin{bsmallmatrix}-0.5\\-0.5\end{bsmallmatrix}\in \textup{cone}(\{\ba^i\}_{i\in\bar\cA})$. We can see that $\bar\bc$ is inverse-feasible because $(\bar\bc,\{\bar\bepsilon^1, \ldots, \bar\bepsilon^5\}, \bar\by, \bar\cS)$ where 
$\bar\by=[0.5,0.5,0,0]'$ is feasible for \textbf{QIO}($\cK,\tau,\theta$).
\end{example}
}

\rev{
Furthermore, Proposition 7 implies that the number of columns in the all-one submatrix 
corresponds to the objective function value of $\mathbf{MQIO}(\cK,\tau,\theta)$
, which is to be maximized. Thus, it is important to find a submatrix with as many columns as possible. Building on these observations, in the next section, we propose specific algorithms to find an all-one submatrix of $\bar\bD$ with at least $\lceil \theta K\rceil$ rows and as many columns as possible.
}

\subsection{\rev{Algorithms}}\label{sec:Algorithms}
\rev{
While the matrix $\bar\bD$ is constructed using $\{\bv^k\}_{k\in\cK}$ obtained by solving problem~\eqref{eq:IO_MIP_Reduced} for each data point $k$, because there may be multiple solutions for \eqref{eq:IO_MIP_Reduced}, there may be multiple $\bar\bD$ matrices achievable. Thus, it is important to find a $\bar\bD$ matrix that can lead to an all-one submatrix of a desirable size as described above, which eventually leads to a desirable feasible solution to $\mathbf{MQIO}(\cK,\tau,\theta)$. In this subsection, we now present our main algorithm that finds such a $\bar\bD$ matrix and corresponding inverse-feasible cost vectors.}

\rev{
Our first idea is that we solve problem~\eqref{eq:IO_MIP_Reduced} with a weighted objective function for each $k$ such that resulting solutions $\{\bv^k\}_{k\in\cK}$ are similar to each other, thus leading to similar rows of $\bar\bD$. For the first data point, i.e., $k=1$, we solve problem~\eqref{eq:IO_MIP_Reduced} with respect to $\hat\bx^1$ with the objective function $\displaystyle\sum_{i\in \cI} v_i$, which will find a solution $\bv^1$ with as many ones as possible. For $k\ge 2$, we solve problem~\eqref{eq:IO_MIP_Reduced} with respect to $\hat\bx^k$ with the weighted objective function $\displaystyle\sum_{i\in \cI}w_i^{k} v_i$ where $\bw^{k}=\bvartheta^{k}/\|\bvartheta^{k}\|$ and  $\displaystyle\vartheta^{k}_i=\sum_{\ell=1}^{(k-1)} v^{\ell}_i$; i.e., we assign a higher weight to column $i$  for which $v_i^{\ell}$ was assigned 1 more often for the previous data points $\ell=1,\ldots,k-1$. The resulting vector ${\bv^{k}}'$ is then inserted to the $k$-th row of $\bar\bD$. Note that because the weights are updated in a sequential manner, the order of the data points can impact the resulting $\bar\bD$. We assume without loss of generality that the data points are sorted in the increasing order of the sum of the distances to the rest of the dataset (i.e., $\displaystyle\sum_{k\in\cK} \|\hat\bx^1-\hat\bx^k\|_\ell \le \ldots \le \sum_{k\in\cK} \|\hat\bx^K-\hat\bx^k\|_\ell $); we find this ordering particularly effective because starting with data points that are closer to the rest of the dataset leads to more similar $\bv^k$'s in the subsequent iterations.
}

\rev{
Next, since the construction of $\bar\bD$ depends on the solution $\bv^1$ for problem~\eqref{eq:IO_MIP_Reduced} with respect to the first data point $\hat\bx^1$ and this problem may have multiple solutions, it is important to consider different solutions for this problem to place different ``seeds'' for $\bar\bD$ (i.e., different first rows of $\bar\bD$), which allows to explore different inverse solutions. To do so, once an optimal solution to problem~\eqref{eq:IO_MIP_Reduced} (with the added objective function) for the first data point $\hat\bx^1$ is found, we impose a cut to exclude the solution from the feasible region and re-solve the problem repeatedly. 
}

\rev{
Finally, once a certain $\bar\bD$ matrix is found, its all-one submatrix 
(i.e., a biclique in the bipartite graph corresponding to $\bar\bD$) with at least $\lceil \theta K\rceil$ rows and as many columns as possible can be found (if exists) by solving the associated maximum clique problem. While any existing algorithm for a general maximum clique problem would work, we use an MIP formulation to solve the problem exactly, which we call Clique($\bar\bD,\theta$); this formulation can be found in Appendix~\ref{sec:clique}.
}

\rev{
This algorithm is formally presented in Algorithm~\ref{alg:Dbar-alg_Exact}. In Lines 4-9, problem~\eqref{eq:IO_MIP_Reduced} is solved with the weighted objective function and cuts. To generate different first rows of $\bar\bD$ (i.e., different $\bv^1$'s), a cut defined in Line 11 is added to the problem in Line 5. The following result shows that, under a certain condition, this algorithm finds inverse-feasible cost vectors that correspond to the optimal solution to \textbf{MQIO}($\cK,\tau,\theta$), hence exact. We call this algorithm $\bar\bD$-Alg-Exact.
}
\rev{
\begin{proposition}\label{prop:exactness2}
Let $(\bv^*,\bu^*,\{\bepsilon^{k*}\}_{k\in \cK})$ be an optimal solution to $\mathbf{MQIO}(\cK,\tau,\theta)$ and $z^*=\sum_{i\in \cI}v^*_i$. Let $\cV_z=\{\bv^1\in \{0,1\}^{m} \,|\, \exists(\bv^1,\bepsilon^1)$ satisfying \eqref{eq:IO_MIP_Reduced2}--\eqref{eq:IO_MIP_Reduced4} with respect to $\hat\bx^{1}\!, \; {\bv^1}'\beee=z\}$. If $u^*_1=1$, then Algorithm~\ref{alg:Dbar-alg_Exact} returns $z^{Alg}=z^*$ and $\bv^{Alg}=\bv^*$ in $r\le \displaystyle\sum_{z= z^*}^n |\cV_z|$ iterations. 
\end{proposition}
%
Once the algorithm is run, the condition $u_1^*=1$ can be easily checked by whether the first row of $\bar\bD_{r}$ is also included in its all-one submatrix where $r$ is the iteration at which the most recent $\bv^{Alg}$ was achieved. Our numerical results show that the condition is met in all instances for which the algorithm was run completely (not interrupted by time limit). If this condition is not met, although the algorithm may still generate a reasonably good solution, the user may reorder the data points to place a different data point in the position of $k=1$ and rerun the algorithm. 
}

\rev{
{\OneAndAHalfSpacedXI
\begin{algorithm}[t]
\rev{	\caption{$\bar\bD$-Alg-Exact}
	\label{alg:Dbar-alg_Exact}	\SetAlgoLined
	\KwResult{$z^{Alg}$, $\bv^{Alg}$, $\cC^{Alg}$}
	\KwIn{$\bA, \bbb, \hat\cX,
	 \tau,\theta$}  
	$r\leftarrow 1, \; \cG_1 \leftarrow \emptyset, \; z^{Alg}\leftarrow -\infty$\\
	\While{$|\cG_r| \ge z^{Alg}$}{
    $ \bw^1 \leftarrow \beee, \; \bar\bD_r \leftarrow [\ ]$\\
	\For{$k=1,\ldots,K$}{ 
	$\text{Find} \; \begin{cases}
	\bv^{k}\in \underset{\bv}{\text{argmax}} \{{\bw^k}'\bv \,|\, \eqref{eq:IO_MIP_Reduced2}-\eqref{eq:IO_MIP_Reduced4}\}, \quad \text{if}\ r=1 \\
	\bv^{k}\in \underset{\bv}{\text{argmax}} \{{\bw^k}'\bv \,|\, \eqref{eq:IO_MIP_Reduced2}-\eqref{eq:IO_MIP_Reduced4},\displaystyle \sum_{i\in\cG_{r'}} v_i\le |\cG_{r'}|-1, \forall r'=2,\ldots,r\}, \quad \text{otherwise} 
    \end{cases}$
with $\hat\bx^{k}$ as input data
  \\
		Insert ${\bv^{k}}'$ to the $k$-th row of $\bar\bD_r$ \\
		$\vartheta^{k}_i \leftarrow \sum_{\ell=1}^{k} v^{\ell}_i$ for each $i\in \cI$\\
		Update weights $\bw^{k+1}\leftarrow \bvartheta^{k}/\|\bvartheta^{k}\|$
	}
	Solve Clique($\bar\bD_r,\theta)$ to find the all-one submatrix of $\bar\bD_r$. 
	Define $\bar\bv^r\in\{0,1\}^{m}$ where 1's correspond to columns of the submatrix. 
	\\
	$\cG_{r+1} \leftarrow \{i\ |\ v_i^1 = 1 \;\text{ (i.e., $(1,i)$-th entry of $\bar\bD_r$ is 1)\,\}}$ 
	\\
	\If{$\bar\bv^{r'}\beee >z^{Alg}$}
	{$z^{Alg}\leftarrow \bar\bv^{r'}\beee$\\
	 $\bv^{Alg}\leftarrow \bar\bv^{r}$}
	$r\leftarrow r+1$
	}
	Construct the set of inverse-feasible cost vectors $\cC^{Alg} = \textup{cone}(\{\ba^i\}_{i: v^{Alg}_i=1})$\\ }
\end{algorithm} 
}
}

\rev{
Algorithm \ref{alg:Dbar-alg_Exact} without cuts, i.e., finding only one $\bar\bD$, would be faster, albeit not exact, which we call $\bar\bD$-Alg-Heuristic. We also propose an even faster heuristic called $\tilde\bD$-Alg-Heuristic in Appendix~\ref{appx:algs}, which solves a linear relaxation of \eqref{eq:IO_MIP_Reduced} and uses weighted $\ell_1$ minimization to make $\mathbf{1}\!-\!\bv$ sparse where $\bv\in [0,1]^m$. The proposed algorithm framework can lead to many different variants depending on how the weights are updated and how cuts are generated. For example, the algorithms can be modified to accommodate a huge dataset by applying the sequential weight updating only to the first few MIPs (say $\bar K$ of them), which leads to the weight parameters $\bw^{\bar K}$, and solving the remaining ($K\!-\!\bar K$) MIPs in parallel with the same weights $\bw^{\bar K}$. Also, the algorithms can be modified to infer cost vectors in an online manner based on datasets that become available through different time points and $(\bA,\bbb)$ parameters changing over time; Appendix~\ref{sec:model_online} shows more details about this algorithm.}

\section{Numerical Results} \label{sec:result}

We examine the performance of the algorithms proposed for the MQIO model using various-sized LP instances. We also assess forward and inverse stability of the MQIO model as well as the previous inverse LP method using the measures proposed in Section \ref{sec:prelim}. Finally, we demonstrate the stability performance of the MQIO model in the diet recommendation and transshipment applications.

\subsection{Performance of the Algorithms}
We evaluate the performance of the 
proposed algorithms for randomly generated LP instances with $n\in \{15, 50\}$, $m\in \{100,300\}$, $K \in \{35,200, 500\}$, and $\tau\in \{ 3, 3.5,4\}$. An instance is defined by a tuple $(n,m,K,\tau)$ and each instance was solved twice---with $\theta=0.75$ and $0.85$. For each instance, as a pre-processing step, we excluded $\ba^i$'s that were not inverse-feasible \emph{a priori} from the construction of $\bc$ (see Algorithm~\ref{alg:IO1}). 
\begin{table}[t]
\begin{threeparttable}
\centering
\caption{\rev{MQIO results via MIP solver, $\bar\bD$-Alg-Exact, $\bar\bD$-Alg-Heuristic, and $\tilde\bD$-Alg-Heuristic.}}
\label{tab:AlgResults}
\renewcommand{\arraystretch}{0.75}
\begin{tabular}{crrrrrrrrrr}
    \toprule
          &       & \multicolumn{4}{c}{Objective Function Value} &       & \multicolumn{4}{c}{Time (s)} \\
\cmidrule{3-6}\cmidrule{8-11}    $(n,m,K)$ & $\tau$     & \multicolumn{1}{c}{MIP}   & $\bar\bD$-Exact & $\bar\bD$-Heur & $\tilde\bD$-Heur  &       & \multicolumn{1}{c}{MIP}  & $\bar\bD$-Exact &  $\bar\bD$-Heur & $\tilde\bD$-Heur \\
    \midrule
    \multirow{3}[2]{*}{(15,100,35)} & 3.0     & *14.00 & 14.50 & 11.50 & 11.00 &       & 5000.00 & 12.42 & 1.32  & 0.18 \\
          & 3.5   & 15.00 & 15.00 & 13.00 & 10.50 &       & 2441.31 & 26.61 & 0.72  & 0.20 \\
          & 4.0     & 15.00 & 15.00 & 12.00 & 9.50  &       & 28.62 & 5.77  & 0.60  & 0.20 \\
    \midrule
    \multirow{3}[2]{*}{(15,100,200)} & 3.0     & *12.50 & 12.50 & 11.00 & 11.50 &       & 5000.00 & 103.03 & 2.77  & 0.90 \\
          & 3.5   & 15.00 & 15.00 & 14.00 & 10.50 &       & 766.59 & 77.94 & 2.35  & 1.05 \\
          & 4.0     & 15.00 & 15.00 & 15.00 & 13.00 &       & 426.08 & 2.18  & 2.18  & 0.88 \\
    \midrule
    \multirow{3}[2]{*}{(50,300,35)} & 3.0     & *26.50 & 37.50 & 30.50 & 25.50 &       & 5000.00 & 2000.15 & 43.26 & 0.77 \\
          & 3.5   & *28.00 & 44.50 & 36.50 & 27.00 &       & 5000.00 & 3511.62 & 34.28 & 0.83 \\
          & 4.0     & *24.50 & 49.00 & 41.00 & 28.50 &       & 5000.00 & 2181.47 & 26.31 & 0.85 \\
    \midrule
    \multirow{3}[2]{*}{(50,300,200)} & 3.0     & *0.00  & 33.00 & 27.00 & 20.00 &       & 5000.00 & 1393.53 & 145.27 & 4.23 \\
          & 3.5   & *0.00  & 40.50 & 29.00 & 26.00 &       & 5000.00 & 2811.66 & 182.90 & 4.75 \\
          & 4.0     & *0.00  & 47.50 & 33.00 & 27.00 &       & 5000.00 & 2705.12 & 140.94 & 5.33 \\
    \midrule
    \multirow{3}[2]{*}{(50,300,500)} & 3.0     &       & 27.00 & 27.00 & 10.50 &       & 5000.00 & 384.28 & 384.28 & 10.23 \\
          & 3.5   &       & 35.00 & 29.50 & 13.50 &       & 5000.00 & 775.12 & 375.80 & 12.21 \\
          & 4.0     &       & 43.50 & 37.00 & 15.50 &       & 5000.00 & 3027.96 & 186.59 & 12.84 \\
    \bottomrule
    \end{tabular}
\begin{tablenotes}
	\footnotesize
	\item *Instances not solved to optimality within time limit (5000 seconds). Empty cells indicate instances where no feasible solution was found within time limit. 
  \end{tablenotes}
 \end{threeparttable}
\end{table}
Table \ref{tab:AlgResults} shows the results for the MQIO problem \eqref{eq:IO_mip_maxV} obtained by the Gurobi solver \citep{gurobi} (labeled ``MIP''), \rev{$\bar\bD$-Alg-Exact, $\bar\bD$-Alg-Heuristic, and $\tilde\bD$-Alg-Heuristic,} averaged over two results with $\theta=0.75$ and $0.85$. Recall that the objective function value represents the number of basis vectors ($\ba^i$'s) to construct the set of inverse-feasible cost vectors. All optimization problems were solved by Gurobi 7.5 with a 4-core 3.6 GHz processor and 32 GB memory.

For small instances where $(n,m)=(15,100)$, the exact MIP formulation found an optimal solution within the time limit of 5,000 seconds for most cases. Note that the solution time for MIP decreases as $\tau$ goes up---we conjecture that this is because the problem becomes easier as the threshold distance $\tau$ becomes more ``generous.'' 
\rev{
The exact algorithm achieved the same objective values as the MIP values within significantly less time. Both heuristics found objective values reasonably close to the exact values and were even faster. The computational benefit of the proposed algorithms becomes more clear for larger instances; the MIP solver failed to solve any instance with $(n,m) = (50,300)$. 
As expected, there is a clear trade-off between the algorithms. In all instances, $\bar\bD$-Alg-Exact returned the greatest objective values but was the slowest, whereas $\tilde\bD$-Alg-Heuristic was the fastest algorithm yet the solution quality was not as good as the other two algorithms. Although the results suggest that both heuristics can be used for very large instances, since $\tilde\bD$-Alg-Heuristic solves LP relaxations for each $\hat\bx^k$ instead of MIPs, we suggest using this algorithm for instances with very large $n$ and $m$.
}

\subsection{Stability of MQIO}
\subsubsection{Forward Stability.}
To assess the forward stability of our proposed MQIO model, we created an instance with $(n,m,K,\tau)=(15,100,35,3)$ and $\theta=0.75$, for which we knew the true optimal objective value was $\sum_{i \in \cI}v^*_i=15$. To see the effect of the number of $\ba^i$'s used to form an inverse solution (i.e., $\sum_{i \in \cI}v_i$) on forward stability, we solved the inverse problem on this instance repeatedly, each with a constraint $\sum_{i \in \cI}v_i \le h$, $h \in \{1,2,\ldots,15\}$. Each problem was solved exactly via the MQIO formulation. Note that the case where $h=1$ can be considered similar to the previous inverse LP model \eqref{eq:IO_prev} in that it finds a single cost vector identical to one of the $\ba^i$'s. For each $h$, 
once the set of inverse-feasible cost vectors, $\cC_h$, was found, we randomly generated 50 cost vectors $\bc\in\cC_h$ and computed $d(\hat\cX,\bx^*(\bc))= \underset{k\in\cS_h}{\max} \, \{\|\hat\bx^k - \bx^*(\bc)\|_{\infty}\}$, where $\cS_h$ denotes the set of chosen data points in iteration $h$ and $\bx^*(\bc)\in\argmin\textbf{FO}(\bc)$, as an estimate for the forward stability measure.

\begin{figure}[b]
\begin{subfigure}[b]{.48\textwidth}
	\includegraphics[width=3.1in]{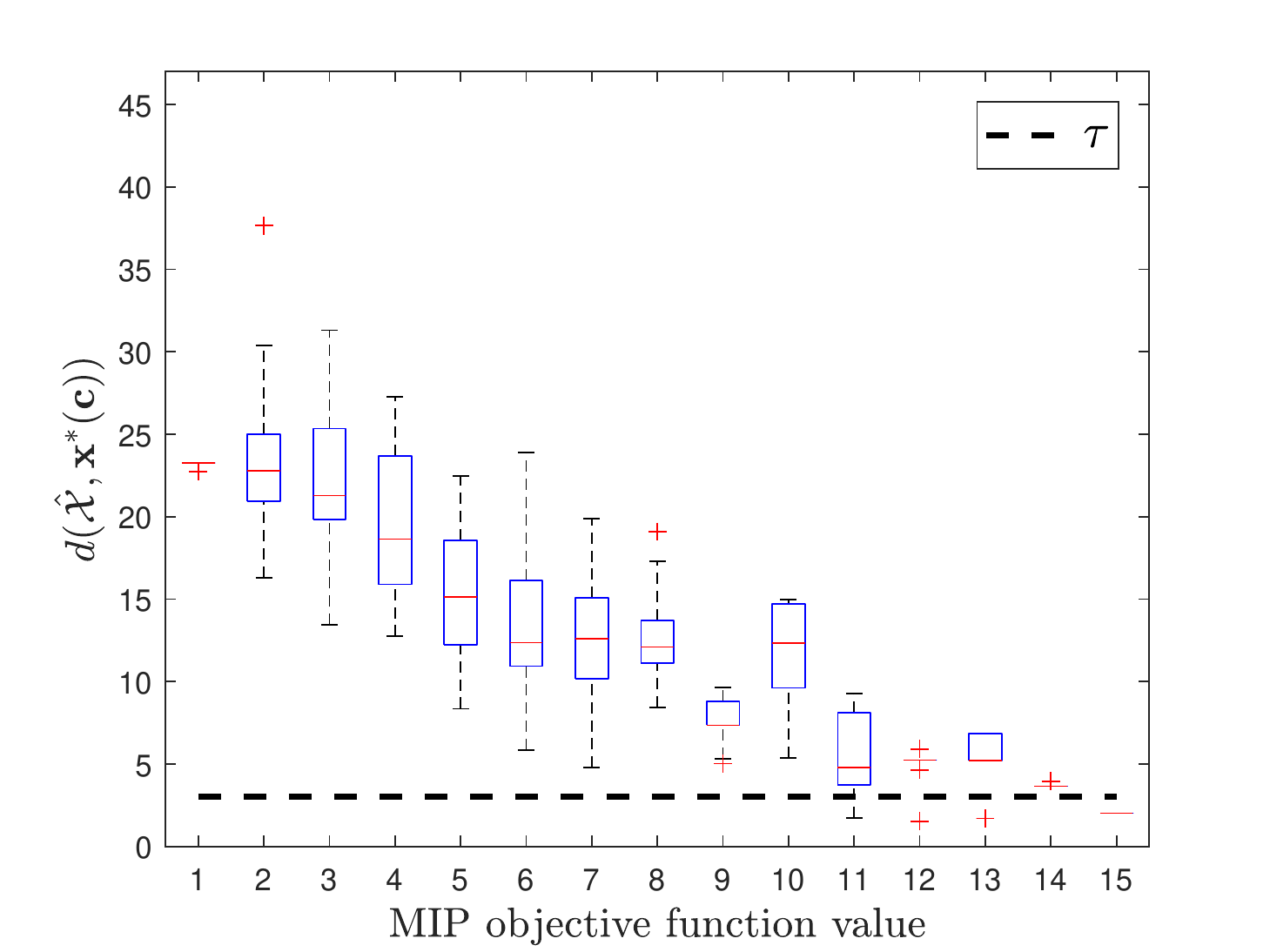}
	\caption{Distances between the dataset and solutions achieved by inverse-feasible cost vectors.}
	\label{fig:FWStab}
\end{subfigure}
\hspace{0.12in}
\begin{subfigure}[b]{.41\textwidth}
\includegraphics[width=3.1in]{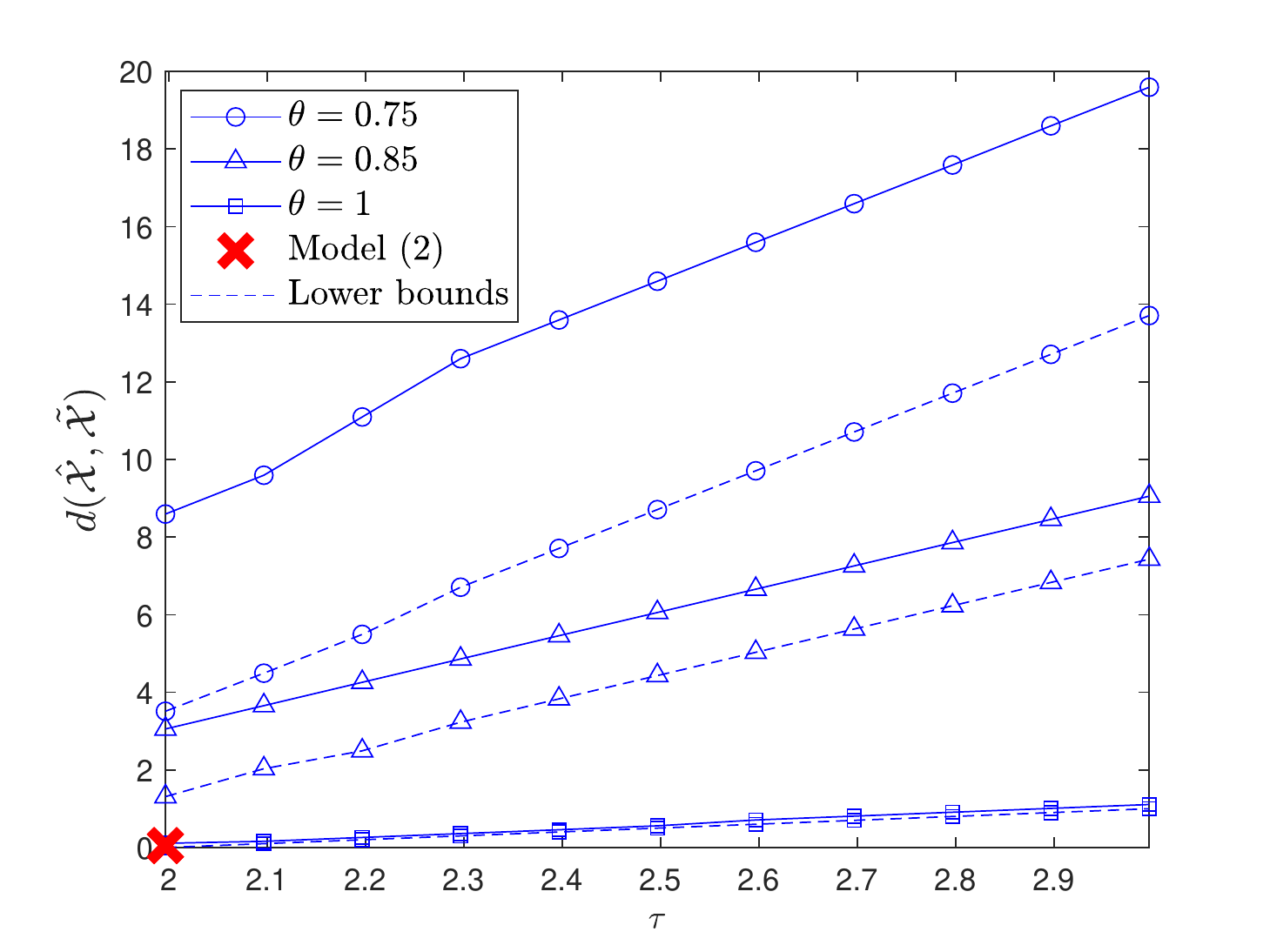}
	\caption{\rev{Data shift distance until no initial inverse-feasible solution remains feasible.}}
	\label{fig:InvStab}
\end{subfigure}	
\caption{
Forward and inverse stability performance of MQIO.
}
\label{fig:stability}
\end{figure}

Figure~\ref{fig:FWStab} shows the resulting distances for each $h \in \{1,2,\ldots,15\}$: each box-plot represents 50 values of $d(\hat\cX,\bx^*(\bc))$ from 50 different cost vectors. We observe that when $h$ is low (i.e., few $\ba^i$'s are used for generating a $\bc$ vector) the distance $d$ is significantly greater than $\tau$, hence  forward-unstable. In other words, if this cost vector is used for the forward problem, the resulting solutions may not be close to the input data. On the other hand, as $h$ increases (i.e., as more $\ba^i$'s are used) the cost vector produces a solution that is closer to the data, eventually within the threshold distance $\tau$ (indicated by the horizontal dashed line). The overall decreasing trend of the distance $d$ supports the idea of maximizing $\sum_{i\in \cI }v_i$ as a surrogate objective function to maximize forward stability. 
%

\vspace{-0.1in}
\subsubsection{Inverse Stability.}\label{sec:IS_Experiment}
To evaluate the inverse stability of MQIO, we investigated the effect of $\theta$ and 
$\tau$ on how sensitive the inverse solution set is to changes in the dataset $\hat\cX$. Again, an instance with $(n,m,K)=(15,100,35)$ was considered. First, we considered the MQIO model with different $\tau$ values increasing from 2 to 3 by 0.1 (the objective value of MQIO increases monotonically as $\tau$ increases). For each $\tau$ value, we first solved the MQIO model with the given dataset $\hat\cX$ to find the initial inverse-feasible set $\hat\cC$. 
\rev{
We then shifted $\lfloor(1-\theta)K\rfloor+1$ data points in $\hat\cX$: $\{\hat\bx^k+\bgamma^k\}_{k\in\cK}$ 
where $\bgamma^k$ is uniformly distributed in $[0,\Gamma]^n$ 
if $\hat\bx^k$ is chosen to be shifted, and $\bgamma^k=\bzero$ otherwise. Given this shifted data, denoted by $\tilde \cX(\Gamma)$, we re-solved the MQIO problem to find the new inverse-feasible set $\tilde\cC$. 
}
We repeated this process while increasing $\Gamma$ until $\tilde\cC$ had no common cost vector shared by the initial set $\hat\cC$, i.e., $\hat\cC \cap \tilde \cC = \emptyset$. When such a $\Gamma$ value was reached, we computed the distance $\displaystyle d(\hat\cX,\tilde\cX(\Gamma))=\sum_{k\in\cK} \|\hat\bx^k - \tilde\bx^k(\Gamma)\|_{\infty}$ as an estimate for the inverse stability measure in Section~\ref{sec:prelim}. 


Figure \ref{fig:InvStab} shows the result of this experiment for all $\tau$ values and $\theta=0.75, 0.85$ and $1$. For a fixed value of $\theta$, the distance $d(\hat\cX,\tilde\cX(\Gamma))$ increases as $\tau$ increases (i.e., more $\ba^i$ vectors are included in the inverse-feasible set), which means the model becomes more tolerant to changes in the dataset, hence more inverse-stable. For a fixed value of $\tau$, Figure \ref{fig:InvStab} also shows that increasing $\theta$ decreases inverse stability, reinforcing the idea that allowing for no or few outliers can lead to an ``impatient'' inverse model that can easily switch its solutions affected by such data imperfection. 

\rev{
Note that the estimated inverse stability values in Figure \ref{fig:InvStab} are in fact upper bounds on the true inverse stability values. To show this, we refer the reader to formulation \eqref{eq:Exact_IS} in Appendix \ref{appx:LB}, whose optimal value is the true inverse stability value for the MQIO model. Since no initial cost vector in $\hat\cC$ remains inverse-feasible given the shifted data $\tilde\cX$, $\tilde\cX$ is a feasible solution to formulation \eqref{eq:Exact_IS}, with the objective value of $d(\hat\cX,\tilde\cX)$. Because \eqref{eq:Exact_IS} is a minimization problem, its optimal value (i.e., the true inverse stability value) cannot be greater than $d(\hat\cX,\tilde\cX)$. Figure \ref{fig:InvStab} shows that the lower bounds on the inverse stability measure (based on Proposition \ref{prop:IS_LB}) are non-trivial, and since the estimated values are upper bounds, the difference between the lower bounds and the true inverse stability values must be even less than shown in Figure \ref{fig:InvStab}. 
}

Finally, we conducted a similar experiment to assess the inverse stability of the previous inverse model \eqref{eq:IO_prev}. We solved model \eqref{eq:IO_prev} with $\hat\cX$ to find an initial inverse solution $\hat\bc$, re-solved the problem repeatedly with $\tilde\cX(\Gamma)$ 
with increasing $\Gamma$ values until the resulting cost vector was different from $\hat\bc$, and computed the distance between the original dataset and the shifted one, which is shown as a thick $\boldsymbol{\times}$ marker in Figure \ref{fig:InvStab}. The result is almost identical to the lowest inverse stability achieved by the MQIO model (i.e., with $\theta=1$ and smallest $\tau$).

\vspace{-0.1in}
\subsection{Diet Recommendation}
Literature on diet prediction/recommendation using historical data largely focuses on a ``direct'' replication of diets where the goal is to create diets that are closest to the old diets in terms of some distance function (e.g., $\ell_2$ distance between new and old diet vectors) \citep{perignon2016low,darmon2006impact}. However, such an approach may not work if there is any change in the constraints of the underlying diet problem such as changes in nutritional requirements or available foods. In this case, learning the objective function that represents one's preferences instead and creating new diets ``indirectly'' is a more robust and transferable way of learning the individual's diet behavior (similar motivations can be found in learning driving behaviors in \cite{Abbeel_2004}). Moreover, by doing so the model can find an \textit{optimal} diet with respect to the inferred preferences, whereas replicating the diets directly can inherit undesirable diet patterns from the past.

We use the database from the National Health and Nutrition Examination Surveys (NHANES), which includes nutritional requirements and nutrition facts per serving for each food type, and build a diet problem for a subset of the foods and nutrients from the database to keep the experiment simple. We classify the foods into food ``types'' (see Table \ref{tab:DietData} in Appendix~\ref{appx:Diet_Data}). We assume that once the number of servings for each food type is determined, more detailed decisions (e.g., specific menus) can be made by dietitians based on this. We believe this is a realistic consideration as there exist a myriad of different menus. We first solve the diet problem with some arbitrary cost vector, assumed to reflect true preferences of an individual, and find $\bx^*$ that represents the number of servings for food type $i$ per day. \rev{In practice, if the data come from some unknown system, the assumed forward LP feasible region may not ``fit'' the data well. In this case, the user needs to tweak some parameters of the forward LP (e.g., constraint parameters or the list of variables), similar to linear regression.}

\rev{Given $\bx^*$ with no noise, both the MQIO and previous models can find the true cost vector and generate the same $\bx^*$ as a forward solution. However, as one's diet behavior is assumed inconsistent and noisy, to examine the impact of data noise on both models we generate multiple diets $\bx^*+\bgamma$} where $\bgamma$ is a normal-distributed random noise vector $\bgamma \sim \cN(\bzero,\sigma^2 \bI)$ where $\bI$ is the identity matrix of an appropriate dimension, which form the input dataset $\hat\cX$.

\begin{figure}[t]
	\centering
	\begin{subfigure}[t]{0.50\textwidth}
		\centering
		\includegraphics[width=3.2in]{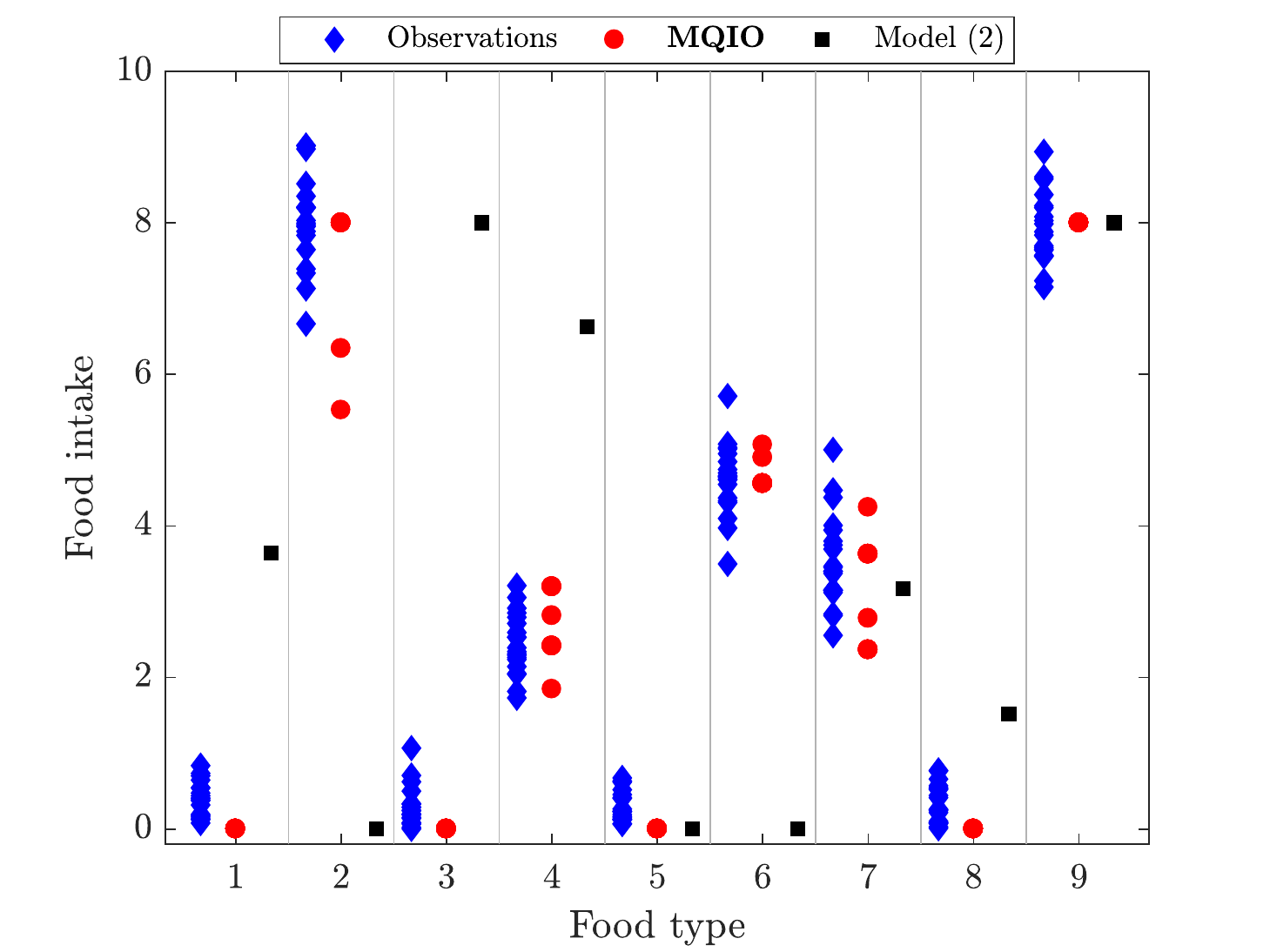}
		\caption{
		$\sigma=0.2$}
		\label{fig:Diet1}
	\end{subfigure}%
	\begin{subfigure}[t]{0.48\textwidth}
		\centering
		\includegraphics[width=3.2in]{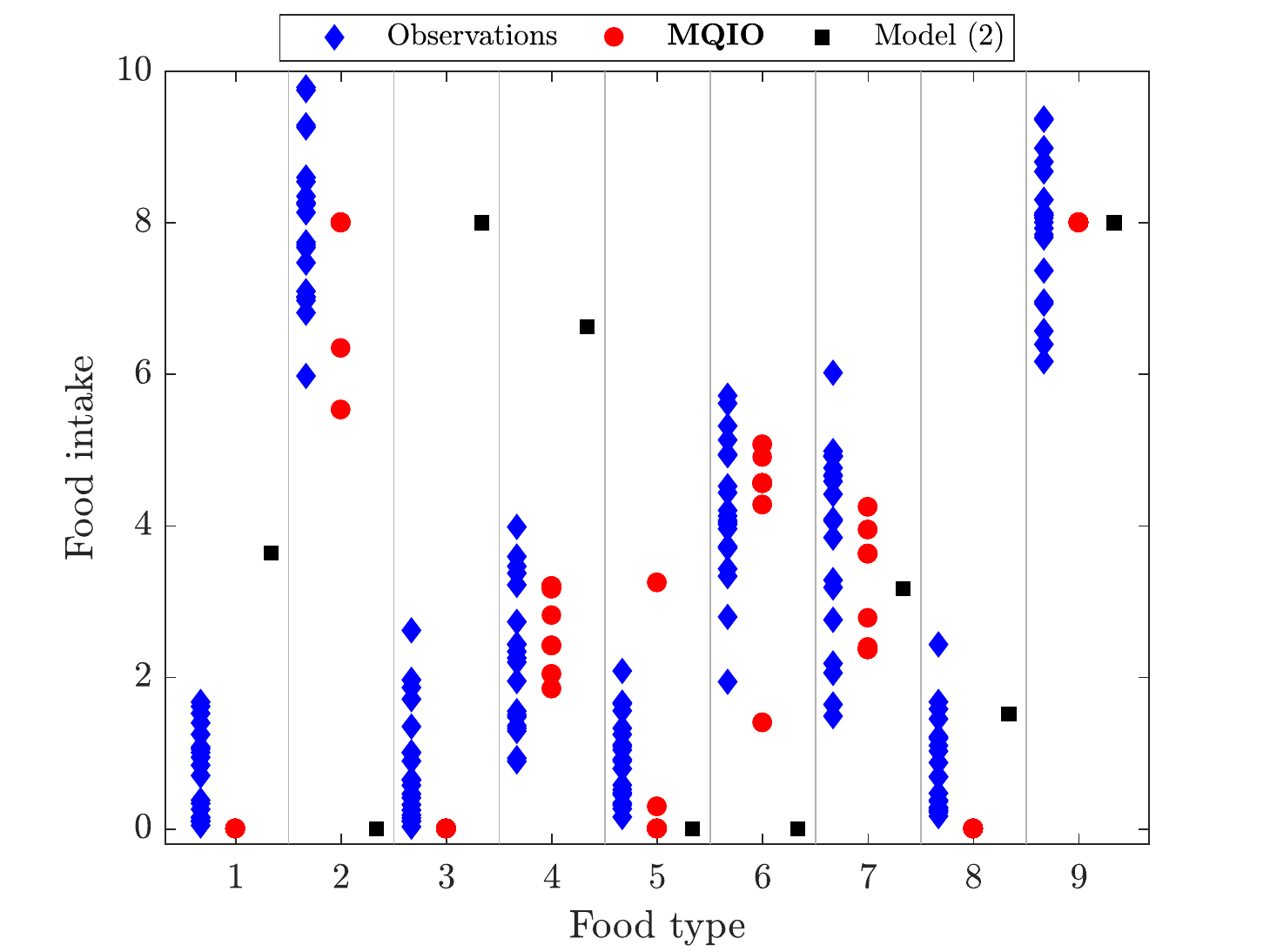}
		\caption{
		$\sigma=1.0$
		}
		\label{fig:DietVariability}
	\end{subfigure}%
	\caption{Comparison of diets recommended by the MQIO model and the previous inverse model.}
	\label{fig:DietReplication}
\end{figure}

We apply the MQIO model as well as the previous inverse model~\eqref{eq:IO_prev} with the same dataset $\hat\cX$ as input and compare the resulting recommended diets. For MQIO, we find a set of inverse-feasible cost vectors $\cC^*$ and solve the forward diet problem 20 times each with a different $\bc$ vector  selected from $\cC^*$ to generate a set of diet recommendations. Since the previous inverse model finds a single cost vector, we solve the diet problem with this vector only once. Figure \ref{fig:Diet1} shows diet recommendations obtained by the MQIO model (circular markers) are much closer to the input data (diamond-shaped markers; $\sigma=0.2$ to create the noisy dataset $\hat\cX$) than those from the previous inverse model are (squared markers), reinforcing the improved stability in the MQIO model. Also note that the MQIO model leads to multiple such diet recommendations based on the multiple cost vectors.

In Figure \ref{fig:DietVariability}, we further increase the noise in the input data by increasing $\sigma$ from 0.2 to 1. While the previous model generates the same recommendation, the MQIO model ``adapts to'' the increased variability in the observations and generates more diverse recommendations for some food types; e.g., 
food types 5--7
. In summary, the MQIO model generates recommendations that are consistent with the individual's past behavior and stable in the face of data noise; thus, the objective function found by the model can better predict one's eating behavior.
\begin{figure}[t]
	\centering
	\includegraphics[width=3.2in]{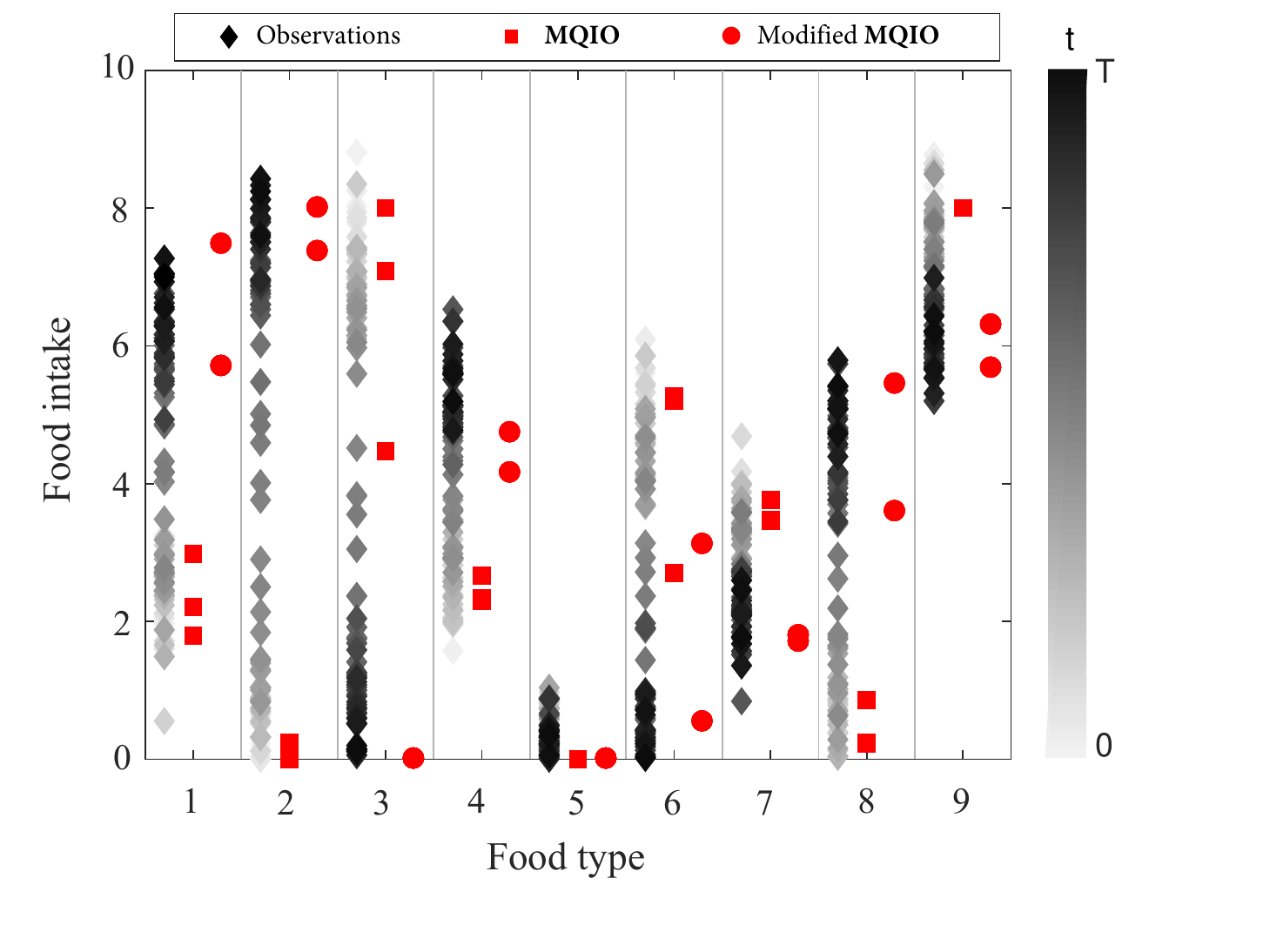}
	\caption{Results of the MQIO model for diet behavior changing over time.}
	\label{fig:DietTime}
\end{figure}

Finally, we show that the MQIO framework is amenable to the case where data point are collected at different time points and the goal of the inverse problem is to infer cost vectors based on more recent data points. Let $\hat \bx^t$ be a data point collected at time $t =1,\ldots,T$. We let the threshold parameter $\tau_t$ vary over time: the more recent $\hat\bx^t$ is, the smaller $\tau_t$ we impose (i.e., the closer we want to make this data point to optimality). Figure \ref{fig:DietTime} shows the result of such a modified model with time-stamped diet observations. 
Darker diamond markers represent more recent observations. 
As more recent data were collected, the modified MQIO model generates recommendations (circular markers) that are closer to the more recent ones. This result suggests that our MQIO framework can offer modeling flexibility that can lend itself to adaptive settings where the preference function can be updated over time.

\vspace{-0.1in}
\subsection{Transshipment Problem}\label{sec:results_transshipment}
We demonstrate the use of our inverse model for online learning in the transshipment application, where production and shipment costs are inferred adaptively as batches of decisions are observed sequentially over time. We consider a transshipment problem with one distribution node, two supply nodes, and two demand nodes (see Appendix~\ref{appx:Transportaion}). Suppose at each time $t$,  demand $\bd^t$ is revealed and the DM makes a set of decisions on production ($\bx^{(p)}$) and transshipment ($\bx^{(t)}$) to minimize the total cost  ${\bc^{(p)}}'\bx^{(p)}+{\bc^{(t)}}'\bx^{(t)}$ where $\bc^{(p)} $ and $\bc^{(t)}$ denote production and transshipment costs, respectively. Let $\bc = [{\bc^{(p)}}; {\bc^{(t)}}]$ and $\bx = [{\bx^{(p)}}; {\bx^{(t)}}]$. Let $\hat\cX_t$ denote the set of decisions observed at time $t$. We assume that there is a true cost vector $\bc^{true}$ (shown in Appendix \ref{appx:Transportaion}), and generate a dataset $\hat\cX_t$ by adding noise to an optimal solution with respect to $\bc^{true}$ and randomly generated demand $\bd^t\sim U[0,1.1]$. 

We use an online learning extension to 
$\bar\bD$-Alg-Heuristic (see Appendix~\ref{sec:model_online}) to infer the set of cost vectors $\cC_t$ at each time $t$ from the datasets collected through time $t$. We consider $T=600$ (i.e., 600 sequential batches of data) and $\theta = 0.75$. At each iteration $t=1,\ldots,600$, we randomly select 20 cost vectors $\bc^t \in \cC_t$ and solve the corresponding forward problem with demand $\bd^t$ to obtain a forward solution $\bar\bx^t$. We assess the performance of the algorithm at each time $t$ in terms of the distance between our forward solution $\bar \bx^t$ and data batch $\hat\cX_t$, i.e., $d(\hat\cX_t,\bar\bx^t) = \underset{\bx\in\hat\cX_t}{\max} \, \|\bx-\bar\bx^t\|_{\infty}$. 

\begin{figure}[b]
\begin{subfigure}{0.5\textwidth}
\includegraphics[width=3.3in]{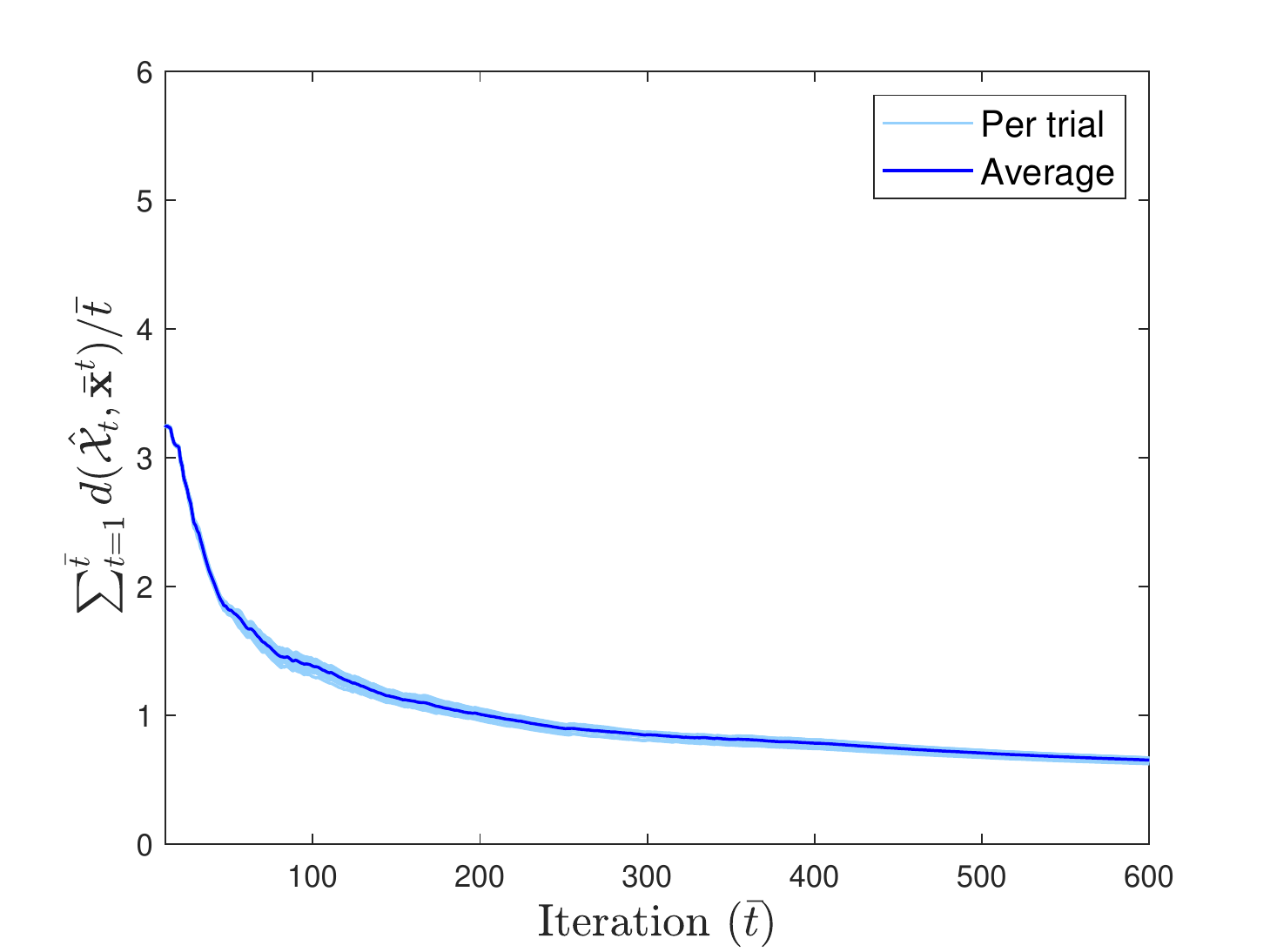}
\caption{Distance between inversely-optimized forward solution and the dataset.}
\label{fig:AvgDistRound1}
\end{subfigure}
\hspace{0.08in}
\begin{subfigure}{0.43\textwidth}
\includegraphics[width=3.3in]{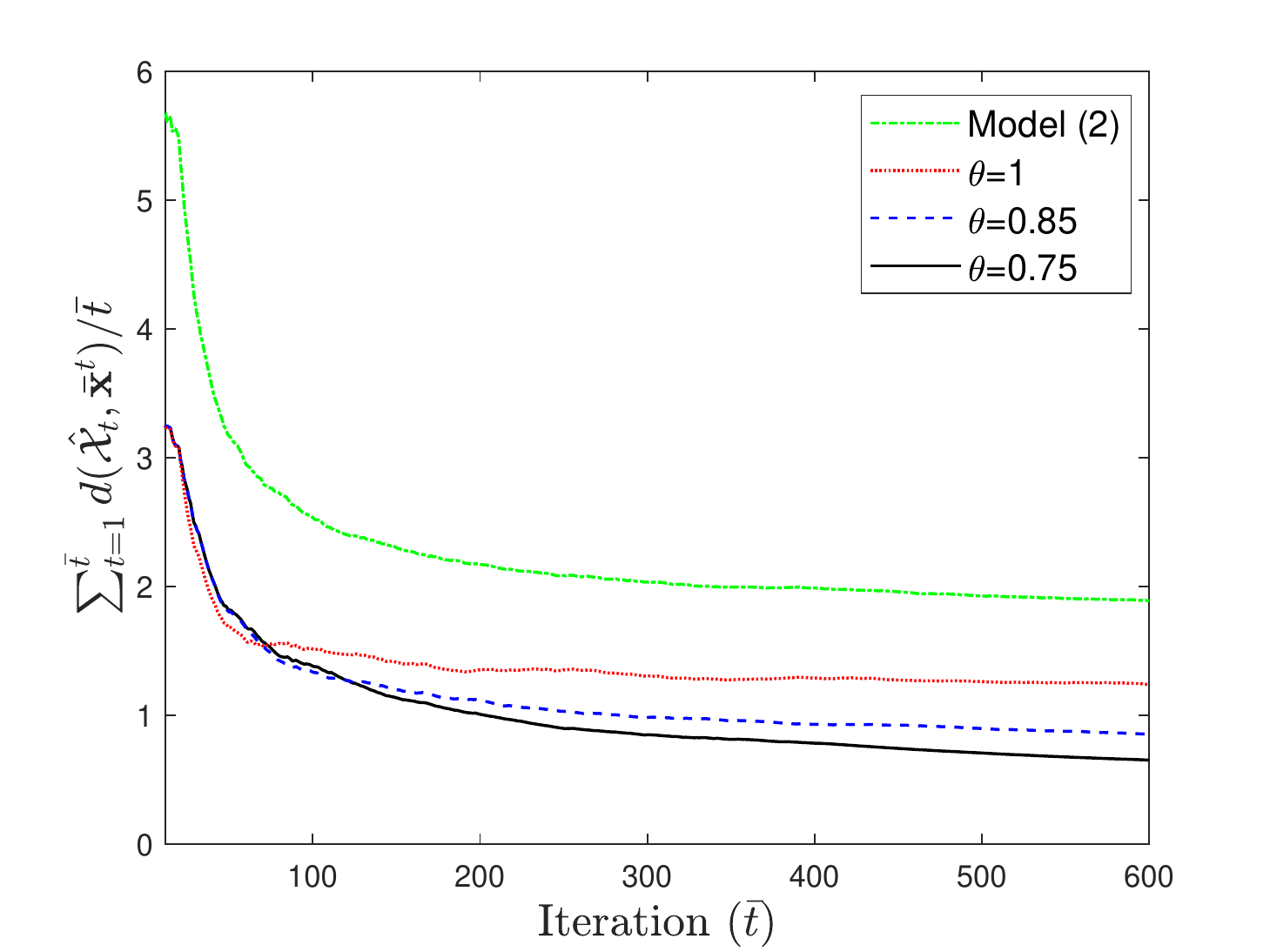}
\caption{Performance of the algorithm with different $\theta$ values and model \eqref{eq:IO_prev}.}
\label{fig:ThetaDiscussionSmodel2}
\end{subfigure}%
\caption{Results of the online learning extension for the transshipment problem.}
\label{fig:transshipment}
\end{figure}

Figure~\ref{fig:AvgDistRound1} shows the average distance $\sum_{t=1}^{\bar t} d(\hat\cX_t,\bar\bx^t)/\bar t$ at iteration $\bar t$ over 20 trials, i.e.,  20 randomly selected cost vectors. Each light-color line shows each distance per trial and the dark line shows the average. As $\bar t$ grows, the average distance between the dataset and the forward solution converges. In fact, the inferred set of cost vectors remains unchanged from iteration 415 onwards, and importantly, this set includes $\bc^{true}$
. That is, as more data are collected, the algorithm successfully finds a set of cost vectors that represent the preferences encoded in 
the true cost vector and generate forward solutions close to the given data (hence forward-stable).
\rev{
Additional results on the comparison between the objective function values achieved by our forward solution $\bar\bx^t$ and that from the given data are available in Appendix~\ref{appx:Transportaion}.
}

%
We also test the algorithm with different values of $\theta \in \{0.75,0.85,1\}$. Figure \ref{fig:ThetaDiscussionSmodel2} shows that as $\theta$ decreases, the algorithm takes more iterations to converge, which is intuitive as less data points are considered relevant at each iteration, requiring the algorithm to collect more data until it finally arrives at a certain set of cost vectors. When $\theta=1$, although the distance converges more quickly, the resulting inverse set does not include the true cost vector because the algorithm with $\theta=1$ is more susceptible to data noise. This again supports the need to explicitly account for potential outliers in inverse LP. Finally, Figure \ref{fig:ThetaDiscussionSmodel2} also shows the performance of the previous inverse model \eqref{eq:IO_prev} in this setting: at each iteration $\bar t$ we find a cost vector by solving \eqref{eq:IO_prev} with all data points up to time $\bar t$. This model fails to find the true cost vector at all iterations. Also, the distance between its corresponding forward solution and the data is greater than that from the MQIO model.  

The performance of the online learning algorithm (i.e., Algorithm \ref{alg:Online1}) also depends on the variability of demand parameters. 
When the interval of the demand distribution increases, the true vector is found more slowly. We conjecture that this is because increased variability in demands causes the feasible region of the forward problem to vary more,  rendering some constraints redundant and thus the algorithm cannot detect some of the $\ba^i$'s needed to construct the true cost vector.

\section{Conclusion} \label{sec:conclusion}
In this paper we developed a new inverse LP method that can capture noise, errors, and uncertainty in the input data and infer cost vectors that are more stable than those obtained by the previous methods. We formulated the model as a large-scale MIP and developed efficient algorithms by exploiting its connection to the well-known biclique problems. Our inverse method was demonstrated in the diet recommendation and transshipment applications where past data used as input can be noisy and inaccurate. Many directions for future research exist. For example, preferences of a large group of DMs can be clustered into smaller representative groups using the quantile statistic-based inverse optimization approach. 
\rev{
A goodness-of-fit measure that is amenable to the quantile inverse model can be useful for assessing the validity of the assumed forward LP, which remains our future work. Extending the inverse nonlinear programming techniques to accommodate the quantile-based framework is also an interesting future work, which can allow to assess our proposed inverse model in a broader context.
}

\clearpage

\begin{APPENDICES}
\section{Supplemental Materials for the Inverse Model}
\subsection{\rev{Finding the Threshold Parameter $\tau$}} \label{appx:finding_E}
The threshold parameter $\tau$ for the QIO model (i.e., \eqref{eq:IO1}) can be set by adding a reasonable margin to the minimum possible value that keeps the model feasible, which can be found by solving the following problem where $\tau$ is now a variable:
\begin{equation}\label{eq:IO1_minE}
\underset{\tau, \bc,\{\bepsilon^k\}_{k\in\cK},\by,\cS}{\text{minimize}} \big\{ \tau \  \big{|} \ \text{\eqref{eq:IO1.2}--\eqref{eq:IO1.9}} \big\}.
\end{equation}
The above problem~\eqref{eq:IO1_minE} finds a cost vector such that the $\theta$ quantile optimality error is minimized; this is similar to the structure of the least quantile method in the regression context in \cite{bertsimas2014least}. 
The following result shows the solution structure of an optimal cost vector for the above problem (proof can be found in Appendix~\ref{appx:proofs}).
\begin{proposition}\label{prop:IO1_minE}
\rev{
If problem \eqref{eq:IO1_minE} is feasible, then there exists an optimal ${\bc}$ for \eqref{eq:IO1_minE} such that ${\bc}=\ba^i$ for some $i \in \cI$.
}
\end{proposition}
Proposition~\ref{prop:IO1_minE} suggests that the search for an optimal cost vector for \eqref{eq:IO1_minE} can be done by evaluating each hyperplane $i=1,\ldots,m$, i.e., solving formulation \eqref{eq:Min_Fizbl_Dist} for each $i$ and each $k$ (see proof of Proposition \ref{prop:IO1_minE} in Appendix \ref{appx:proofs}) and finding $i$ that induces the minimum $\tau^*$; \rev{formulation~\eqref{eq:Min_Fizbl_Dist} is a convex optimization problem and thus is straightforward to solve.}

\rev{To further illustrate the implication of Proposition~\ref{prop:IO1_minE}, consider the forward feasible region and the initial data points in Example~\ref{Example:Inv_Instability_shift} (Figure \ref{fig:instab_inv_shift}) and problem~\eqref{eq:IO1_minE} with $\theta=0.5$ and $\|\cdot\|_\ell=\|\cdot\|_\infty$. Cost vector $\bc=\ba^1$ is the optimal solution to problem~\eqref{eq:IO1_minE}, leading to the objective value of $\tau^*=0.2$ with the optimality errors $\|\bepsilon^1\|_\infty=\|\bepsilon^2\|_\infty=0.2$; the objective function values for \eqref{eq:IO1_minE} achieved by $\bc=\ba^2$, $\bc=\ba^3$, and $\bc=\ba^4$ are $0.3$, $2$, and $2$, respectively.  
}


\subsection{Finding Big $M$ Parameters} \label{appx:Big-M}  
For each $k\in \cK$, consider the following formulation that finds $\bx\in \cX$ that has maximum $\ell_\infty$ distance from $\hat\bx^k$: 
\begin{equation}\label{eq:Max_Fizbl_Dist}
\underset{\bx}{\text{maximize}}\; \Big\{\,\lVert \bx-\hat\bx^k\rVert_{\infty} \, \Big{|} \, \bx\in \cX \Big\}.
\end{equation}
Let $d_k$ denote the optimal value of \eqref{eq:Max_Fizbl_Dist} with respect to $\hat\bx^k$. Given that $\cX$ is bounded, $d_k$ is also bounded. To find $d_k$, we first reformulate \eqref{eq:Max_Fizbl_Dist} as 
$\displaystyle \displaystyle\max_{j\in\cJ}\; \underset{\bx\in\cX}{\max}\;\{|x_j - \hat x_j^k|\} $; for each $j\in \cJ$, if $\bar \bx$ denotes an optimal solution to the inner max problem, then ${\max}\{|\bar x_j - \hat x_j^k|\} = \max\Big\{ \bar x_j-\hat x_j^k, \hat x_j^k - \bar x_j \Big\}$.
As a result, we have $\displaystyle d_k=\max_{j\in\cJ}\; \max \Big\{ \max_{\bx\in\cX}\; \{x_j - \hat x_j^k\}, \max_{\bx\in\cX}\; \{\hat x_j^k - x_j\} \Big\}$, i.e., $d_k$ can be obtained by solving $2n$ LPs and selecting the largest objective value. Note that $M_2=\underset{k\in \cK}{\max}\; \{d_{k}\}$ is a sufficiently large number valid for constraint \eqref{eq:IO_mip.3}. More efficiently, one may use a constraint-specific parameter $M_{2k}=d_k$ for each $k\in \cK$.

To find $M_1$, we use Cauchy-Schwarz and Minkowski inequalities as follows. For each $i\in \cI$ and $k\in \cK$, by Cauchy-Schwarz inequality, ${\ba^i}'(\bxhat^k-\bepsilon^k) \le \lVert \ba^i \rVert_2(\lVert \bxhat^k - \bepsilon^k \rVert_2)$. Using Minkowski inequality, we have $\lVert \bxhat^k - \bepsilon^k \rVert_2 \le \lVert \bxhat^k \rVert_2+\lVert \bepsilon^k \rVert_2$ for each $k\in \cK$. Moreover, for each vector $\bepsilon^k\in \mathbb{R}^n$ we have $\lVert\bepsilon^k\rVert_2 =
\sqrt{\sum_{i=1}^n \lvert\epsilon_i^k\rvert^2} \le
\sqrt{\sum_{i=1}^n (\underset{i}{\max}\{ \lvert\epsilon_i^k\rvert\})^2} = \sqrt{n\times(\underset{i}{\max}\{ \lvert\epsilon_i^k\rvert\})^2} =\sqrt{n}\lVert\bepsilon^k\rVert_\infty$. It follows that, ${\ba^i}'(\bxhat^k-\bepsilon^k) -b_i \le \lVert \ba^i \rVert_2(\lVert \bxhat^k - \bepsilon^k \rVert_2) -b_i \le \lVert \ba^i \rVert_2(\lVert \bxhat^k \rVert_2+\lVert \bepsilon^k \rVert_2)-b_i \le \lVert \ba^i \rVert_2(\lVert \bxhat^k \rVert_2+\sqrt{n} \lVert \bepsilon^k \rVert_\infty) -b_i$. From the previous paragraph, we have $\lVert \bepsilon^k \rVert_\infty \le 
d_k$ for each $k\in \cK$. Hence, choosing $M_1=\underset{i\in\cI}{\max}\{\lVert \ba^i \rVert_2\}\big(\underset{k\in\cK}{\max}\{\lVert\hat\bx^k\rVert_2\}+\sqrt{n}\,\underset{k\in\cK}{\max}\{d_k\}\big)+\underset{i\in\cI}{\max}\{-b_i\})$ is valid for each constraint of $\eqref{eq:IO_mip.4}$. As a more efficient and constraint-specific alternative, one may choose $M_{1ik}=\lVert \ba^i \rVert_2(\lVert \bxhat^k \rVert_2+\sqrt{n} d_k)-b_i$ for constraint $\eqref{eq:IO_mip.4}$ for each $i\in \cI$ and $k\in \cK$.

The $M$ parameter in formulation \eqref{eq:IO_MIP_Reduced}, with respect to each $k\in \cK$, can be calculated similarly to $M_1$ with $d_k$ replaced by $\tau$, i.e., $M=\underset{i\in\cI}{\max}\{\lVert \ba^i \rVert_2\}\big(\lVert\hat\bx^k\rVert_2+\tau\,\sqrt{n}\big)+\underset{i\in\cI}{\max}\{-b_i\})$.

\vspace{0.1in}

\subsection{\rev{Computing the Inverse Stability Measure}}\label{appx:LB}
\rev{
Recall that inverse stability measure \eqref{eq:inv_stab} is defined as the minimum distance of data shift such that all inverse-feasible cost vectors are no longer feasible for the inverse model with respect to the shifted data. Given a dataset $\hat\cX,$ let  $\hat\cC$ be the set of inverse-feasible cost vectors for \textbf{QIO}($\hat\cX,\tau,\theta$). As Theorem~\ref{prop:IO_mip} suggests, every inverse-feasible cost vector $\bc\in\hat\cC$ is a conic combination of some $\ba^{i}$'s that are also inverse-feasible. The following result shows that, to check whether or not all inverse-feasible cost vectors in $\hat\cC$ are no longer feasible under the shifted data, it is sufficient to check whether or not all inverse-feasible $\ba^{i}$'s whose conic hull forms $\hat\cC$ are no longer feasible.
\begin{lemma}\label{lemma:Sufficient_condition}
Let $\hat\cC$ be the set of inverse-feasible cost vectors for \textbf{\emph{QIO}}$(\hat\cX,\tau,\theta)$ and $\hat\cI = \{i\in\cI \,|\, \ba^{i}\in\hat\cC\}$. Let $\tilde\cX$ denote the shifted dataset and $\tilde\cC$ the set of inverse-feasible cost vectors for \textbf{\emph{QIO}}$(\tilde\cX,\tau,\theta)$. If no $\ba^i, i\in\hat{\cI}$, is inverse-feasible for \textbf{\emph{QIO}}$(\tilde\cX,\tau,\theta)$, then no $\bc\in\hat\cC$ is inverse-feasible for \textbf{\emph{QIO}}$(\tilde\cX,\tau,\theta)$.
\end{lemma}
Lemma \ref{lemma:Sufficient_condition} suggests that to compute the inverse stability measure for $\hat\cC$ given the original dataset $\hat\cX$, all we need is to find a minimally shifted set of the data points, $\tilde\cX$, such that no $\ba^{i},i\in\hat\cI$, is inverse-feasible for \textbf{QIO}($\tilde\cX,\tau,\theta$). Given $\hat\cI$, to make $\bc=\ba^{i}$ infeasible for \textbf{QIO}($\tilde\cX,\tau,\theta$) for all $i\in\hat\cI$, we need to ensure that there are more than $\lfloor (1-\theta)K \rfloor$ data points in $\tilde\cX$ whose minimum distance from each $\cX_i$, $i\in\hat\cI$, is strictly greater than $\tau$. This observation leads to the following formulation.
\begin{subequations}\label{eq:Exact_IS}	
\begin{align}
\underset{\{\bx^{k}\}_{k\in \cK},\{\cS^i\}_{i\in 
{\hat\cI}}}{\text{minimize}} &\quad  d(\hat\cX,\cX)=\sum_{k\in \cK}\|\bx^k-\hat\bx^k \|_\ell \label{eq:Exact_IS_1}\\
\text{subject to}
&\quad D^{i}(\bx^{k})=\underset{\bz\in\cX_i}{\min} \{\|\bx^k-\bz\|_{\ell}
\}, \quad \; \forall i\in 
{\hat\cI},  k\in \cK,\label{eq:Exact_IS_2}\\
& \quad D^{i}(\bx^{k}) \ge \tau+\sigma,  \quad \; \forall k\in \cS^i, i\in 
{\hat\cI},\label{eq:Exact_IS_3}\\& \quad |\cS^i|\ge \lfloor (1-\theta)K\rfloor +1, \; i\in 
{\hat\cI}. \quad \; \label{eq:Exact_IS_4}
\end{align}	
\end{subequations}
In the above formulation, variable $\bx^k$ is considered a shifted data point from $\hat\bx^k$. In constraint \eqref{eq:Exact_IS_2}, $D^{i}(\bx^{k})$ computes the minimum distance between some shifted data point $\bx^k$ and $\cX_i$, which is then ensured to be no less than $\tau+\sigma$ for at least $\lfloor (1-\theta)K\rfloor +1$ data points by constraints \eqref{eq:Exact_IS_3} and \eqref{eq:Exact_IS_4} where $\sigma$ denotes an infinitesimal positive constant. The objective function minimizes the sum of distances between the original data points and shifted ones. As computing the exact inverse stability value would have required $D^{i}(\bx^{k}) > \tau$ in place of constraint~\eqref{eq:Exact_IS_3}, the optimal value of the above problem is only infinitesimally greater than the exact inverse stability value, and we treat this value as the true inverse stability value (the inverse stability measure is defined on the one-dimensional open set). Note that problem~\eqref{eq:Exact_IS} is always feasible as one can let $\bx^{k},k\in\cK$, be arbitrarily far from each $\cX_i, i\in\hat\cI$. Problem \eqref{eq:Exact_IS} is non-convex and thus it is still challenging to compute its optimal value exactly. However, we use this formulation to compute the lower bound for the true inverse stability value, which is easy to compute (see Proposition \ref{prop:IS_LB} and its proof).
}

\vspace{0.1in}

\subsection{\rev{Upper Bound for the Forward Stability Measure}}\label{appx:FS_upperbound}
\rev{
In this section, we show that an optimal solution to the MQIO model leads to an upper bound on forward stability measure \eqref{eq:for_stab}. We make the following assumption: for each $i\in\cI$, there exists a constant $\rho_i$ such that $\|\bx^1-\bx^2\|_\ell\le \rho_i$ for any $\bx^1,\bx^2\in \cX_i$ and $\ell\ge 1$.
The following result shows that if such $\rho_i$ values are known \textit{a priori}, an upper bound for the forward stability measure for a cost vector found by the MQIO model can be computed. 
%
\begin{proposition}\label{lemma:FS_bound}
Let $(\bar\bv,\bar\bu,\{\bar\bepsilon\}_{k\in \cK})$ be a feasible solution to $\mathbf{MQIO}(\cK,\tau,\theta)$, $\bar\cI=\{i\in\cI \,|\, \bar v_i=1\}$, $\bar\cS=\{k\in\cK \,|\, \bar u_k=1\}$, and $\bar\bc=\sum_{i\in \bar\cI}\lambda_i\ba^i$ where $\lambda_i>0$ for all $i\in \bar\cI$. Then, $\displaystyle \max_{\bx\in \cX^*(\bar\bc)}\{\|\hat\bx^k-\bx\|_\ell\}\le \underset{\bx\in {\cap}_{i\in \bar\cI} \cX_i}{\min} \{\|\hat\bx^k-\bx\|_\ell\}+\displaystyle\min_{i\in\bar\cI}\;  
 \{\rho_i\}$ for all $k\in\bar\cS$.
\end{proposition}
}

\rev{
Because $\displaystyle \max_{\bx\in \cX^*(\bar\bc)}\{\|\hat\bx^k-\bx\|_\ell\}$ is bounded above for each $k\in\bar\cS$, the forward stability measure given any distance function $d$ in \eqref{eq:for_stab}, e.g., 
$\displaystyle \max_{k\in\bar\cS}\max_{\bx\in \cX^*(\bar\bc)}\{\|\hat\bx^k-\bx\|_\ell\}$, is also bounded above.
}

\rev{
Given the same inverse-feasible cost vector $\bar\bc$ in Proposition \ref{lemma:FS_bound}, if $\ell=\infty$, the exact forward stability measure \eqref{eq:for_stab} can be written as $\displaystyle \underset{k\in\bar\cS}{\max}\; \Big\{ \displaystyle\max_{\bx\in\cX^*(\bar\bc)}\; \{\|\bx-\hat\bx^k\|_\infty\} \Big\}$. The optimal value of the inner max problem for each $k$ can be obtained by solving $2n$ LPs: $z^{(1)}_{kj} = \underset{\bx\in\cX^*(\bar\bc)}{\max}\;\{ x_j-\hat x^k_j\}$ and $ z^{(2)}_{kj}=\underset{\bx\in\cX^*(\bar\bc)}{\max}\; \{ \hat x^k_j-x_j\}$ for each $j\in\cJ=\{1,\ldots,n\}$ and selecting the largest value among them, i.e., $\displaystyle\max_{j\in \cJ}\, \{ z^{(1)}_{kj}, z^{(2)}_{kj}\}$. Thus, we have the exact forward stability measure:  $\displaystyle \underset{k\in\bar\cS}{\max}\; \Big\{ \displaystyle\max_{\bx\in\cX^*(\bar\bc)}\; \{\|\bx-\hat\bx^k\|_\infty\} \Big\}
 = \displaystyle\max_{k\in \bar\cS} \; \max_{j\in \cJ}\, \{ z^{(1)}_{kj}, z^{(2)}_{kj}\}.$
 }
 
\vspace{0.1in}

\subsection{\rev{Extensions to the MQIO model}}\label{appx:Extensions}
In Section~\ref{sec:model}, we assume the constraint parameters $(\bA,\bbb)$ remain the same over different observations. However, the proposed formulations can also apply when $\bbb$ varies over $k$ (i.e., $\bbb^k$ for observation $k$) as they only exploit the structure of $\bA$. Should $\bA$ vary over $k$ (i.e., $\bA^k$ for observation $k$), an inverse-feasible cost vector can be obtained by replacing constraint~\eqref{eq:IO_mip.2} by $b_i^k \leq {\ba^{ki}}'(\bxhat^k-\bepsilon^k) \leq b_i^k +M_1(1-v_i^k), \forall i  \in \cI,\forall k \in \cK$, constructing observation-specific conic hulls $\cC_{k}$ using $\bv^k$ similarly as in Section~\ref{sec:model}, and finding the intersection of the multiple conic hulls (e.g., via the following additional constraint: $\bc=\sum_{i\in\cI} \lambda_i^k \ba^{ki}$ and $\bzero\le \blambda^k\le \bv^k$ for all $k\in \cS$). 

\rev{
Additionally, we note that the MQIO model is also more stable than the previous inverse LP method when the parameters ($\bA$,$\bbb$) are subject to uncertainty (i.e., the realized parameters may turn out different from the assumed ones). Suppose that both the previous inverse LP and our MQIO models are used to infer cost vectors given some ``nominal'' parameters $(\bar\bA,\bar\bbb)$. Since the previous model finds a cost vector that is identical to one of the $\ba^i$'s, a slight change in the corresponding $\ba^i$ can make the inferred cost vector no longer inverse-feasible. On the other hand, 
if we let $\bar\cC$ be the set of inverse-feasible cost vectors from the original MQIO model, there exists $\bc\in\bar\cC$ that is still inverse-feasible for the \textit{realized} MQIO problem (with $(\tilde\bA,\tilde\bbb)\neq (\bar\bA,\bar\bbb)$) as long as there exists $\bx^k$ within $\tau$-distance from at least 
$\lceil \theta K \rceil$ data points such that 
$\bx^k\in\cap_{i\in\tilde\cI}\cX_i$ for some set $\tilde\cI\subseteq\cI$ 
and $\bar\cC\cap \textup{cone}(\{\ba^i\}_{i\in \tilde\cI})\neq \emptyset$. In general, if $(\bA,\bbb)$ is subject to uncertainty, the MQIO problem lends itself to traditional robust optimization or stochastic programming techniques as it is a linear MIP.
}

Finally, should the cost vector be constrained, say $\bD\bc \le \bd$, the MQIO formulation can be rewritten as follows:
\begin{equation}\label{eq:IO_mip_ConstrainedC}
\underset{\bv,\bu,\bc,\blambda,\{\bepsilon^k\}_{k\in\cK}}{\text{maximize}} \; \Big\{\,\sum_{i\in\cI} v_i \, \Big{|} \, \bD\bc\le \bd , \, \bc=\sum_{i \in \cI}\lambda_i \ba^i, \, \bzero\le\blambda \le \bv, \,  \text{\eqref{eq:IO_mip.2}--\eqref{eq:IO_mip.5}}\Big\},
\end{equation}
where the first three constraints are added to ensure that the resulting $\bc$ satisfies the cost vector constraint and is still a conic combination of some $\ba^i$'s. 

\vspace{0.1in}

\section{Supplemental Materials for the Algorithms}

\subsection{The Biclique Problem Formulation}\label{sec:clique}
Let $\bar\bD \in \{0,1\}^{ K\times m}$ be the bi-adjacency matrix defined as \eqref{eq:D_bar}. The following MIP finds an all-one submatrix of $\bar\bD$ (i.e., a biclique) with at least $\lceil \theta K \rceil$ rows (i.e., data points) and as many columns ($\ba^i$'s) as possible.
\begin{subequations} \label{eq:Biclique}
	\begin{align}
	\text{Clique($\bar\bD,\theta$)}:\quad  \underset{\bv^{clq},\bu^{clq}}{\text{maximize}} 
	&\quad \sum_{i \in \cI} v_i^{clq} & \label{eq:Biclique1}\\ 
	\text{subject to}
	&\quad  u_k^{clq} + v_i^{clq}  \le 1, \quad \forall (k,i) \in \cK\times\cI \text{ such that } \bar\bD_{ki}=0,& \label{eq:Bicliqe2}\\ 
	& \quad \sum_{k\in \cK} u_k^{clq} \ge \theta K,& \label{eq:Bicliqe3}\\ 
	&\quad  v_i^{clq},u_k^{clq} \in \{0,1\}, 
	\quad \forall i  \in \cI ,\forall k \in \cK.&\label{eq:Bicliqe4}
	\end{align}
\end{subequations}
\rev{
In \eqref{eq:Biclique}, $u_k^{clq}$ and $v_i^{clq}$ are 1 if node $k\in V_1$ and $i\in V_2$ are selected to be in the clique, respectively, and 0 otherwise. The objective function maximizes the number of nodes selected from $V_2$ (i.e., columns of the all-one submatrix of $\bar\bD$). Constraint \eqref{eq:Bicliqe2} ensures that at most one of the nodes can be selected from $k\in V_1$ and $i\in V_2$ if there is no edge corresponding to $(k,i)$, and constraint \eqref{eq:Bicliqe3} guarantees that at least $\lceil \theta K\rceil$ nodes from $V_1$ (i.e., rows of the submatrix) are selected. Once problem \eqref{eq:Biclique} is solved, the submatrix can be retrieved by selecting rows and columns of $\bar\bD$ corresponding to $\{k\in\cK \,|\, u_k^{clq}=1\}$ and $\{i\in\cI \,|\, v_i^{clq}=1\}$, respectively. 
}

\vspace{0.1in}
\subsection{Alternative Heuristic}\label{appx:algs}
The algorithms proposed in Section~\ref{sec:Algorithms} can still be computationally burdensome because the construction of each $\bar\bD$ requires solving the MIP problem \eqref{eq:IO_MIP_Reduced} $K$ times. To address this, we propose an LP relaxation of problem~\eqref{eq:IO_MIP_Reduced} using the idea of weighted $\ell_1$ minimization \citep{candes2008enhancing}. For each $k \in \cK$ we solve the following LP formulation:
\begin{subequations}\label{eq:IO_Alpha}	
	\begin{align}	
	\underset{\balpha, \bepsilon}{\text{minimize}} 
	&\quad \sum_{i\in\cI} w_i^k \alpha_i
	\nonumber\\	
	\text{subject to} 
	&\quad b_i \leq {\ba^i}'(\bxhat^k-\bepsilon) \leq b_i +M\alpha_i,
	\quad \forall i  \in \cI ,
	\label{eq:IO_Alpha2}\\	
	&\quad \lVert \bepsilon \rVert_{\ell}\leq \tau,
	\label{eq:IO_Alpha3}\\
	&\quad  \balpha \ge \bzero. 
	\label{eq:IO_Alpha4}
	\end{align}	
\end{subequations}
Then we construct a matrix $\tilde\bD \in \{0,1\}^{ K\times m}$:
\begin{equation}\label{D_tilde}
\tilde\bD_{ki} = \begin{cases} 1, 
& \mbox{if } \alpha^k_i=0, \\ 0, 
& \mbox{otherwise,} 
\end{cases}
\end{equation} 
where $\balpha^k$ is a feasible solution to \eqref{eq:IO_Alpha} with respect to $\bxhat^k$. Note that the construction of $\tilde\bD$ is more efficient than that of $\bar\bD$. 
\rev{
The following result shows the relationship between a solution to problem~\eqref{eq:IO_Alpha}	and a solution to problem~\eqref{eq:IO_MIP_Reduced}.
}
\rev{
\begin{lemma}\label{lemma:Alpha_vs_ReducedMIP}
Given an observation $\hat\bx^k$, let $(\bar\balpha^k,\bar\bepsilon^k)$ be a feasible solution to problem~\eqref{eq:IO_Alpha}  and $\bar\cI = \{i\in\cI\,|\,\bar\alpha^k_i=0\}$. Then, there exists a feasible solution $(\bar\bv^k,\bar\bepsilon^k)$ to problem~\eqref{eq:IO_MIP_Reduced} with respect to $\hat\bx^k$ where $\bar v^k_i=1$ for $i\in \bar\cI$ and $\bar v^k_i=0$ for $i\in \cI\setminus\bar\cI$.
\end{lemma} 
Next, similar to Proposition~\ref{prop:Dbar-IO}, the following result shows how a feasible solution to the QIO (hence MQIO) model can be retrieved from $\tilde\bD$.
\begin{proposition}\label{prop:Dtilde-IO}
If there exists $\tilde\bD$ that satisfies \eqref{D_tilde} and has an all-one submatrix whose rows and columns correspond to $\bar\cS$ and $\bar\cA$, respectively, where $\bar\cS \subseteq \cK$, $|\bar \cS| \ge \theta K$, and $\bar\cA \subseteq \cI$, then there exists a solution $(\bar{\bc},\{\bar{\bepsilon}^k\}_{k\in \cK},\bar\by,\bar{\cS})$ feasible for $\mathbf{QIO}(\cK,\tau,\theta)$ such that $\bar\bc \in \textup{cone}(\{\ba^i\}_{ i \in \bar\cA})$.
\end{proposition}
}

\rev{
We use the same instance provided in Example \ref{ex:InvSol_Clique} to illustrate the result of Proposition \ref{prop:Dtilde-IO}. 
\begin{example}\label{ex:Dtilde}
Let $(\bar\balpha^k,\bar\bepsilon^k)$ denote the optimal solution to problem \eqref{eq:IO_Alpha} associated with each data point $k=1,\ldots,5$: we have $\bar\balpha^k = [0,\,0,\,\bar\alpha_3^k>0,\,\bar\alpha_4^k>0]'$ for $k=1,2,3,4,$ (with $\bar\bepsilon^1=\begin{bsmallmatrix}-0.5\\-0.2\end{bsmallmatrix}$, $\bar\bepsilon^2=\begin{bsmallmatrix}-0.3\\-0.2\end{bsmallmatrix}$, $\bar\bepsilon^3=\begin{bsmallmatrix}-0.3\\-0.5\end{bsmallmatrix}$, and $\bar\bepsilon^4=\begin{bsmallmatrix}-0.5\\-0.5\end{bsmallmatrix}$), and $\bar\balpha^5=[\bar\alpha_1^5>0,\,0,\,0,\,\bar\alpha_4^5>0]'$ (with $\bar\bepsilon^5=\begin{bsmallmatrix}-0.3\\0.3\end{bsmallmatrix}$). 
We then construct a bipartite graph by defining node $i$ for $\ba^i$, node $k$ for $\hat\bx^k$, and edge $(k,i)$ where $\bar\alpha_i^k=0$, which leads to the same bipartite graph in Figure~\ref{fig:Prop5-Bipartite}. As a result, the corresponding bi-adjacency matrix $\tilde\bD$ is identical to the $\bar\bD$ matrix in Figure \ref{fig:Prop5-Matrix}. Similar to Example \ref{ex:InvSol_Clique}, we can obtain an inverse-feasible cost vector from the all-one submatrix of $\tilde\bD$.
\end{example} 
}

Similar to the algorithms in Section~\ref{sec:Algorithms}, there can be multiple $\tilde\bD$'s because there may be multiple optimal solutions to problem~\eqref{eq:IO_Alpha}. We propose a heuristic, which we call $\tilde\bD$-Alg-Heuristic, via a similar weighting approach to find a matrix $\tilde\bD$ and the same clique problem to find the all-one submatrix of $\tilde\bD$ (i.e., Clique($\tilde\bD,\theta$)). The pseudo-code for this heuristic can be found in Algorithm \ref{alg:Dtilde-alg}.

\begin{algorithm}[H]\label{alg:Dtilde-alg}
	\caption{$\tilde\bD$-Alg-Heuristic}
	\SetAlgoLined
	\KwResult{$\cC^*$}
	\KwIn{$\bA, \bbb, \hat\cX, \tau,\theta$}
	$ \bw^1 \leftarrow \beee, \ \tilde\bD \leftarrow [\ ] $ \\
	\For{$k=1,\ldots,K$}{ 
		Find $\balpha^{k}\in \underset{\balpha}{\text{argmin}
		} \{{\bw^k}'\balpha \,|\, \text{\eqref{eq:IO_Alpha2}--\eqref{eq:IO_Alpha4}} \}$ with $\hat\bx^{k}$ as input data\\
		$\tilde\bD_{ki}\leftarrow 1$ if $\alpha^{k}_i=0$; $\tilde\bD_{ki}\leftarrow 0$ otherwise \quad // $\tilde\bD_{ki}\coloneqq (k,i)$ entry of $\tilde\bD$\\
		$\vartheta^{k}_i \leftarrow \sum_{\ell=1}^{k} \tilde\bD_{\ell i}$ for each $i\in \cI$\\ 
		Update weights $\bw^{k+1}\leftarrow \bvartheta^{k}/\|\bvartheta^{k}\|$
	}
		Solve Clique($\tilde\bD,\theta)$ to find the all-one submatrix of $\tilde\bD$. Define $\bv^*\in\{0,1\}^{m}$ where 1's correspond to columns of the submatrix.\\
		Construct the set of inverse-feasible cost vectors $\cC^* = \textup{cone}(\{\ba^i\}_{i: v_{i}^*=1})$\\
\end{algorithm}



\vspace{0.1in}

\subsection{Application to Online Learning}\label{sec:model_online}
The assumption we made in Section \ref{sec:model} is that the data points are available all at once in advance from a fixed forward LP formulation. However, in many application domains, data may become available through different time points in separate batches. Moreover, some parameters of the underlying forward LP formulation may change at each time point as well. In this subsection, we show how our inverse method can be extended for learning the cost vectors adaptively over time in an online manner.

Suppose at each time $t \in \cT = \{1,\ldots,T\}$, the DM observes an external signal $\bbb=\bbb^t$ as the right-hand-side of the forward problem and makes a set of decisions $\hat\cX_t$ with the index set $\cK_t$ of cardinality $K_t$. We assume that the decisions can be noisy. The following formulation is a modification to MQIO that finds inverse-feasible cost vectors for the entire collection of batches $\bigcup_{t\in\cT}\hat\cX_t$: 
\begin{subequations}\label{eq:IO_mipOnline}   	
	\begin{align}   	
	\underset{\bv,\bu,\{\bepsilon^{kt}\}_{k\in\cK_t,t\in\cT}}{\text{maximize}}	
	&\quad \sum_{i\in \cI} v_i \label{eq:IO_mipOnline.1}\\   	
	\text{subject to}	
	&\quad \bbb^t \leq \bA(\bxhat^{kt}-\bepsilon^{kt}) \leq \bbb^t +M_1(\mathbf{1}-\bv), \quad   \forall k \in  \cK_t, \; \forall t\in\cT, \label{eq:IO_mipOnline.2}\\	
	&\quad \lVert\bepsilon^{kt}\rVert_{\ell}\leq \tau_t+M_2(1-u_{kt}), \quad \forall k \in \cK_t,  \; \forall t\in\cT, \label{eq:IO_mipOnline.3}\\	
	&\quad \sum_{t\in\cT}\sum_{k\in\cK_t} u_{kt} \ge \theta \sum_{t\in\cT}K_t, \label{eq:IO_mipOnline.4}\\	
	&\quad  v_i,u_{kt} \in \{0,1\},	
	\quad \forall i  \in \cI,\; \forall k \in \cK_t, \; \forall t\in\cT. \label{eq:IO_mipOnline.5}	
	\end{align}   	
\end{subequations}
The above formulation still aims to solve the inverse problem in the traditional setting, which means it has to wait until all data points are collected, leading to a large-size MIP. Instead, we propose an online learning approach using a modified version of $\bar\bD$-Alg-Heuristic in Algorithm~\ref{alg:Online1}. 

In this algorithm, we start with an empty $\bar\bD$ matrix and update it at each time $t$ by finding $\bv^*\in\underset{\bv}{\text{argmax}}\{{\bw^t}'\bv\,|\,\eqref{eq:IO_mip.2}$--$\eqref{eq:IO_mip.5} \}$ with $\hat\cX_t$ as input data (where $\bw^t$ is an appropriate weight vector at time $t$) and inserting ${\bv^*}'$ into the next empty row of $\bar\bD$. At each time $t$, given the updated $\bar\bD$ we find its largest all-one submatrix and let $\bv^t$ be a binary vector with 1's corresponding to the columns of the submatrix. Recall that each column of this matrix corresponds to each $\ba^i$; thus, $\bv^t$ indicates which $\ba^i$'s should be used for creating the set of inverse-feasible cost vectors at time $t$, i.e., $\cC_t = \textup{cone}(\{\ba^i\}_{i: v_i^t=1})$. We show in the numerical results section that with a sufficiently large $T$, there exists $\bc\in\cC_T$ that is also inverse-feasible for the above MIP \eqref{eq:IO_mipOnline}. Again, with such a large $T$ and an excessive number of data points, the MIP problem \eqref{eq:IO_mipOnline} would have been computationally extremely challenging, whereas the online learning approach  is more efficient as it attempts to solve the problem in a distributed manner.

\begin{algorithm}[H]\label{alg:Online1}
	\caption{Online algorithm to find inverse-feasible solutions}
	\SetAlgoLined
	\KwResult{$\cC_T$}
	\KwIn{$\bA, \tau,\theta$}
	$ \bw^1 \leftarrow \textbf{1}, \ \bar\bD \leftarrow [\ ] $ \\
	\For{$t=1,\ldots,T$}{ 
		Observe $(\bbb^t,\hat\cX_t)$\\
		Find $\bv^{*}\in \underset{\bv}{\text{argmax}} \{{\bw^t}'\bv \,|\, \text{\eqref{eq:IO_mip.2}--\eqref{eq:IO_mip.5}} \}$ with $\hat\cX_t$ as input data and $\bbb=\bbb^t$\\
		Insert ${\bv^{*}}'$ to the $t$-th row of $\bar\bD$ \\
		Solve Clique($\bar\bD,\theta)$ to find the all-one submatrix of $\bar\bD$. Define $\bv^t\in\{0,1\}^{m}$ where 1's correspond to columns of the submatrix.\\
		Construct the set of inverse-feasible cost vectors $\cC_t = \textup{cone}(\{\ba^i\}_{i: v_i^t=1})$\\
		Update weights $\bw^{t+1} \leftarrow (\textbf{1}+\bv^t)/\lVert \textbf{1}+ \bv^t \rVert$ 
	}
\end{algorithm}

\vspace{0.2in}


\section{Proofs} \label{appx:proofs}
\rev{
\proof{Proof of Proposition \ref{prop:ISLB_Nondecreasing}:}
Given $\hat\cX$, $\tau$, $\theta_1$, and $\theta_2$, let $\hat\cC(\theta_1)$, $\hat\cC(\theta_2)$, $\tilde\cC(\theta_1)$, and $\tilde\cC(\theta_2)$ denote the set of cost vectors obtained by $\mathbf{QIO}(\hat\cX,\tau,\theta_1)$, $\mathbf{QIO}(\hat\cX,\tau,\theta_2)$, $\mathbf{QIO}(\tilde\cX,\tau,\theta_1)$, and $\mathbf{QIO}(\tilde\cX,\tau,\theta_2)$, respectively, where $\tilde\cX$ denotes some shifted dataset. To show that inverse stability measure \eqref{eq:inv_stab} for $\hat\cC(\theta)$ is non-decreasing as $\theta$ decreases, we need to show that given $\theta_1$ and $\theta_2$, if $\theta_2\le \theta_1$, then $\underset{\tilde\cX\subset\mathbb{R}^{K\times n}}{\min}\{d(\hat\cX,\tilde\cX)|\ \hat\cC(\theta_2)\cap \tilde\cC(\theta_2)=\emptyset\}\ge \underset{\tilde\cX\subset\mathbb{R}^{K\times n}}{\min}\{d(\hat\cX,\tilde\cX)|\ \hat\cC(\theta_1)\cap \tilde\cC(\theta_1)=\emptyset\}$.
\
Given $\theta_2\le \theta_1$
, we have $\hat\cC(\theta_1)\subseteq \hat\cC(\theta_2)$ and $\tilde\cC(\theta_1)\subseteq \tilde\cC(\theta_2)$. Let $\hat\cC_{\text{diff}}=\hat\cC(\theta_2)\setminus\hat\cC(\theta_1)$ and $\tilde\cC_{\text{diff}}=\tilde\cC(\theta_2)\setminus\tilde\cC(\theta_1)$. Then,
%
%
\begin{subequations}\nonumber
    \begin{align}
      \underset{\tilde\cX\subset\mathbb{R}^{K\times n}}{\min}\{d(\hat\cX,\tilde\cX)|\ \hat\cC(\theta_2)\cap \tilde\cC(\theta_2)=\emptyset\}
      &=\underset{\tilde\cX\subset\mathbb{R}^{K\times n}}{\min}\{d(\hat\cX,\tilde\cX)|\ \big(\hat\cC_{\text{diff}}\cup\hat\cC(\theta_1)\big)\cap \big(\tilde\cC_{\text{diff}}\cup\tilde\cC(\theta_1)\big)=\emptyset\}\\
    &=\underset{\tilde\cX\subset\mathbb{R}^{K\times n}}{\min}\{d(\hat\cX,\tilde\cX)|\ \big(\tilde\cC_{\text{diff}}\cup\tilde\cC(\theta_1)\big)\cap \big(\hat\cC_{\text{diff}}\cup\hat\cC(\theta_1)\big)=\emptyset\}\\
      &= \underset{\tilde\cX\subset\mathbb{R}^{K\times n}}{\min}\{d(\hat\cX,\tilde\cX)|\ \big(\tilde\cC_{\text{diff}}\cap(\hat\cC_{\text{diff}}\cup\hat\cC(\theta_1))\big)\cup\\
      & \hspace{1.8in} \big(\tilde\cC(\theta_1)\cap (\hat\cC_{\text{diff}}\cup\hat\cC(\theta_1))\big)=\emptyset \}\\
      &\ge \underset{\tilde\cX\subset\mathbb{R}^{K\times n}}{\min}\{d(\hat\cX,\tilde\cX)|\ \tilde\cC(\theta_1)\cap (\hat\cC_{\text{diff}}\cup\hat\cC(\theta_1))=\emptyset \}\\
      &= \underset{\tilde\cX\subset\mathbb{R}^{K\times n}}{\min}\{d(\hat\cX,\tilde\cX)|\ \big(\tilde\cC(\theta_1)\cap \hat\cC_{\text{diff}}\big)\cup \big(\tilde\cC(\theta_1) \cap \hat\cC(\theta_1)\big)=\emptyset \}\\
      &\ge \underset{\tilde\cX\subset\mathbb{R}^{K\times n}}{\min}\{d(\hat\cX,\tilde\cX)|\ \tilde\cC(\theta_1) \cap \hat\cC(\theta_1)=\emptyset \}, 
    \end{align}
\end{subequations}
which completes the proof.~\Halmos
\endproof
}

\vspace{0.1in}

\proof{Proof of Proposition \ref{prop:IO1}:}
Let $(\bar{\bc},\{\bar{\bepsilon}^k\}_{k\in \cK},\bar{\by},\bar{\cS})$ be a feasible solution to $\textrm{\textbf{QIO}}(\cK,\tau,\theta)$. Note that $\bar{y}_i >0$ for at least one $i \in \cI$. Let $\bar{\cI}=\{ i \in \cI \  \arrowvert \ \bar{y}_{i}>0
\}$. From \eqref{eq:IO1.2}, $\bar{\bc}=\sum_{i\in\bar{\cI}}\bar{y}_i\ba^i$. From \eqref{eq:IO1.5}, for each $k \in \bar{\cS}$, $\sum_{i\in\bar{\cI}}\bar{y}_i{\ba^i}'(\bxhat^k-\bar\bepsilon^k)=\sum_{i\in\bar{\cI}}b_i\bar{y}_i$ and  thus $\sum_{i\in\bar{\cI}}\bar{y}_i({\ba^i}'(\bxhat^k-\bar\bepsilon^k)-b_i)=0$. 
Because $\bar{y}_i >0$, we have ${\ba^i}'(\bxhat^k-\bar\bepsilon^k)-b_i=0, \forall i \in \bar{\cI}$. Now, pick any arbitrary $\tilde{i} \in \bar{\cI}$ and let $\tilde{\by}=\beee_{\tilde{i}}$ where $\beee_{\tilde{i}}$ denotes the $\tilde{i}$-th unit vector. Then, from \eqref{eq:IO1.2}, we 
can construct a new cost vector $\tilde{\bc}=\ba^{\tilde{i}}$. By replacing $\bar\bc$ with this $\tilde{\bc}$ in \eqref{eq:IO1.5}, \eqref{eq:IO1.5} becomes ${\ba^{\tilde{i}}}'(\bxhat^k-\bar\bepsilon^k)-b_{\tilde{i}}=0, \forall k\in\bar\cS$, which still holds because $\tilde i \in \bar\cI$. It is clear that $(\tilde{\bc}=\ba^{\tilde{i}},\{\bar{\bepsilon}^k\}_{k\in \cK},\tilde{\by}=\beee_{\tilde{i}},\bar{\cS})$ also satisfies constraints \eqref{eq:IO1.2}--\eqref{eq:IO1.9} and thus is feasible for \eqref{eq:IO1}. \Halmos
\endproof 
\vspace{0.1in}
\rev{
\proof{Proof of Lemma \ref{lemma:Feasible_conic}:}
Given some sets $\bar\cA\subseteq\cA$ and $\bar\cS$ such that $|\bar\cS|\ge \theta K$, let $\bar\bc\in \textup{cone}(\{\ba^i\}_{i\in \bar\cA})$ such that $\|\bar\bc\|_p=1$. Recall that $\cX_i=\{\bx\in\cX\,|\,{\ba^i}'\bx = b_i\},i\in\cI$. We want to show that if there exists $\bar\bx^k\in \cX_i$ for all $i\in \bar\cA$ {such that $\|\bar\bx^k-\hat\bx^k\|_{\ell} \le \tau$ for all $k\in\bar\cS$}, then there exists $\{\bar\bepsilon^k\}_{k\in \cK}$ and $\bar\by$ such that $(\bar{\bc},\{\bar{\bepsilon}^k\}_{k\in \cK},\bar{\by},\bar\cS)$ is feasible for \textbf{QIO}$(\cK,\tau,\theta)$ {(i.e., formulation~\eqref{eq:IO1})}.  Note that $\bar\bc$ can be written as $\bar\bc=\sum_{i\in \bar\cA}\bar\lambda_i \ba^i$ {where $\bar\lambda_i>0$ for at least one $i\in \bar\cA$}. Let $\bar\lambda_i=0$ for $i\in \cI\setminus \bar\cA$. Because $\bar\bx^k\in\cX^i,\forall i\in\bar\cA$, for each $k\in \bar\cS$ we have ${\ba^i}'\bar\bx^k=b_i$ for all $i\in \bar\cA$. 
Multiplying both sides by $\bar\lambda_i$ and summing over all $i \in \bar\cA$, we have $\sum_{i\in \bar\cA}\bar\lambda_i {\ba^i}'\bar\bx^k=\sum_{i\in \bar\cA}\bar\lambda_i b_i$, 
and because $\bar\lambda_i=0$ for $i\in \cI\setminus \bar\cA$ we further have $\sum_{i\in \cI}\bar\lambda_i {\ba^i}'\bar\bx^k=\sum_{i\in \cI}\bar\lambda_i b_i$, hence $\bar\bc'\bar\bx^k=\bbb'\bar\blambda$ for all $k\in\bar\cS$. By letting $\bar{\by}=\bar\blambda$, we have $\bA'\bar\by = \bar\bc$ and $\bar\by\ge\bzero$, which satisfy constraints~\eqref{eq:IO1.2} and \eqref{eq:IO1.3}, respectively. For $k\in \bar\cS$, let $\bar\bepsilon^k=\hat\bx^k-\bar\bx^k$; for $k\in \cK\setminus\bar\cS$, let $\bar\bepsilon^k$ be an arbitrary vector in $\mathbb{R}^n$. Then constraint~\eqref{eq:IO1.4} is satisfied because $\hat\bx^k-\bar\bepsilon^k=\bar\bx^k\in\cX,\forall k\in\bar\cS$, and constraint~\eqref{eq:IO1.5} is satisfied because $\bar\bc'(\hat\bx^k - \bar\bepsilon^k) = \bar\bc'\bar\bx^k=\bbb'\bar\blambda =\bbb'\bar\by,\forall k\in\bar\cS$. Constraint~\eqref{eq:IO1.6} is satisfied because $\|\bar\bepsilon^k\|_\ell = \|\hat\bx^k-\bar\bx^k\|_\ell\le \tau, \forall k\in \bar\cS$. Finally, \eqref{eq:IO1.7}, \eqref{eq:IO1.8}, and \eqref{eq:IO1.9} are satisfied as it was assumed that $\bar\cS\subseteq\cK$, $|\bar\cS|\ge \theta K$ and $\|\bar\bc\|_p=1$. As a result, $(\bar{\bc},\{\bar{\bepsilon}^k\}_{k\in \cK},\bar{\by},\bar\cS)$ is feasible for \textbf{QIO}$(\cK,\tau,\theta)$.\Halmos
\endproof
}
\
\proof{Proof of Theorem \ref{prop:IO_mip}:}
$(\Rightarrow)$ Assume $(\bar{\bc} \ne \mathbf{0},\{\bar{\bepsilon}^k\}_{k\in \cK},\bar{\by},\bar\cS)$ is feasible for $\mathbf{QIO}(\cK,\tau,\theta)$, i.e., model~\eqref{eq:IO1}. Then from \eqref{eq:IO1.2}, $\bar{\bc}=\sum_{i\in\bar{\cI}}\bar{y}_i\ba^i$ where $\bar\cI=\{i\in \cI \,|\, \bar{y}_i>0\}$. Construct $\bar\bv\in\mathbb{R}^{|\cI|}$ by setting $\bar v_i=1$ for all $i \in \bar \cI$ and $\bar v_i=0$ for all $i \in \cI \setminus \bar\cI$. Also define $\bar\bu\in\mathbb{R}^{|\cK|}$ such that $\bar u_k=1$ for $k \in \bar\cS$ and $\bar u_k=0$ for $k \in \cS\setminus\bar\cS$. Let $\bepsilon^k=\bar\bepsilon^k$ for $k\in \bar\cS$ and $\bepsilon^k=\hat\bx^k-\bx$ for $k\in \cK\setminus \bar\cS$ where $\bx\in\cX$ satisfying ${\ba^i}'\bx=b_i, \forall i \in \bar\cI$. We want to show that $(\bar\bv,\bar\bu,\{{\bepsilon}^k\}_{k\in \cK})$ is feasible for model \eqref{eq:IO_mip}. From \eqref{eq:IO1.2} and \eqref{eq:IO1.5}, we have $\sum_{i\in\bar{\cI}}\bar{y}_i{\ba^i}'(\hat\bx^k-\bepsilon^k)=\sum_{i\in\bar{\cI}}b_i\bar{y}_i$, or equivalently $\sum_{i\in\bar{\cI}}\bar{y}_i({\ba^i}'(\hat\bx^k-\bepsilon^k)-b_i)=0$ for each $k\in \bar\cS$. Since $\bar y_i>0$ for all $i\in \bar\cI$, for each $k\in \bar\cS$ we have ${\ba^i}'(\hat\bx^k-\bepsilon^k) = b_i$ for $i\in \bar\cI$ (i.e., $i:\bar{v}_i=1$) and ${\ba^i}'(\hat\bx^k-\bepsilon^k) \ge b_i$ for $i\in \cI\setminus\bar\cI$ (i.e., $i:\bar{v}_i=0$). Also, for each $k \in \cK\setminus \bar\cS$, based on how $\bepsilon^k$ is defined above, we also have ${\ba^i}'(\hat\bx^k-\bepsilon^k) = b_i$ for all $i\in \bar\cI$ and ${\ba^i}'(\hat\bx^k-\bepsilon^k)\ge b_i$ for all $i\in \cI$. Therefore, $\{{\bepsilon}^k\}_{k\in \cK}$ and $\bar\bv$ together satisfy \eqref{eq:IO_mip.2}. Because $\bar\bu$ is constructed based on whether or not $\|\bar\bepsilon^k\|_{\ell}\leq \tau$, $\{{\bepsilon}^k\}_{k\in \cK}$ and $\bar\bu$ clearly satisfy \eqref{eq:IO_mip.3}. Also, $\sum_{k=1}^{ K} \bar{u}_k = \sum_{k\in\bar\cS} \bar{u}_k = |\bar\cS| \ge \theta K$, satisfying \eqref{eq:IO_mip.4}. Thus, $(\bar\bv,\bar\bu,\{{\bepsilon}^k\}_{k\in \cK})$ is feasible for model \eqref{eq:IO_mip}. Finally, from \eqref{eq:IO1.2}, $\bar\bc=\bA'\bar\by=\sum_{i\in \bar\cI}\bar y_i \ba^i=\sum_{i:\bar{v}_i=1} \bar y_i \ba^i$, and therefore $\bar\bc \in \textup{cone}(\{\ba^i\}_{ i:  \bar{v}_i=1})$. Because $\bar\bc$ is feasible for $\mathbf{QIO}(\cK,\tau,\theta)$, $\lVert \bar\bc\rVert_p=1$.

\noindent$(\Leftarrow)$ Assume $(\bar\bv \ne \mathbf{0},\bar\bu,\{\bar{\bepsilon}^k\}_{k\in \cK})$ is feasible for model \eqref{eq:IO_mip}. Let $\bar\cI =\{i\in \cI\ |\ \bar{v}_i=1\}$ and $\bar\cS=\{k\in\cK\ | \ \bar u_k=1\}$, and consider $\bar\bc\in \textup{cone}(\{\ba^i\}_{i \in \bar\cI})$ such that $\|\bar\bc\|_p=1$. We can write $\bar{\bc}=\sum_{i\in \bar\cI}\bar{\lambda}_i\ba^i$ where $\bar{\lambda}_i \ge 0$ $\forall i\in\bar\cI$. If we let $\bar{y}_i=\bar{\lambda}_i$ for all $i \in \bar\cI$ and $\bar{y}_i=0$ for all $i \in \cI \setminus \bar\cI$ then $\bar{\bc}=\sum_{i\in \bar\cI}\bar{y}_i\ba^i=\bA'\bar\by$ satisfies \eqref{eq:IO1.2}; because $\bar\by \ge \bzero$,  \eqref{eq:IO1.3} is also satisfied. Because ${\ba^i}'(\hat\bx^k-\bar\bepsilon^k)\ge b_i$ holds for all $i\in\cI,k\in\cK$ by \eqref{eq:IO_mip.2},  \eqref{eq:IO1.4} is satisfied for all $k\in\cK$. 
Also because model \eqref{eq:IO_mip} is feasible, ${\ba^i}'(\bxhat^k-\bar{\bepsilon}^k)=b_i$ for all $i \in \bar\cI$ and all $k$ for which $\bar u_k=1$, i.e., $k\in\bar\cS$; multiplying both sides by $\bar{y}_i$ and summing over all $i \in \cI$, we have $\sum_{i \in \cI}\bar{y}_i{\ba^i}'(\bxhat^k-\bar{\bepsilon}^k)=\sum_{i \in \cI}\bar{y}_i b_i$ or equivalently $\bar\bc'(\bxhat^k-\bar{\bepsilon}^k)=\bbb'\bar{\by}$, satisfying \eqref{eq:IO1.5}. Clearly, $\{\bar\bepsilon_k\}_{k\in\cK}$ and $\bar\cS$ satisfy constraints \eqref{eq:IO1.6}--\eqref{eq:IO1.8}. Therefore, $(\bar\bc,\{\bar{\bepsilon}^k\}_{k\in \cK},\bar{\by},\bar{\cS})$ is feasible for $\mathbf{QIO}(\cK,\tau,\theta)$.~\Halmos
\endproof

\vspace{0.1in}

\rev{
\proof{Proof of Proposition~\ref{prop:IS_LB}:}
For this proof, we first refer the reader to Appendix~\ref{appx:LB} where we show a formulation to compute the true inverse stability measure; this formulation is non-convex but we use it for proving Proposition~\ref{prop:IS_LB}. Note that the true inverse stability measure $\xi^*$ is the optimal value of problem~\eqref{eq:Exact_IS} in Appendix~\ref{appx:LB}.
}

\rev{
We first present a variant of \eqref{eq:Exact_IS} for each $i\in\bar\cI$ separately (i.e., $|\bar\cI|$ of them) and show that the optimal objective value of \eqref{eq:Exact_IS}, i.e., $\xi^*$, is bounded below by the largest optimal objective value among these $|\bar\cI|$ smaller problems. 
Consider the following problem for each $i\in \bar\cI$: 
\begin{subequations}\label{eq:Exact_IS_Reform_k}	
\begin{align}
\underset{\{\bx^{k}\}_{k\in \cK},\cS^i}{\text{minimize}} &\quad  \sum_{k\in \cK}\|\bx^k-\hat\bx^k \|_\ell \label{eq:Exact_IS_Reform_k_1}\\
\text{subject to}
&\quad D^{i}(\bx^{k})=\underset{\bz\in\cX_i}{\min} \{\|\bx^k-\bz\|_{\ell}\}, \quad \; \forall  k\in \cK,\label{eq:Exact_IS_Reform_k_2}\\
& \quad D^{i}(\bx^{k}) \ge \tau+\sigma,  \quad \; \forall k\in \cS^i,\label{eq:Exact_IS_Reform_k_3}\\& \quad |\cS^i|\ge \lfloor (1-\theta)K\rfloor +1. \quad \; \label{eq:Exact_IS_Reform_k_4}
\end{align}
\end{subequations}
Let $\zeta_i^{*}$ be the optimal value of the above problem for 
each $i\in\bar\cI$. Note that unlike in formulation \eqref{eq:Exact_IS} where constraints of $D^i(\bx^k)\ge \tau+\sigma$ in \eqref{eq:Exact_IS_3} must be satisfied for all $i\in\hat\cI$, problem~\eqref{eq:Exact_IS_Reform_k} requires constraint~\eqref{eq:Exact_IS_Reform_k_3} be satisfied for only one $i\in\bar\cI\subseteq\hat\cI$. Hence, we have $\zeta^*_i \le \xi^{*}$ for all $i\in \bar\cI$, or equivalently, $\displaystyle \max_{i\in \bar\cI}\{\zeta^*_i\} \le \xi^{*} $. 
Problem~\eqref{eq:Exact_IS_Reform_k} is still non-convex and thus finding $\zeta_i^*$ is still a challenge. We next show that its optimal objective value $\zeta_i^{*}$ is equal to the sum of the optimal objective values of even smaller problems that can be efficiently solved.
Given $i\in\bar\cI$ and $k\in\cK$, consider the following problem:
\begin{subequations}\label{eq:Exact_IS_Reform_i_k}	
\begin{align}	
\underset{\bx^{k}}{\text{minimize}} &\quad  \|\bx^k-\hat\bx^k \|_\ell \label{eq:Exact_IS_Reform_i_k_1}\\	
\text{subject to} 
&\quad D^{i}(\bx^{k})=\underset{\bz\in\cX_i}{\min} \{\|\bx^k-\bz\|_\ell 
\}, \label{eq:Exact_IS_Reform_i_k_2}\\
& \quad D^{i}(\bx^{k})\ge \tau+\sigma.\label{eq:Exact_IS_Reform_i_k_3}
\end{align}	
\end{subequations} 
Let $\tilde\zeta_{ik}$ denote the optimal value of the above problem. For a given $i$, let $\tilde\zeta_{i[k]}$ be the $k$-th smallest value of $\tilde\zeta_{ik},k\in\cK$. Then we have $\displaystyle \zeta_i^{*} = \sum_{k=1}^{\lfloor(1-\theta)K\rfloor+1}  \tilde\zeta_{i[k]}$. Thus, $\displaystyle \max_{i\in \bar\cI}\bigg\{{\sum_{k=1}^{\lfloor(1-\theta)K\rfloor+1}  \tilde\zeta_{i[k]}}\bigg\} = \max_{i\in \bar\cI}\{\zeta_i^{*}\} \le \xi^*.$
}

\rev{
Next, we show how we can find the optimal value for formulation \eqref{eq:Exact_IS_Reform_i_k}, i.e., $\tilde\zeta_{ik}$. Let $d_{ik}^* = \underset{\bx\in\cX_i}{\min}\{\|\hat\bx^k-\bx\|_\ell\}$ and $\bar \tau=\tau+\sigma$. Note that when $d_{ik}^*\ge \bar \tau$, we have $\tilde\zeta_{ik}=0$ because $\bx^k=\hat\bx^k$ satisfies both \eqref{eq:Exact_IS_Reform_i_k_2} and \eqref{eq:Exact_IS_Reform_i_k_3}. Now we want to show that if $d_{ik}^*< \bar \tau$, $\tilde\zeta_{ik}=\bar \tau-d_{ik}^*$. To do so, we first show that $\tilde\zeta_{ik}\ge \bar \tau-d^*_{ik}$. Let $\bar\bx\in\underset{\bx\in\cX_i}{\text{argmin}}\{\|\hat\bx^k-\bx\|_\ell
\}$. Suppose to the contrary that there exists an optimal solution $\tilde\bx^k$ that attains $\tilde\zeta_{ik}<\bar \tau-d^*_{ik}$.  By the definition of $d^*_{ik}$ we can rewrite this inequality as $\|\tilde\bx^k-\hat\bx^k\|_\ell<\bar \tau-\|\hat\bx^k-\bar\bx\|_\ell$ or equivalently $\|\tilde\bx^k-\hat\bx^k\|_\ell+\|\hat\bx^k-\bar\bx\|_\ell<\bar \tau$. Furthermore, by triangular inequality the following is true: $\|\tilde\bx^k-\bar\bx\|_\ell \le \|\tilde\bx^k-\hat\bx^k\|_\ell+\|\hat\bx^k-\bar\bx\|_\ell<\bar \tau$. However, since $\bar\bx\in \cX_i$, we have $D^i(\tilde\bx^k)=\underset{\bz\in\cX_i}{\min} \{\|\tilde\bx^k-\bz\|_\ell\}\le \|\tilde\bx^k-\bar\bx\|_\ell <\bar \tau$, which violates \eqref{eq:Exact_IS_Reform_i_k_3} and thus contradicts the feasibility of $\tilde\bx^k$. Thus, it must be that $\tilde\zeta_{ik}\ge \bar \tau-d^*_{ik}$. All we need now is to show that there exits a feasible solution $\tilde\bx^k$ to problem~\eqref{eq:Exact_IS_Reform_i_k} that attains the objective function value of $\bar \tau-d^*_{ik}$. Note that if a given $\tilde\bx^k$ is feasible for problem~\eqref{eq:Exact_IS_Reform_i_k}, we must have $D^{i}(\tilde\bx^{k}) =\underset{\bz\in\cX_i}{\min} \left\{\|\tilde\bx^k-\bz\|_\ell\right\}$ bounded below by $\bar\tau$ (due to \eqref{eq:Exact_IS_Reform_i_k_3}). Consider $\tilde\bx^k = \frac{\bar \tau}{d^*_{ik}}\hat\bx^k - \frac{\bar \tau}{d^*_{ik}}\bar\bx + \bar\bx$. Clearly, $\bx^k=\tilde\bx^k$  and $\bz=\bar\bx$ together satisfy all constraints in \eqref{eq:Exact_IS_Reform_i_k} with equality because we have $\bar\bx\in \cX_i$ and 
\begin{subequations}\label{eq:Dik}
\begin{align}
\|\tilde\bx^k-\bar\bx\|_\ell&= \|\frac{\bar \tau}{d^*_{ik}}\hat\bx^k- \frac{\bar \tau}{d^*_{ik}}\bar\bx + \bar\bx-\bar\bx\|_\ell \label{eq:Dik_b}\\ 
& = \frac{\bar \tau}{d^*_{ik}} \|\hat\bx^k- \bar\bx\|_\ell \label{eq:Dik_c}\\
&=\frac{\bar \tau}{d^*_{ik}} d^*_{ik}= \bar \tau.
\end{align}	
\end{subequations} 
Moreover, $\bx^k=\tilde\bx^k$ leads to the objective value  $\tilde\zeta_{ik}=\|\hat\bx^k-\tilde\bx^k\|_\ell=\|\hat\bx^k+ \frac{\bar \tau}{d^*_{ik}}\hat\bx^k - \frac{\bar \tau}{d^*_{ik}}\bar\bx + \bar\bx\|_\ell=\frac{\bar \tau-d^*_{ik}}{d^*_{ik}}\|\hat\bx^k-\bar\bx\|_\ell=\bar \tau-d^*_{ik}$ and thus is an optimal solution to \eqref{eq:Exact_IS_Reform_i_k}. 
}

\rev{
As a result, $\displaystyle \max_{i\in\bar\cI}\bigg\{ \sum_{k=1}^{(\lfloor(1-\theta)K\rfloor+1)} \!\!\!\!\!  \max(0,\bar \tau-d^*_{i[k]})\bigg\} = \max_{i\in\bar\cI}\bigg\{ \sum_{k=1}^{(\lfloor(1-\theta)K\rfloor+1)} \!\!\!\!\! \tilde\zeta_{i[k]}\bigg\} = \max_{i\in\bar\cI}\{ \zeta_i^*\} \le \xi^*$. Clearly, because $\tau<\bar \tau$ we have $\displaystyle \max_{i\in\bar\cI}\bigg\{ \sum_{k=1}^{(\lfloor(1-\theta)K\rfloor+1)} \!\!\!\!\! \max(0,\tau-d^*_{i[k]})\bigg\} \le \xi^*$, as desired.~\Halmos
\endproof
}

\vspace{0.1in}


\rev{
\proof{Proof of Proposition \ref{prop:FS_bound_E}:}
Given a feasible solution $(\bar\bv,\bar\bu,\{\bar\bepsilon^k\}_{k\in\cK})$ for \eqref{eq:IO_mip}, let $\bar\cI=\{i\in \cI \,|\, \bar v_i=1\}$ and $\bar\cS=\{k\in \cK \,|\, \bar u_k=1\}$. Due to constraints \eqref{eq:IO_mip.2} and \eqref{eq:IO_mip.3}, for each $k\in \bar\cS$ we have $\hat\bx^k-\bar\bepsilon^k \in \{\bx\in \cX \,|\, {\ba^i}'\bx^k=b_i, \forall i\in \bar\cI\} = \underset{i\in \bar\cI}{\cap}\cX_i$ and $\|\bar\bepsilon^k\|_\ell\le \tau$. 
Now consider $\bar\bc\in\textup{cone}_+(\{\ba^i\}_{i:\bar{v}_i=1})$. This $\bar\bc$ vector leads to $\cX^*(\bar\bc)=\arg\min\;\textbf{FO}(\bar\bc)=\underset{i\in \bar\cI}{\cap}\cX_i$
. Recall that $\hat\bx^k-\bar\bepsilon^k \in \underset{i\in \bar\cI}{\cap}\cX_i$; this means that $\hat\bx^k-\bar\bepsilon^k \in\cX^*(\bar\bc)$ for all $k \in \bar\cS$. Moreover, because $\cX$ is an $n$-dimensional polytope and $|\bar\cI|=n$, $\cX^*(\bar\bc)$ is a singleton, i.e., $\cX^*(\bar\bc)=\{\hat\bx^k-\bar\bepsilon^k\},\forall k\in\bar\cS$ and $\hat\bx^k-\bar\bepsilon^k = \hat\bx^{k'}-\bar\bepsilon^{k'}$ for all $k, k' \in \cK$. To complete the proof, for all $k\in \bar\cS$ we have $\displaystyle \max_{\bx\in \cX^*(\bar\bc)}\{\|\hat\bx^k-\bx\|_\ell\}=\|\hat\bx^k-(\hat\bx^k-\bar\bepsilon^k)\|_\ell=\|\bar\bepsilon^k\|_\ell\le \tau$.\Halmos
\endproof
}
\
\proof{Proof of Proposition~\ref{prop:MQIO-FWStab}:}
(i) \rev{Let $\bar\cH=\{\ba^i: \bar v_i=1\}$ and $\tilde\cH=\{\ba^i: \tilde v_i=1\}$. Because $\bar\bv \ge \tilde \bv$, we have $\tilde \cH \subseteq \bar\cH$ and thus $\tilde \cC\subseteq \bar\cC$.} Given the dataset $\hat\cX$, we want to show that for any $\tilde\bc\in \tilde\cC$ there exists $\bar\bc\in\bar\cC$ such that $\underset{\bx \in \cX^*\!(\bar\bc)}{\max}  \{\|\hat\bx^k-\bx\|_\ell\} \le \underset{\bx \in \cX^*\!(\tilde\bc)}{\max}  \{\|\hat\bx^k-\bx\|_\ell\}$, $\forall k\in\cK$. Suppose to the contrary that there exists $\bc^*\in\tilde\cC$ such that 
$\underset{\bx \in \cX^*\!(\bc^*)}{\max}  \{\|\hat\bx^k-\bx\|_\ell\} < \underset{\bx \in \cX^*\!(\bc)}{\max} \{\|\hat\bx^k-\bx\|_\ell\}, \forall \bc\in\bar\cC$ for some $k\in\cK$. 
This means $\bc^*\in\tilde\cC\setminus\bar\cC$ and thus $\tilde\cC\setminus\bar\cC\ne\emptyset$, which contradicts $\tilde\cC \subseteq \bar\cC$. Thus, for any $\tilde\bc\in \tilde\cC$ there exists $\bar\bc\in\bar\cC$ such that $\underset{\bx \in \cX^*\!(\bar\bc)}{\max}  \{\|\hat\bx^k-\bx\|_\ell\} \le \underset{\bx \in \cX^*\!(\tilde\bc)}{\max}  \{\|\hat\bx^k-\bx\|_\ell\}$, $\forall k\in\cK$, as desired. 

\noindent(ii) Because $\bar\bv\ge \tilde\bv$, we have $\tilde\cH \subseteq \bar\cH$ and $\tilde\cC \subseteq \bar\cC$. %
Because $\tilde\cC \subseteq \bar\cC$, 
\begin{subequations}\nonumber
\begin{align}
\underset{\cX\subset\mathbb{R}^{K\times n}}{\min}\{d(\hat\cX,\cX) \,|\, \bar\cC\cap \cC(\cX)=\emptyset\}
    & = \underset{\cX\subset\mathbb{R}^{K\times n}}{\min}\{d(\hat\cX,\cX)\,|\, (\tilde\cC\cap \cC(\cX)) \cup((\bar\cC\setminus\tilde\cC) \cap \cC(\cX)) =\emptyset\} \\
& = \underset{\cX\subset\mathbb{R}^{K\times n}}{\min}\{d(\hat\cX,\cX) \,|\, \tilde\cC\cap \cC(\cX) =\emptyset  \text{ and } (\bar\cC\setminus\tilde\cC) \cap \cC(\cX) =\emptyset\} \\
&\ge \underset{\cX\subset\mathbb{R}^{K\times n}}{\min}\{d(\hat\cX,\cX)\,|\, \tilde\cC\cap \cC(\cX)=\emptyset\}, \text{\; as desired.~\Halmos} 
\end{align}
\end{subequations}
\endproof
\vspace{0.15in}
\proof{Proof of Lemma \ref{lemma 1}:}
	 Assume that there exists an optimal solution $(\bar\bv,\bar\bu,\{\bar{\bepsilon}^k\}_{k\in \cK})$ for $\mathbf{MQIO}(\cK,\tau,\theta)$ such that $\sum_{k=1}^{ K} \bar u_k > \lceil \theta K \rceil$. Without loss of generality let $\sum_{k=1}^{ K} \bar u_k = \lceil \theta K \rceil +r$ where $r>0$. Pick any $r$ observations for which $\bar u_k=1$ and change their values to $0$, and denote this new vector by $\tilde\bu$. Then we can construct a new solution $(\bar\bv,\tilde\bu,\{\bar{\bepsilon}^k\}_{k\in \cK})$, which is still feasible for $\mathbf{MQIO}(\cK,\tau,\theta)$ because $\|\bar\bepsilon^k\|_{\ell}\leq \tau+M_2(1-\bar u_k) \le \tau+M_2(1-\tilde u_k), \forall k \in \cK$ and $\sum_{k=1}^{ K} \tilde u_k = \lceil\theta K\rceil \ge \theta K$, and generates the same objective value, hence also an optimal solution.~\Halmos
\endproof
\vspace{0.15in}
\proof{Proof of Proposition \ref{prop:Dbar-IO}:} 
%
\rev{
Consider a matrix $\bar\bD$ constructed as in \eqref{eq:D_bar}. By definition, each row $k\in \cK$ of $\bar\bD$ is a binary vector ${\bv^k}'$ and there exists $\bar\bepsilon^k$ such that $(\bv^k,\bar\bepsilon^k)$ is feasible for \eqref{eq:IO_MIP_Reduced}. Assume that there exists an all-one submatrix in $\bar\bD$ whose rows and columns correspond to $\bar\cS$ and $\bar\cA$, respectively, such that $\bar\cS \subseteq \cK$, $|\bar \cS| \ge \theta K$, and $\bar\cA \subseteq \cI$. Define $\bar\bv$ such that $\bar v_i=1$ for all $i \in \bar\cA$ and $\bar v_i=0$ otherwise. Then we have $\bar \bv \le \bv^k, \forall k \in \bar\cS$, and thus $(\bar\bv,\bar\bepsilon^k)$ is also feasible for \eqref{eq:IO_MIP_Reduced} for all $k \in \bar\cS$. Define $\bar\bu$ such that $\bar u_k=1$ for all $k \in\bar\cS$ and $\bar u_k=0$ otherwise; note that $\sum_{k=1}^K \bar u_k \ge \theta K$. Let $\bepsilon^k=\bar\bepsilon^k$ for $k\in \bar\cS$ and $\bepsilon^k=\hat\bx^k-\bx^k$ for $k\in \cS\setminus\bar\cS$ where $\bx^k$ is arbitrarily chosen from $\{\bx\in \cX\ | \ {\ba^i}'\bx=b_i, \forall i\in \bar\cA\}$. We will first show that the solution $(\bar\bv,\bar\bu,\{\bepsilon^k\}_{k\in\cK})$ constructed as above is feasible for problem~\eqref{eq:IO_mip}. For all $i\in\bar\cA=\{i\,|\,\bar{v}_i = 1\}$, we have ${\ba^i}'(\hat\bx^k-\bepsilon^k)={\ba^i}'\bx^k=b_i$ for each $k\in\cK$. Also, for all $i\in\cI\setminus\bar\cA$ (i.e., $\{i\,|\,\bar{v}_i = 0\}$), we have ${\ba^i}'(\hat\bx^k-\bepsilon^k)={\ba^i}'\bx^k\ge b_i$ for each $k\in\cK$. Thus, $\bar\bv$ and $\{\bepsilon^k\}_{k\in\cK}$ satisfy \eqref{eq:IO_mip.2}.  Constraint \eqref{eq:IO_mip.3} is satisfied because $\|\bepsilon^k\|_\ell = \|\bar\bepsilon^k\|_\ell \le \tau$ for all $k\in\bar\cS = \{k\,|\,\bar{u}_k = 1\}$. Constraint \eqref{eq:IO_mip.4} is satisfied as we already know  $\sum_{k=1}^K \bar u_k \ge \theta K$. As a result, $(\bar\bv,\bar\bu,\{\bepsilon^k\}_{k\in\cK})$ is feasible for \eqref{eq:IO_mip}. Then, by Theorem \ref{prop:IO_mip}, there must also exist a solution $(\bar{\bc},\{{\bepsilon}^k\}_{k\in \cK},\bar\by,\bar{\cS})$ feasible for $\mathbf{QIO}(\cK,\tau,\theta)$ where $\bar\bc \in \textup{cone}(\{\ba^i\}_{i: \bar v_i = 1}) = \textup{cone}(\{\ba^i\}_{i \in \bar\cA})$, as desired.~\Halmos
\endproof
}
\vspace{0.15in}
\rev{
\proof{Proof of Proposition \ref{prop:exactness2}:} 
Let $P(r,k)$ denote the optimization problem solved in Line 5 of Algorithm~\ref{alg:Dbar-alg_Exact} with respect to $\hat\bx^k$, and $(\bv^{(r,k)},\bepsilon^{(r,k)})\in \underset{(\bv,\bepsilon)}{\text{argmax}} \ P(r,k)$. We first show that there exists  $\bar r\le \displaystyle\sum_{z= z^*}^n |\cV_z|$ such that $\bv^{(\bar r,1)}=\bv^*$ for problem $P(\bar r,1)$ (i.e., with respect to the first data point $\hat\bx^1$). For each iteration $r\ge2$, each cut only prevents the solution found in the previous iteration $r-1$, i.e., $\bv^{(r-1,1)}$, from being feasible for constraints \eqref{eq:IO_MIP_Reduced2}--\eqref{eq:IO_MIP_Reduced4}. This means that $P(r,1)$ can be solved $\displaystyle\sum_{z=1}^n |\cV_z|$ times with iteratively added cuts to collect all solutions satisfying \eqref{eq:IO_MIP_Reduced2}--\eqref{eq:IO_MIP_Reduced4} with respect to $\hat\bx^{1}$. Note that because $(\bv^*,\bepsilon^{k*})$ satisfies \eqref{eq:IO_mip.2} for all $k\in \cK$ and $u^*_1=1$, $(\bv^*,\bepsilon^{1*})$ satisfies constraints \eqref{eq:IO_MIP_Reduced2}--\eqref{eq:IO_MIP_Reduced4} with respect to $\hat\bx^{1}$. 
Thus, it is guaranteed that for some $\bar r\le \displaystyle\sum_{z=1}^n |\cV_z|$ we have $\bv^{(\bar r,1)}=\bv^*$. Furthermore, because the value function of $P(r,1)$ is non-increasing in $r$, for any $r>\displaystyle\sum_{z= z^*}^n |\cV_z|$ the optimal value of $P(r,1)$ is less than $z^*$ and thus $\bv^*$ is not achievable. Therefore, $\bar r$ must be no greater than $\displaystyle\sum_{z= z^*}^n |\cV_z|$.
}

\rev{
Next, we need to show that when $\bv^{(\bar r,1)}=\bv^*$ at iteration $\bar r$, the algorithm achieves $\bv^{Alg}= \bv^*$ and $z^{Alg}= z^*$. To do so, we first show that when $\bv^{(\bar r,1)}=\bv^*$, Lines 4--9 of the algorithm lead to $\bar\bD_{\bar r}$ with an all-one submatrix with $z^*$ columns and at least $\lceil \theta K\rceil$ rows. Let $\cI^*=\{i\in \cI  \,|\, v^*_i=1\}$ and $\cS^*=\{k\in \cK \,|\, u^*_k=1\}$. For some given $r$ and $k$, $P(r,k)$ has an optimal solution with $v_i^{(r,k)}=0$ if $w^k_i=0$; also, based on how the weight parameters  $\bw^k$ are updated (see Section \ref{sec:Algorithms}), if $v_i^{(r,k)}=0$ for $k=1,\ldots,k'-1$, we have $w^{k'}_i=0$. Generalizing this observation, if $\bv^{(\bar r,1)}=\bv^*$ 
then for all $k=2,\ldots,K$ 
we have $w^k_i=0$ for all 
$i\in\cI\setminus\cI^*$ and $w^k_i>0$ for all 
$i\in\cI^*$, and thus the maximum of the objective function of $P(\bar r,k)$ for each $k$, i.e., ${\bw^k}'\bv$, can be obtained when $\bv=\bv^*$. 
Because $(\bv^*,\bepsilon^{k*})$ satisfies all constraints of problem~\eqref{eq:IO_MIP_Reduced} as well as all cuts in $P(\bar r,k)$ generated up to iteration $\bar r$ for all $k\in \cS^*$, we have $(\bv^*,\bepsilon^{k*})\in \underset{\bv,\bepsilon}{\text{argmax}}\ P(\bar r,k)$. This means if ${\bv^*}'$ is assigned to the first row of $\bar\bD_{\bar{r}}$, it is guaranteed that ${\bv^*}'$ will also be assigned to rows $k\in \cS^*$ of $\bar\bD_{\bar{r}}$. Clearly, the submatrix of this $\bar\bD_{\bar{r}}$ with columns and rows corresponding to $i\in \cI^*$ and $k\in \cS^*$, respectively, is all-one and we have $|\cI^*|=z^*$ and $|\cS^*|\ge \theta K$. The clique problem in Line 10 of the algorithm (also see formulation~\eqref{eq:Biclique}) then takes this $\bar\bD_{\bar{r}}$ as input to return its optimal solution $\bar\bv^{\bar{r}}=\bv^*$ with ${\bar\bv^{r'}}\beee = z^*$. Since $z^*$ is the maximum achievable value (it is the optimal value of the exact MQIO problem), it satisfies the condition in Line 12 of Algorithm~\ref{alg:Dbar-alg_Exact} and thus we have $\bv^{Alg} = \bar\bv^{\bar{r}} = \bv^*$ and $z^{Alg}= z^*$.
}

\rev{
Finally, we show that the stopping criterion is valid, i.e., the while loop condition is not violated before $\bv^{Alg}=\bv^*$ is observed. Consider a certain iteration $\bar r\le \displaystyle\sum_{z= z^*}^n |\cV_z|$ and suppose we have $\bv^{Alg}\neq \bv^*$ and $z^{Alg}<z^*$ at the end of iteration $\bar r$ (i.e., Lines 13--14 of Algorithm~\ref{alg:Dbar-alg_Exact}). For the while loop condition to be violated in the next iteration, i.e., iteration $\bar r+1$, we must have $|\cG_{\bar r+1}|<z^{Alg}<z^*,$ which means ${\beee}'\bv^{(\bar r,1)}<z^*$. However, this cannot happen because we showed earlier in this proof that for all $r\le \displaystyle\sum_{z= z^*}^n |\cV_z|$ we have ${\beee}'\bv^{(r,1)}\ge z^*$.
~\Halmos
\endproof
}
\vspace{0.1in}
\proof{Proof of Proposition \ref{prop:IO1_minE}:} 
 Proposition \ref{prop:IO1} implies that given a solution $(\bar\bc,\{\bar\bepsilon^k\}_{k\in \cK},\bar\by,\bar\cS)$ feasible for \eqref{eq:IO1_minE}, $(\tilde\bc=\ba^{\tilde{i}},\{\bar\bepsilon^k\}_{k\in \cK},\tilde\by=\beee_{\tilde{i}},\bar\cS)$ is also feasible for any $\tilde{i} \in \bar{\cI}=\{ i \in \cI \  \arrowvert \ \bar y_i>0 \}$. For an arbitrary $\tilde i\in\bar\cI$, consider the following formulation defined for each $k \in \cK$:
\begin{align}\label{eq:Min_Fizbl_Dist}
 		\underset{\bepsilon^k}{\text{minimize}} \; \Big\{\, \lVert \bepsilon^k\rVert_{\ell} \, \Big{|} \,  \bA(\bxhat^k-\bepsilon^k) \geq \bbb, \; {\ba^{\tilde{i}}}'(\bxhat^k-\bepsilon^k)  = b_{\tilde{i}}  \Big\}.
 \end{align}
 Let $\tilde\bepsilon^k$ be the optimal solution to \eqref{eq:Min_Fizbl_Dist} for each $k \in \cK$ (here we suppress the index $\tilde{i}$ for brevity). Note that a solution constructed as $(\tilde\bc=\ba^{\tilde i},\{\tilde\bepsilon^k\}_{k\in \cK},\tilde\by=\beee_{\tilde i },\bar\cS)$ is also feasible for \eqref{eq:IO1_minE} because $\lVert \tilde\bepsilon^k\rVert_{\ell} \le \lVert \bar\bepsilon^k\rVert_{\ell} \le \tau$ for each $k \in \bar\cS$. 
 The objective of \eqref{eq:IO1_minE}, i.e., minimizing $\tau$, is equivalent to minimizing $\underset{k \in \bar\cS}{\max}\{ \lVert \bepsilon^k\rVert_{\ell} \}$. Because $|| \tilde\bepsilon^k||_{\ell} \le || \bar\bepsilon^k||_{\ell} \le \tau$ for each $k \in \bar\cS$, we have 
 $\underset{k \in \bar\cS}{\max}\{ \lVert \tilde\bepsilon^k\rVert_{\ell} \}  \leq  \underset{k \in \bar\cS}{\max}\{ \lVert \bar\bepsilon^k\rVert_{\ell} \}
$. That is, for any given feasible solution $(\bar\bc,\{\bar\bepsilon^k\}_{k\in \cK},\bar\by,\bar\cS)$ for \eqref{eq:IO1_minE} where $\bar\bc$ is not identical to $\ba^i$ for any $i\in\cI$, we can always construct another feasible solution $(\tilde\bc=\ba^{\tilde i},\{\tilde\bepsilon^k\}_{k\in \cK},\tilde\by=\beee_{\tilde i },\bar\cS)$ 
without increasing the objective value of \eqref{eq:IO1_minE}, which completes the proof.~\Halmos
\endproof
\vspace{0.1in}
\rev{
\proof{Proof of Lemma \ref{lemma:Sufficient_condition}:}
Let $\hat\cI = \{i\in\cI \,|\, \ba^{i}\in\hat\cC\}$. Given shifted data $\tilde\cX$, assume $\bc=\ba^i$ is not inverse-feasible for \textbf{QIO}($\tilde\cX,\tau,\theta$) for any $i \in \hat\cI$. We want to show that there is no $\bc\in\hat\cC$ that is inverse-feasible for \textbf{QIO}($\tilde\cX,\tau,\theta$). 
Suppose to the contrary that there is $\bar\bc\in\hat\cC$ feasible for \textbf{QIO}($\tilde\cX,\tau,\theta$) and let $\bar\cI$ be such that $\bar\bc\in \textup{cone}(\{\ba^{i}\}_{i\in\bar\cI})$. Then, from Proposition~\ref{prop:IO1}, there must be some $\bar i\in\bar\cI$ such that $\bc=\ba^{\bar i}$ is also feasible for \textbf{QIO}($\tilde\cX,\tau,\theta$). Furthermore, since $\bar\bc\in\hat\cC$, this means that $\bc=\ba^{\bar i}$ must be also feasible for \textbf{QIO}($\hat\cX,\tau,\theta$) (i.e., it must be that $\bar i \in \hat\cI$), which is a contradiction.\Halmos
\endproof
}
\vspace{0.1in}
\rev{
\proof{Proof of Proposition \ref{lemma:FS_bound}:}
Since $\bar\bc$ is a strict conic combination of the selected $\ba^i$ vectors, we have $\cX^*(\bar\bc)=\underset{i\in \bar\cI}{\cap} \cX_i$. For each $k\in \bar\cS$, let $\tilde\bx^k\in \underset{\bx\in \underset{i\in \bar\cI}{\cap}\cX_i}{\text{argmin}} \{\|\hat\bx^k-\bx\|_\ell\}$. Then, for all $k\in\bar\cS$ and $i\in \bar\cI$, we have  $\displaystyle\max_{\bx\in \cX_i}\{\|\hat\bx^k-\bx\|_\ell\}\le  \|\tilde\bx^k-\hat\bx^k\|_\ell+\displaystyle\max_{\bx\in \cX_i}\{\|\tilde\bx^k-\bx\|_\ell\}\le \|\tilde\bx^k-\hat\bx^k\|_\ell+\rho_i$. Hence, for each $k\in \bar\cS$,
\begin{subequations}\nonumber
\begin{align}
\displaystyle \max_{\bx\in \cX^*(\bar\bc)}\{\|\hat\bx^k-\bx\|_\ell\}
 &=\displaystyle \max_{\bx\in \underset{i\in \bar\cI}{\cap} \cX_i} \; \{\|\hat\bx^k-\bx\|_\ell\} \\
 &\le  \displaystyle\min_{i\in\bar\cI}\;  
 \big\{\displaystyle\max_{\bx\in \cX_i}\{\|\hat\bx^k-\bx\|_\ell\}\big\}\\
  &\le   \displaystyle\min_{i\in\bar\cI}\;  
 \big\{\|\tilde\bx^k-\hat\bx^k\|_\ell+\rho_i \big\}\\
  &= \displaystyle\min_{i\in\bar\cI}\;  
 \big\{\underset{\bx\in \underset{i\in \bar\cI}{\cap}\cX_i}{\text{min}} \{\|\hat\bx^k-\bx\|_\ell\}+\rho_i \big\},\\
 &= \underset{\bx\in \underset{i\in \bar\cI}{\cap}\cX_i}{\text{min}} \{\|\hat\bx^k-\bx\|_\ell\}+ \displaystyle\min_{i\in\bar\cI}\;  
 \{\rho_i \},
\end{align}
\end{subequations} 
where the first inequality holds because we have $\underset{i\in \bar\cI}{\cap}\cX_i\subseteq \cX_i$ for all $i\in \bar\cI$, which leads to $\displaystyle \max_{\bx\in \underset{i\in \bar\cI}{\cap} \cX_i}\{\|\hat\bx^k-\bx\|_\ell\}\le \max_{\bx\in \cX_i}\{\|\hat\bx^k-\bx\|_\ell\}, \forall i\in \bar\cI$.\Halmos
\endproof
}
\vspace{0.05in}
\rev{
\proof{Proof of Lemma \ref{lemma:Alpha_vs_ReducedMIP}:}
Because $\bar\alpha^k_i=0$ for all $i\in \bar\cI$ and $(\bar\balpha^k,\bar\bepsilon^k)$ is feasible for \eqref{eq:IO_Alpha}, we have ${\ba^i}'(\bxhat^k-\bar\bepsilon^k)=b_i$ for all $i\in \bar\cI$, ${\ba^i}'(\bxhat^k-\bar\bepsilon^k)\ge b_i$ for all $i\in \cI\setminus\bar\cI$, and $\rVert\bar\bepsilon^k\lVert_{\ell}\le \tau$. Let $\bar v_i^k=1-\bar\alpha_i^k$ for all $i\in \bar\cI$ and $\bar v_i^k = 0$ otherwise. Clearly, $(\bar\bv^k,\bar\bepsilon^k)$ satisfies all constraints of \eqref{eq:IO_MIP_Reduced}.\Halmos
\endproof 
}
\vspace{0.05in}
\rev{
\proof{Proof of Proposition \ref{prop:Dtilde-IO}:}
Assume that $\tilde\bD$ satisfies \eqref{D_tilde} and has an all-one submatrix whose rows and columns correspond to $\bar\cS$ and $\bar\cA$, respectively, where $\bar\cS\in\cK, |\bar\cS|\ge \theta K$, and 
$\bar\cA\in\cI$. By definition, for each $k\in \cK$ there exists a feasible solution $(\bar\balpha^k,\bar\bepsilon^k)$ for \eqref{eq:IO_Alpha} with respect to $\hat\bx^k$ where $\bar\alpha_i^k=0$ for all $i\in\bar\cA$ and $k\in \bar\cS$, and $\bar\alpha_i^k>0$ otherwise. This means that, by Lemma \ref{lemma:Alpha_vs_ReducedMIP}, for each $k\in\cK$ there exists a corresponding feasible solution $(\bar\bv^k,\bar\bepsilon^k)$ for \eqref{eq:IO_MIP_Reduced} with respect to $\hat\bx^k$ where $\bar v_i^k=1$ for all $i\in\bar\cA$ and $k\in \bar\cS$, and $\bar v_i^k=0$ otherwise. Construct a matrix $\bar\bD\in\{0,1\}^{K\times m}$ such that $\bar\bD_{ki}=1$ if $\bar v_i^k=1$ and $\bar\bD_{ki}=0$ otherwise. Note that $\bar\bD$ satisfies \eqref{eq:D_bar} and has an all-one submatrix whose rows and columns correspond to $\bar\cS$ and $\bar\cA$, respectively, where $\bar\cS\subseteq\cK$, $|\bar\cS|\ge \theta K$, and $\bar\cA\subseteq\cI$. Thus, by Proposition \ref{prop:Dbar-IO}, there exists a solution $(\bar{\bc},\{\bar{\bepsilon}^k\}_{k\in \cK},\bar\by,\bar{\cS})$ feasible for $\mathbf{QIO}(\cK,\tau,\theta)$ where $\bar\bc \in \textup{cone}(\{\ba^i\}_{i \in \bar\cA})$
. 
~\Halmos
\endproof
}

\section{Supplemental Materials for the Numerical Results}
\subsection{The Diet Problem: Data}\label{appx:Diet_Data}

\begin{threeparttable}
\caption{Food items and nutrient data per serving.}
\label{tab:DietData}%
{\fontsize{9.5}{10} \selectfont
\begin{tabular}{lrrrrrrrrrrr}
     			\toprule
     			& \multicolumn{9}{c}{Food Type}                               &   {Lower}    &  {Upper} \\
     			\cmidrule{2-10}    \multicolumn{1}{c}{ } & 1  & 2 & 3 & 4  & 5 & 6  & 7  & 8 & 9 &  {limit}  & {limit} \\
     			\midrule
     			Energy (KCAL) & 91.53 & 68.94 & 23.51 & 65.49 & 110.88 & 83.28 & 80.50 & 63.20 & 52.16 & 1800.00 & 2500.00 \\
     			Total\_Fat (g) & 4.95  & 0.71  & 1.80  & 3.48  & 6.84  & 4.41  & 5.80  & 0.94  & 0.18  & 44.00 & 78.00 \\
     			Carbohydrate (g) & 6.89  & 12.16 & 0.25  & 0.00  & 5.44  & 4.68  & 0.56  & 11.42 & 13.59 & 220.00 & 330.00 \\
     			Protein (g) & 4.90  & 3.68  & 1.59  & 7.99  & 6.80  & 5.93  & 6.27  & 2.40  & 0.41  & 56.00 & NA \\
     			Fiber (g) & 0.00  & 0.06  & 0.00  & 0.00  & 0.28  & 0.29  & 0.00  & 1.19  & 1.81  & 20.00 & 30.00 \\
     			Vitamin C (mg) & 0.01  & 1.76  & 0.00  & 0.00  & 0.17  & 0.16  & 0.00  & 0.02  & 11.19 & 90.00 & 2000.00 \\
     			Vitamin B6 (mg) & 0.06  & 0.03  & 0.01  & 0.09  & 0.11  & 0.06  & 0.06  & 0.03  & 0.10  & 1.30  & 100.00 \\
     			Vitamin B12 (mcg) & 0.67  & 0.39  & 0.09  & 0.65  & 0.11  & 0.63  & 0.56  & 0.00  & 0.00  & 2.40  & NA \\
     			Calcium (mg) & 172.09 & 125.72 & 46.24 & 2.21  & 5.90  & 15.03 & 29.00 & 27.21 & 6.14  & 1000.00 & 2500.00 \\
     			Iron (mg) & 0.05  & 0.08  & 0.04  & 0.75  & 0.35  & 0.35  & 0.73  & 0.79  & 0.13  & 8.00  & 45.00 \\
     			Copper (mg) & 0.02  & 0.04  & 0.01  & 0.03  & 0.03  & 0.03  & 0.04  & 0.05  & 0.06  & 0.90  & 10.00 \\
     			Sodium (mg) & 61.02 & 48.24 & 65.08 & 72.32 & 211.05 & 128.27 & 223.50 & 125.62 & 1.42  & 1500.00 & 2300.00 \\
     			Vitamin A (mcg) & 42.89 & 22.24 & 12.78 & 0.00  & 1.33  & 9.53  & 81.00 & 0.01  & 13.04 & 900.00 & 3000.00 \\
     			\multirow{2}[1]{*}{Max serving } & \multirow{2}[1]{*}{4} & \multirow{2}[1]{*}{8} & \multirow{2}[1]{*}{4} & \multirow{2}[1]{*}{5} & \multirow{2}[1]{*}{5} & \multirow{2}[1]{*}{5} & \multirow{2}[1]{*}{4} & \multirow{2}[1]{*}{4} & \multirow{2}[1]{*}{8} & \multicolumn{2}{c}{\multirow{2}[1]{*}{}} \\
     			&       &       &       &       &       &       &       &       &       & \multicolumn{2}{c}{} \\
     			\bottomrule
\end{tabular}}
\end{threeparttable}

\vspace{0.2in}

\subsection{The Transshipment Problem: Problem Description and Additional Numerical Results} \label{appx:Transportaion}
Figure~\ref{fig:Graph} and the following formulation show the transshipment problem we consider.
\begin{subequations} \label{eq:Transship}
	\begin{align}
	 \underset{\bx^{(p)},\bx^{(t)}}{\text{minimize}} 
	&\quad \sum_{i\in \cN_s}c^{(p)}_i x_i^{(p)} + \sum_{i\in\cN}\sum_{j \in \cN}c^{(t)}_{ij} x_{ij}^{(t)}& \\ \label{eq:Transship1} 
	\text{subject to} 
	&\quad  \sum_{j\in \cN\setminus\cN_s} x_{ij}^{(t)}=x_i^{(p)}, \quad \forall i \in \cN_s, \\ \label{eq:Transship2}
	&\quad  \sum_{i\in \cN\setminus\cN_d} x_{ij}^{(t)}=d_j, \quad \forall j \in \cN_d, \\ \label{eq:Transship3}
	&\quad  \sum_{j\in \cN} x_{ij}^{(t)} - \sum_{j\in \cN} x_{ji}^{(t)}=0, \quad \forall i \in \cN_t, \\ \label{eq:Transship4}
	& \quad 0\le x_i^{(p)} \le \pi^{(p)}_i, 0\le x_{ij}^{(t)} \le \pi^{(t)}_{ij}, \quad \forall i,j \in \cN, & 
	\end{align}
\end{subequations}
where $\cN_s$, $\cN_d$, and $\cN_t$ denote the set of supply, demand, and transshipment nodes, respectively, and $\cN$ denotes the set of all nodes. Variables $\bx^{(p)}$ and $\bx^{(t)}$ denote the production level at each supply node and transshipment flow on each arc, respectively. Parameters $\boldsymbol{\pi}^{(t)}$ and $\boldsymbol{\pi}^{(p)}$ represent the arc capacity and production capacity, respectively. 

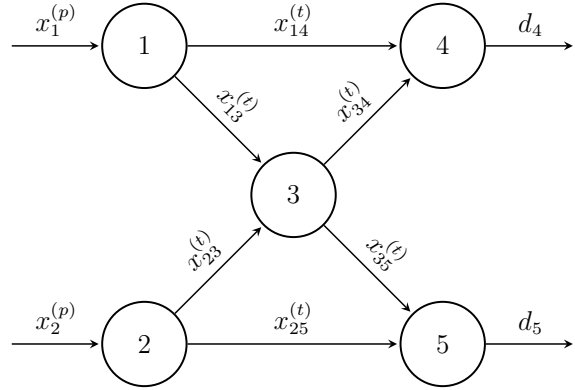
\begin{figure}[h]
    \centering
    \resizebox{0.6\textwidth}{!}{%
    \begin{tikzpicture}[
            > = stealth, 
            shorten > = 1pt, 
            auto,
            node distance = 3cm, 
            semithick 
        ]
        \centering

        \tikzstyle{every state}=[
            draw = black,
            thick,
            fill = white,
            minimum size = 4mm
        ]
        \node[state,minimum size=1.2cm]  (1) {$1$};
        \node[state,minimum size=1.2cm,color=white,node distance = 2.5cm]  (01) [left of=1] {};
        \node[state,minimum size=1.2cm]  (3) [below right of=1] {$3$};
        \node[state,minimum size=1.2cm]  (2) [below left of=3] {$2$};
        \node[state,minimum size=1.2cm,color=white,node distance = 2.5cm]  (02) [left of=2] {};
        \node[state,minimum size=1.2cm]  (4) [above right of=3] {$4$};
        \node[state,minimum size=1.2cm,color=white,node distance = 2.5cm]  (40) [right of=4] {};
        \node[state,minimum size=1.2cm]  (5) [below right of=3] {$5$};
        \node[state,minimum size=1.2cm,color=white,node distance = 2.5cm]  (50) [right of=5] {};
         \path[->,rotate=45,label=above:rotate] (01) edge node {$x_1^{(p)}$} (1);
         \path[->] (02) edge node {$x_2^{(p)}$} (2);
         \path[->] (1) edge node {$x_{14}^{(t)}$} (4);
         \path[->] (1) edge node [above,rotate=-45]{$x_{13}^{(t)}$} (3);
         \path[->] (2) edge node [above,rotate=45] {$x_{23}^{(t)}$} (3);
         \path[->] (2) edge node {$x_{25}^{(t)}$} (5);
         \path[->] (3) edge node [above,rotate=45] {$x_{34}^{(t)}$} (4);
         \path[->] (3) edge node [above,rotate=-45] {$x_{35}^{(t)}$} (5);
         \path[->] (4) edge node {$d_4$} (40);
         \path[->] (5) edge node {$d_5$} (50);
    \end{tikzpicture}
    }
    \caption{Network for the transshipment example \citep{dong2018generalized}. Nodes 1 and 2 are supply nodes and 4 and 5 are demand nodes. The capacity of each arc is 1.3 and production capacities are 3 and 1.5 for nodes 1 and 2, respectively. The true production costs are $c_1^{(p)}=0.2393$ and $c_2^{(p)}=0.1496$ and true transshipment costs are $c_{13}^{(t)}=0.0935$, $c_{14}^{(t)}=0.1232$, $c_{23}^{(t)}=0.1141$, $c_{25}^{(t)}=0.0320$, $c_{34}^{(t)}=0.1615$, and $c_{35}^{(t)}=0.0867$.}
    \label{fig:Graph}
\end{figure}

In addition to the numerical results presented in Section~\ref{sec:results_transshipment}, we also compare the objective function values achieved by our forward solution $\bar\bx^t$ (recall that $\bx = [{\bx^{(p)}}; {\bx^{(t)}}]$ and $\bc = [{\bc^{(p)}}; {\bc^{(t)}}]$) and that from the given data. We adopt the metrics from the previous study on online inverse LP \citep{barmann2017emulating}; this study assumes a single noiseless (i.e., optimal) decision at each time point and shows that both $1/T\sum_{t=1}^{T}{\bc^{true}}'(\bar\bx^t-\hat\bx^t)$ and $ 1/T\sum_{t=1}^{T}{\bc^t}'(\bar\bx^t-\hat\bx^t)$ converge to zero as $T\rightarrow \infty$. In our case, because we consider batches of data at each time we use the following modified metrics: $ 1/T \big[ \sum_{t=1}^{T} \sum_{k=1}^{K_t} {\bc^{true}}'(\bar\bx^t-\hat\bx^k)/K_t \big]$ and $ 1/T \big[\sum_{t=1}^{T}  \sum_{k=1}^{K_t} {\bc^t}'(\bar\bx^t-\hat\bx^k) /K_t \big]$. Note that $\bc^t$ is a cost vector randomly selected from the inverse set at time $t$. Figure \ref{fig:ObjvalRound1} shows these metrics achieved by our algorithm with noisy datasets at each iteration. Both metrics averaged over the data points also converge; however, neither of them reaches 0 because of the noise in the data.

\begin{figure}[h!]
	\centering
	\begin{subfigure}[t]{0.5\textwidth}
		\centering
		\includegraphics[width=3in]{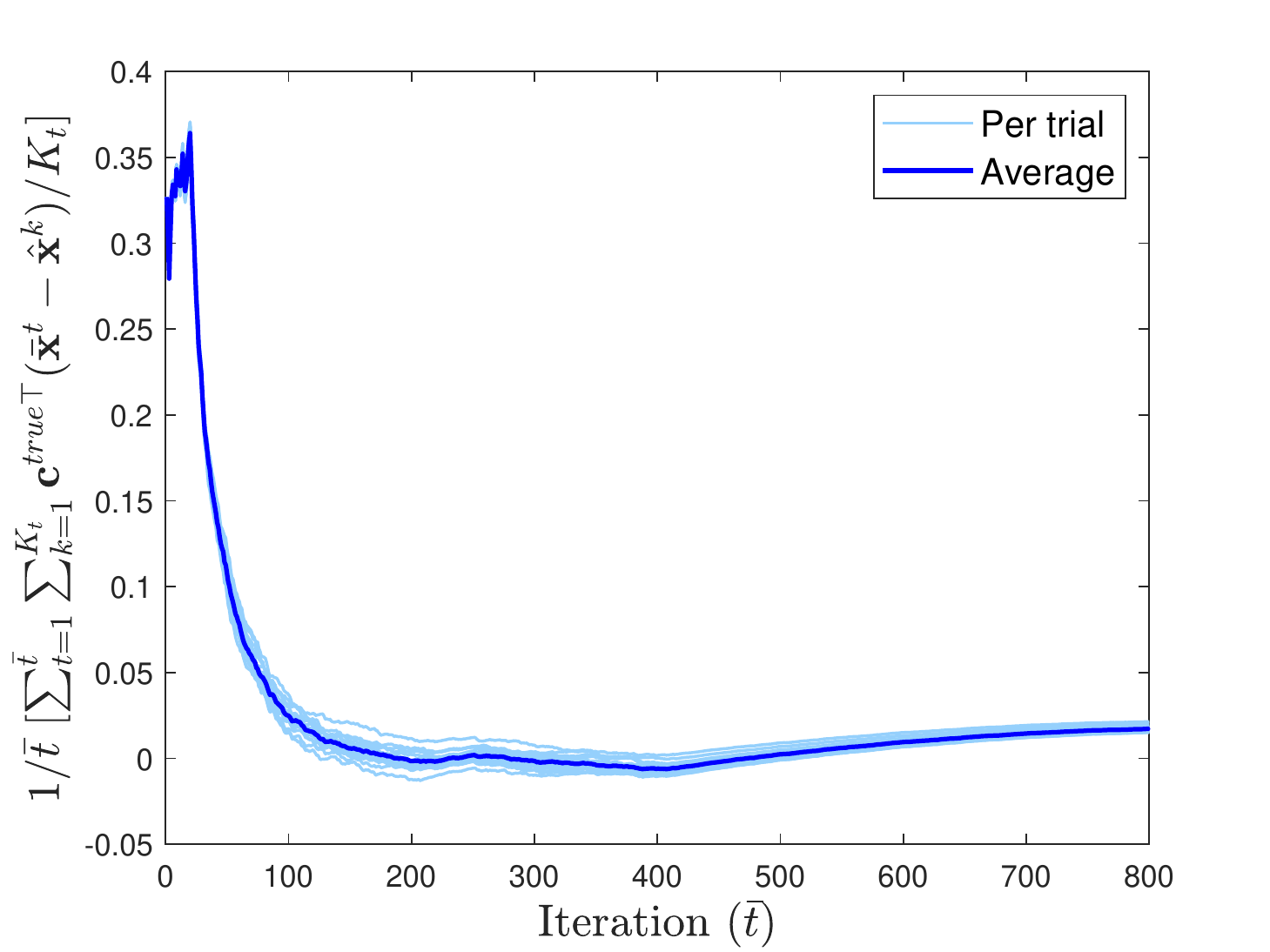}
		\caption{Average cumulative error under $\bc^{true}$.}
		\label{fig:ObjvalTrueRound1}
	\end{subfigure}%
	\begin{subfigure}[t]{0.5\textwidth}
		\centering
		\includegraphics[width=3in]{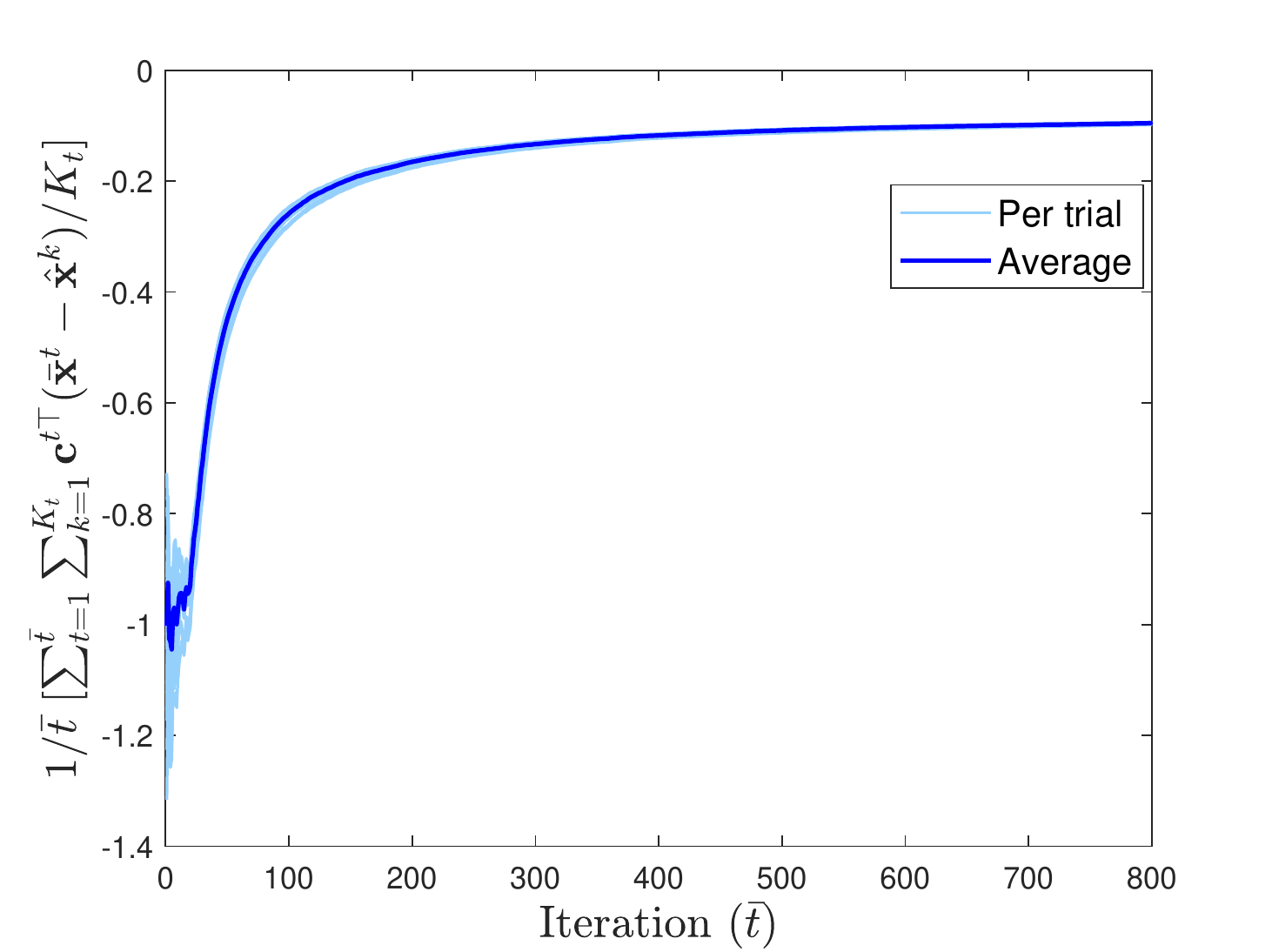}
		\caption{Average cumulative error under $\bc^t\in \cC_t$.}
		\label{fig:ObjvalourCRound1}
	\end{subfigure}%
	\caption{Convergence of average cumulative objective function errors.}
	\label{fig:ObjvalRound1}
\end{figure}

\end{APPENDICES}





\bibliographystyle{informs2014} 
\bibliography{main} 

\end{document}